 \documentclass[12pt,reqno]{amsart}

\usepackage{hyperref}
\usepackage{amsmath}
 \usepackage[dvips]{graphicx}
 \usepackage{amsfonts}
  \usepackage{amssymb}
      \usepackage{amsthm}
 \usepackage{enumerate}
  \usepackage{latexsym}
 \usepackage{color}

 \renewcommand{\proof}{\noindent{\bf Proof. \quad}}
 \newcommand{\eproof}{\hfill \mbox{${\square}$} \vspace{3mm}}

\newcommand{\R}{\mathbb{R}}

\newcommand{\C}{\mathcal{C}}
\newcommand{\com}{\mathbb{C}} 

\def\dif{\mathrm{Diff}}

\newtheorem {teo} {\bf Theorem\,} [section]
\newtheorem {prop} [teo] {\bf Proposition}
\newtheorem {cor} [teo]{\bf Corollary}
\newtheorem {lema} [teo] {\bf Lemma}
\newtheorem {defn} [teo] {\bf Definition}
\newtheorem{rem}[teo]{\bf Remark}

\makeatletter 
\def\namedlabel#1#2{\begingroup
   \def\@currentlabel{#2}%
   \label{#1}\endgroup
}

\begin{document}

\title[{continuity of attractors in Lipschitz}]{Continuity of attractors for a family of $\mathcal{C}^1$ perturbations of the square}

\author[Pricila S. Barbosa,  Ant\^onio L. Pereira and Marcone C. Pereira]{Pricila S. Barbosa$^{\star}$, Ant\^onio L. Pereira$^{\diamond}$ and Marcone C. Pereira$^{\dagger}$}

\thanks{$^\star$Partially supported by CNPq-Brazil grant 142278/2010-6 and 233524/2014-2.}
\thanks{$^\diamond$Partially supported by CNPq-Brazil grant 141882/2003-4.}
\thanks{$^\dagger$Partially supported by CNPq-Brazil grant 302960/2014-7 and 471210/2013-7, FAPESP 2013/22275-1.}

\address{Pricila S. Barbosa \hfill\break
 Departamento de Matem\'atica Aplicada\\
Universidad Complutense de Madrid - Madrid - Spain}
\email{pricila@ime.usp.br}

\address{Ant\^onio L. Pereira \hfill\break
Instituto de Matem\'atica e Estat\'istica \\
Universidade de S\~ao Paulo - S\~ao Paulo - Brazil}
\email{alpereir@ime.usp.br}

\address{Marcone C. Pereira \hfill\break
Instituto de Matem\'atica e Estat\'istica \\
Universidade de S\~ao Paulo - S\~ao Paulo - Brazil}
\email{marcone@ime.usp.br}

\date{}


\begin{abstract}

We consider here the family of semilinear parabolic problems 

\begin{equation*}
\begin{array}{rcl}
\left\{
\begin{array}{rcl}
u_t(x,t)&=&\Delta u(x,t) -au(x,t) + f(u(x,t)) ,\,\,\ x \in \Omega_\epsilon
 \,\,\,\mbox{and}\,\,\,\,\,\,t>0\,,
\\
\displaystyle\frac{\partial u}{\partial N}(x,t)&=&g(u(x,t)), \,\, x \in \partial\Omega_\epsilon \,\,\,\mbox{and}\,\,\,\,\,\,t>0\,,
\end{array}
\right.
\end{array}
\end{equation*}
where $ {\Omega} $ is the unit square, $\Omega_{\epsilon}=h_{\epsilon}(\Omega)$ and $h_{\epsilon}$ is a family of diffeomorphisms converging to the identity in the $C^1$-norm. We show that the problem is well posed for $\epsilon>0$ sufficiently small in a suitable phase space,  the associated semigroup has a global attractor $\mathcal{A}_{\epsilon}$ and the family $\{\mathcal{A}_{\epsilon}\}$  is continuous at $\epsilon = 0$.

\end{abstract}

\maketitle

 \allowdisplaybreaks

\section{Introduction} \label{intro}

\begin{figure}[!h] \label{figura}
\centering
 \includegraphics[width=.70\columnwidth]{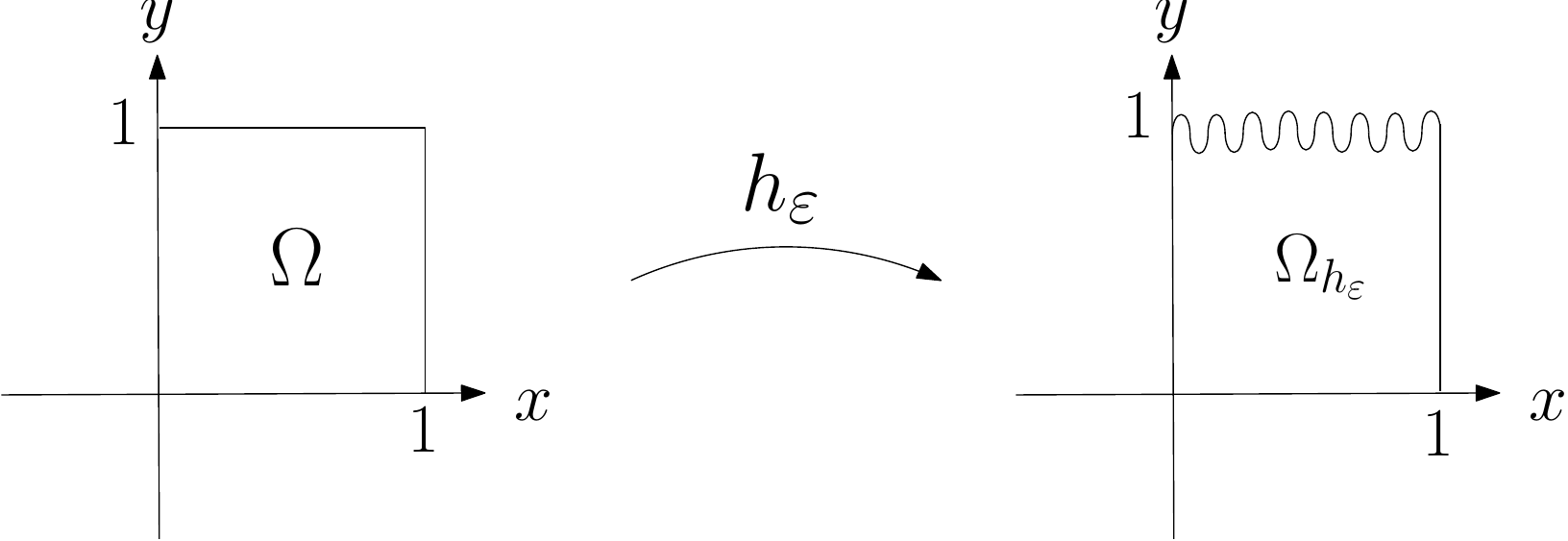}
 \caption{Map $h_\epsilon$.}
\end{figure}

  Let $ { \Omega} \subset \R^2$ be the unit square, $a$ a positive   {number}  $f, g: \R \to \R$ real functions, and consider  the family of semilinear parabolic problems with nonlinear Neumann boundary conditions:
 
  \begin{equation} \label{nonlinBVP} \tag{$P_{\epsilon}$}
\begin{array}{rcl}
\left\{
\begin{array}{rcl}
u_t(x,t)&=&\Delta u(x,t) -au(x,t) + f(u(x,t)) ,\,\,\ x \in \Omega_\epsilon
 \,\,\,\mbox{and}\,\,\,\,\,\,t>0\,,
\\
\displaystyle\frac{\partial u}{\partial N}(x,t)&=&g(u(x,t)), \,\, x \in \partial\Omega_\epsilon \,\,\,\mbox{and}\,\,\,\,\,\,t>0\,,
\end{array}
\right.
\end{array}
\end{equation}
where  $\Omega_{\epsilon}=h_{\epsilon}(\Omega)$ and $h_{\epsilon}$ is the family of diffeomorphisms, given by 

\begin{equation} \label{per}
h_{\epsilon} (x_1,x_2) =  (\,x_1\,,\,x_2 + x_2\,\epsilon\, sen(x_1/\epsilon^\alpha)\,)
\end{equation}
 with $0 <\alpha <1$  and  $\epsilon >0$ is sufficiently small
 (see figure  (\ref{figura})).

  Our aim here is to prove well-posedness, 
  establish the existence of a global attractor $\mathcal{A}_{\epsilon}$
  for sufficiently small $\epsilon$ and prove the continuity of the family of attractors  at  $\epsilon=0$, under appropriate conditions on the nonlinearities. 

 It well known that, under some  smoothness hypotheses on  $\Omega \subset \R^n$, $f$ and $g$, the problem is well posed in appropriate phase spaces. 
 The existence of a global compact attractor has been proved in 
 \cite{COPR} and \cite{OP}, under some additional growth and dissipative conditions on the nonlinearities $f$ and $g$.  In \cite{PP} the authors prove the continuity of the attractors with respect to 
 $C^2$ perturbations of a smooth domain of $\R^n$.
  
  These results do not extend immediately to the case considered here, due to the lack of smoothness of the domains considered and the fact that the perturbations do not
 converge to the inclusion in the $C^2$ norm.
 
 \par The problem of existence and continuity of global attractors for semilinear parabolic problems, with respect to change of domains has also been considered  in  \cite{AC1}, for the problem with homogeneous boundary condition 
\begin{equation*}
\begin{array}{rlr}
\left\{
\begin{array}{rcl}
u_t = \Delta u + f(x,u) \,\,\,\mbox{in}\,\,\, \,\,\,\Omega_\epsilon \\
\displaystyle\frac{\partial u}{\partial N}= 0 \,\,\,\mbox{on}\,\,\,\,\,\,\partial\Omega_\epsilon\,,
\end{array}
\right.
\end{array}
\end{equation*}
where  $\Omega_\epsilon$, $0 \leq \epsilon \leq \epsilon_0$ are bounded domains
with  Lipschitz  boundary in $\R^N$, $N \geq 2$. They prove that, if the perturbations are such that the convergence of the eigenvalues and eigenfunctions of the linear part of the problem can be shown, than the upper semicontinuity of attractors follow. With the additional assumption that the equilibria are all hyperbolic, the lower semicontinuity is also obtained.

\par  The behavior of the equilibria of (\ref{nonlinBVP}) was studied in
  \cite{AB2} and \cite{AB3}. In these papers, the authors consider
 a family of smooth domains  $\Omega_\epsilon \subset \R^N$, $N \geq 2$ and  $0 \leq \epsilon \leq \epsilon_0$ whose boundary oscillates rapidly when the parameter
 $\epsilon \to 0$ and prove that the solutions, as well as the spectra of the linearised problem around them, converge to the solution of a ``limit problem''.

 In this work, we follow the general   approach of \cite{PP}, which consists basically in ``pull-backing'' the perturbed problems to the fixed domain $ {\Omega}$
 and then considering the  family of abstract semilinear problems thus generated.
  

 We observe  that  the results obtained can be easily  extended to 
 \emph{convex} domains, and more general families of $C^1$ perturbations,
 but we have chosen to consider a specific setting, for the sake of clarity. The extension to more general Lipschitz domains is problematic for the lack of appropriate regularity results.

 The paper is organized as follows: in section \ref{prelim} we collect some results needed later.  In section \ref{reduction} we show how the problem can be reduced to a family of problems in the initial domain and, in section
 \ref{sectorialoper}, we show that the perturbed linear operators are sectorial operators in suitable spaces. In section \ref{abstract} we show that the problem can be reformulated as an abstract problem in a scale of Banach spaces  which are shown to be locally
 well-posed in section \ref{wellposed}, under suitable growth assumptions onf
 $f$ and $g$.  In  section \ref{lyapunov}, assuming a dissipative condition for the problem,   using the properties of
 a Lyapunov functional, we prove that the solutions are globally defined.
 In section \ref{existglob} we prove the existence of global attractors.
 Finally, in section \ref{contattractors}, we show first that these attractors
 behave upper semicontinuosly and, with some additional properties on the nonlinearities and on the set of equilibria, we  show that they are also lower semicontinuos at $\epsilon=0$ in the $H^{s}$-norm for $s < 1$.

\section*{Preliminaries} \label{prelim}

We collect here some definitions and results that will be used in the sequel.

 \subsection{Boundary perturbation}

 One of the difficulties encountered in problems of perturbation of the domain
 is   that the function spaces change with the change of the region. One way to  overcome this difficulty is to  effect a ``change of variables'' in order to bring the problem back to a fixed region. This approach was developed by D. Henry in \cite{He1} and is the one we adopt here.

If $\Omega \subset \R^n$ is a bounded region, we denote by $Diff^m(\Omega),
 m \geq 0$, 
the set of $\mathcal{C}^m $ embeddings (=diffeomorphisms from $\Omega$
 to  its  image).
  
If $ h \in Diff^m(\Omega)$,   we may define the  ``pull-back'' map
     $h^*:  C^m(h(\Omega)) \to C^m(\Omega)$ by
$$h^*u(x)= (u\circ h)(x)=u(h(x)), x \in \Omega.$$  

\begin{prop}
 The map
$$
\begin{array}{lclc}
h^*:& C^m(h(\Omega))& \to &C^m(\Omega) \\
&u&\longmapsto & u\circ h
\end{array}
$$
is an isomorphism, with inverse  $h^{*-1}=(h^{-1})^*$, for  $m \geq 0$.
\end{prop}

\noindent \proof See \cite{He1}. \eproof

It is important to observe that the above result is also true in other function
 (e.g. Sobolev) spaces. (see Lemma \ref{isosobolev}).

If $L$ is a (formal)  constant matrix coefficient differential operator
$$Lu(x)= \left(\,u(x), \frac{\partial u}{\partial x_1}(x), \cdots, \frac{\partial u}{\partial x_n}(x), \frac{\partial^2 u}{\partial x_1^2}(x),\frac{\partial^2 u}{\partial x_1 \partial x_2}(x), \cdots\,\right),\,\,x \in \R^n,$$ with values in $\R^k$
and  $f(x,\lambda)$ is a function   defined for $(x, \lambda)$
in some open set  $O \subset \R^n \times \R^k$, we may define the nonlinear differential operator   $F_\Omega$, by 
  $$F_\Omega(u)(x) = f(x, Lu(x)), x \in \Omega, $$ for  sufficiently smooth functions $u$ defined in   $\Omega$ such that  $(x, Lu(x)) \in O$, for any 
 $x \in \overline{\Omega}$.
 
\par If  $h : \Omega \to \R^n \in Diff^m(\Omega)$  we may consider the differential operator  $F_{h(\Omega)}$  with  $D(F_{h(\Omega)}) \subset   C^m(h(\Omega))$, open.
  The map  $$F_{h(\Omega)}: D(F_{h(\Omega)}) \subset C^m(h(\Omega)) \to C^0(h(\Omega))$$
is the \emph{Eulerian} form of the formal nonlinear differential operator 
 $v \longmapsto f(\cdot, Lv(\cdot))$  on $h(\Omega)$ , whereas  
$$h^*F_{h(\Omega)}h^{*-1}: h^* D(F_{h(\Omega)})  \subset C^m(\Omega) \to C^0(\Omega)$$ is the {\em Lagrangian}  form of the same operator. The same notation
 is used in other function spaces.  

\par The {\em Eulerian} form is frequently more natural and simpler for computations, but the {\em Lagrangian} form is more adequate to prove theorems, since it allow us to work in fixed function spaces. However, to use  standard tools such as
 the implicit function theorem, we need to show that the map
\[
(u, h) \longmapsto h^*F_{h(\Omega)}h^{*-1}u 
\]
is smooth  and compute its derivatives.
 The differentiability with respect to $u$ poses no problem, since
 $h^*$ is linear. The differentiability with  respect to $h$,
 is the content of the next   result.
 
\begin{prop}
Let $Diff^m(\Omega)$ denote the open subset of maps in 
$\mathcal{C}^m (\Omega, \R^n)$
which are diffeomorphisms to their images (=embeddings).
Then the   map 
\begin{equation*}
(u, h) \to (h^*Lh^{*-1})u :
 Diff^m(\Omega) \times \mathcal{C}^m (\Omega)
 \to \mathcal{C}^0 (\Omega)  \label{lagrange_form_lin}
\end{equation*}
is analytic. If, in addition $f$ is $\mathcal{C}^k$ or analytic, then 
 \begin{equation*}
(u, h) \to h^*F_{h(\Omega)}h^{*-1}u: 
 Diff^m(\Omega) \times \mathcal{C}^m (\Omega)
 \to \mathcal{C}^0 (\Omega)  
\label{lagrange_form}
\end{equation*}
is $\mathcal{C}^k$ or analytic, respectively. 
\end{prop}

\noindent \proof See \cite{He1}. \eproof


\subsection{Lipschitz domains and fractional Sobolev spaces}

\begin{defn}\label{def_lip_domain}
Let $\Omega \subset \R^n$  be an open set. We say that  $\overline{\Omega}$  is a Lipschitz domain (resp.continuously  differentiable,  of class  $\mathbb{C}^{k,1}$ ($k\geq 1$), $m$-times continuously differentiable)  if, for any  $ x \in \partial\Omega$,   there exists a neighborhood    $V$   of   $x$ in $\R^n$  and a map $\psi: V \to\R^n$  such that 
\begin{itemize}
\item [(a)] $\psi$ is injective;
\item [(b)] $\psi$ and $\psi^{-1}$ (defined in  $\psi(V)$) are
 Lipschitz maps  (resp.continuously  differentiable,  of class  $\mathbb{C}^{k,1}$ ($k\geq 1$), $m$-times continuously differentiable);
\item [(c)] $\Omega \cap V = \big\{\,y \in \Omega \,\,\big|\,\, \psi_n(y) < 0\,\big\}$, where $\psi_n(y)$ denotes the 
$n$-th component of  $\psi(y)$.
\end{itemize}
\end{defn}

We now define the (fractional) Sobolev spaces
$W^{s,p}(\Omega)$ in open sets $\Omega$ of $\R^n$.

\begin{defn}\label{sobolev_open_sets} If $\Omega$ is an open set of $\R^n$,  the 
Sobolev space $W^{s,p}(\Omega)$ is the set of all distributions  $u$  in 
$\Omega$ such that  
\begin{itemize}
\item [(a)] $D^\alpha u \in L^p(\Omega)$, for $|\,\alpha\,| \leq m$, when  $s=m$ is a non-negative integer;
\item [(b)] $u \in W^{m,p}(\Omega)$ and  
$$
\begin{array}{l}
\displaystyle\int\displaystyle\int \frac{|\,D^\alpha u(x) - D^\alpha u(y)\,|^{\,p}}{|\,x-y\,|^{\,n\, +\, \sigma p}}dx\,dy < +\infty \\
\Omega \times \Omega
\end{array}
$$
for  $|\,\alpha\,| = m$, when $s= m + \sigma$ is non-negative and non-integer. 
\end{itemize}
\par We define a norm in  $W^{s,p}(\Omega)$ by
$$||\,u\,||_{W^{m,p}(\Omega)}=\left\{ \sum_{|\,\alpha\,| \leq m}\int_\Omega |\,D^\alpha u\,|^{\,p} dx\right\}^\frac{1}{p}$$
in  case  (a), and by
$$
||\,u\,||_{W^{s,p}(\Omega)}=\left\{
\begin{array}{lll}
||\,u\,||_{W^{m,p}(\Omega)}^{\,p} + &
\displaystyle\sum &\displaystyle\int\displaystyle\int \frac{|\,D^\alpha u(x) - D^\alpha u(y)\,|^{\,p}}{|\,x-y\,|^{\,n\, +\, \sigma p}}dx\,dy \\
&|\,\alpha\,| = m & \Omega \times \Omega
\end{array}
\right\}^\frac{1}{p}
$$
 in case (b).
\end{defn}

In what follows, we indicate by $X \hookrightarrow Y$, the continuous inclusion of  the space $X$ into $Y$.

\begin{teo} \label{imbed_wk}
Suppose $\Omega \subset \R^2$ in an open bounded domain.
 If  $p \geq 1$, $kp < 2$  and  $\displaystyle\frac{1}{q} = \displaystyle\frac{1}{p} - \displaystyle\frac{k}{2}$ then  $W^{k,p}(\Omega) \hookrightarrow L^q(\Omega)$.
If $p \geq 1$, $kp = 2$, then $W^{k,p}(\Omega) \hookrightarrow L^q(\Omega)$ for all $q$, $1 \leq q < \infty$.
\end{teo}
\proof  See \cite{Ne}, pg 66.
\eproof

\begin{teo}\label{c}
Suppose $\Omega \subset \R^2$ an open bounded domain. Suppose $p \geq 1$, $kp > 2$,
and put
$$
\begin{array}{ll}
\mu 
\left\{
\begin{array}{lll}
= k - \displaystyle\frac{2}{p} &\mbox{if}& k - \displaystyle\frac{2}{p} < 1\,,
\\
< 1 &\mbox{if}& k - \displaystyle\frac{2}{p} =1\,,
\\
=1 &\mbox{if}& k - \displaystyle\frac{2}{p} > 1.
\end{array}
\right.
\end{array}
$$
Then $W^{k,p}(\Omega) \hookrightarrow  C^{0, \mu}(\overline{\Omega})$.
\end{teo}
\proof See \cite{Ne}, pg 66.
\eproof

\begin{teo}\label{compact_imb_sob}
If $\Omega \subset \R^2$ is a bounded Lipschitz domain and
 $s > t \geq 0$,  {then} the inclusion 
 $I : W^{s,p}(\Omega) \hookrightarrow W^{t,p}(\Omega)$ is compact.
\end{teo}
See \cite{Gri}, pg 26.



\begin{teo}\label{trace_deriv}
Let  $\Omega \subset \R^n$ be an open bounded domain
of class $C^{k,1}$, $k\geq 0$. If $s - 1/p$ is non integer, $s \leq k + 1$,\, $s - 1/p = l + \sigma$,\, $0 < \sigma < 1$, \,$l\geq 0$ integer, then the trace map
$$u \longmapsto \left\{\,\gamma u, \gamma \frac{\partial u}{\partial \nu}, ... , \gamma \frac{\partial^lu}{\partial \nu^l}\,\right\}$$
which is defined for  $u \in C^{k,1}(\overline{\Omega})$, has a  unique continuous extension as a map
$$W^{s,p}(\Omega) \,\,\,\mbox{to} \,\,\,\prod_{j=0}^l W^{s-j-1/p\,,\,p}(\partial\Omega).$$
\end{teo}
 {See \cite{Gri}, pg 37.}

\begin{teo}\label{trace}
Let  $\Omega \subset \R^n$ be an open bounded Lipschitz domain. Then, the trace
 map
  $u \to \gamma u$ which is defined for $u \in C^{0,1}(\overline{\Omega})$, has a unique continuous extension as a map
from $W^{1,p}(\Omega)$ to $W^{1 - 1/p\,,\,p}(\partial\Omega)$.
\end{teo}
 {See \cite{Gri}, pg 38.}

\begin{teo}\label{trace_hk}
 Let  $\Omega \subset \R^2$  be an open bounded Lipschitz domain
 $ 0 \leq k< 1$ and $q = \displaystyle\frac{1}{1-k}$. Then the trace map
 $u \to \gamma u$ which is defined for $u \in C^\infty(\overline{\Omega})$, 
has a unique continuous extension as a map
from  $ H^k(\Omega)$ to  $ L^q(\partial\Omega)$.
\end{teo}
\proof See \cite{Ne}, pg 81.
\eproof



\begin{teo}
Let  $\Omega \subset \R^N$ be a  Lipschitz domain, $u \in W^{1,p}(\Omega)$, $v \in W^{1,q}(\Omega)$, where $\displaystyle\frac{1}{p} + \displaystyle\frac{1}{q} \leq \frac{N+ 1}{N}$ if  $N > p \geq 1$, $N > q \geq 1$ with $q>1$ if $p \geq N$
and $p> 1$ if $q \geq N$. Then
$$\int_\Omega \frac{\partial u}{\partial x_i}v\,dx = \int_{\partial\Omega}uv\nu_i\,dS - \int_\Omega u \frac{\partial v}{\partial x_i}dx\,,$$
where $(\nu_1, \cdots , \nu_N)$ is the outward unit normal, defined a.e.
\end{teo}
\proof See \cite{Ne},  {pg 117}.
\eproof


 When dealing with Lipschitz domains, we face some difficulties concerning the regularity of solutions for elliptic problems (see \cite{Gri} for details).
An appropriate result for our needs is available   in the case of \emph{convex} domains for \emph{strongly elliptic} operators in the divergence form, defined as follows:

\begin{defn}\label{def_st_ell}
A  second  order differential operator $L$ in the domain $\Omega$  in the
 divergence  form
$ Lu = {\displaystyle\sum_{i,j=1}^n D_i
\left(a_{i,j} D_j u \right)}$ 
 where $a_{i,j}$ are \emph{real} constants satisfying  $a_{i,j} = a_{j,i}
\in \mathcal{C}^{0,1} (\overline{\Omega})$  is \emph{strongly elliptic} if
there exists a constant $\alpha >0$, such that
$-{\displaystyle\sum_{i,j=1}^n a_{i,j} {(x)} \xi_i \xi_j \geq \alpha |\xi|^2}$
for any    $\xi=(\xi_1, \cdots,
  \xi_n) \in \R^n$ and $x \in \overline{\Omega}$.
 \end{defn}

\begin{rem}\label{def_st_ell_complex}
Since $a_{i,j} = a_{j,i}$, we also have 
$-{\displaystyle\sum_{i,j=1}^n a_{i,j}(x) \xi_i \xi_j \geq \alpha \left(|\textrm{Re}\,\xi|^2 +
 |\textrm{Im}\,\xi|^2 \right)
 =\alpha |\xi|^2}$
for any    $\xi=(\xi_1, \cdots,
  \xi_n) \in  {\com^n}$, if $L$  is strongly elliptic.
\end{rem}
(For example: $-\Delta u$ is strongly elliptic.)

\begin{teo}\label{existconvex} 
Let  $\Omega \subset \R^n$ be an open bounded and convex domain and
 $ Lu = {\displaystyle\sum_{i,j=1}^n D_i
\left(a_{i,j} D_j u \right)}$ a strongly elliptic operator in the divergence
 form. Then, for any
 $f \in L^2(\Omega)$  {and} $\lambda >0$, there exists a unique  $u \in H^2(\Omega)$ such that 
\begin{equation*}
 -\sum_{i,j=1}^n \displaystyle\int_\Omega a_{i,j}D_iuD_jv\,dx \,+\, \lambda\int_{\Omega} uv \,dx \,=\, \int_\Omega fv \, dx\,,
\end{equation*}
for any $v \in H^1(\Omega)$.
\end{teo}
\proof See \cite{Gri}, pg 149.
\eproof

\subsection{Sectorial operators and semilinear abstract problems}

The results in this section were taken from  \cite{PP}  {and \cite{He2}, where} more details and 
  proofs can be found.

\begin{defn}\label{defset}
 A linear operator  $A$ in a Banach space  $X$ is called \emph{sectorial}
if it  {is} closed, densely defined and there exists $\theta \in \left(0\,,\,\displaystyle\frac{\pi}{2}\right)$, $M \geq 1$ and $a \in \R$, such that the sector
$$S_{a\,,\,\theta}=\left\{\,\lambda \, \big| \, \theta \leq |\,arg(\lambda - a )\,| \leq \pi, \,\lambda \neq a\,\right\}$$
is in the resolvent set
of  $A$, and  
$$||\,(\,\lambda - A\,)^{-1}||\leq \frac{M}{|\,\lambda - a\,|} \,,\,\,\mbox{ for all } \,\,\lambda \in S_{a\,,\,\theta}.$$
\end{defn}

\begin{defn}
If  $A$ is a sectorial operator with $Re \,\sigma(A) >0$ then,
for any  $\alpha >0$ we define the fractional power $A^{-\alpha}$ of $A$, by 
$$A^{-\alpha} = \frac{1}{\Gamma(\alpha)} \int_0^{\infty} t^{\,\alpha -1}e^{-At}dt.$$
\end{defn}

\begin{teo}\label{142}
If $A$ is a sectorial {operator in}  $X$ with $Re \,\sigma(A) >0$, 
then for any  $\alpha >0$, $A^{-\alpha}$ is a bounded linear operator  {on} $X$
which is injective and satisfies  $A^{-\alpha}A^{-\beta} = A^{-(\alpha + \beta)}$,
 when   $\alpha >0$ and $\beta >0$.  Furthermore, for $0 < \alpha <1$ 
$$A^{-\alpha} = \frac{\sin \pi \alpha}{\pi} \int_0^{\infty} \lambda^{-\alpha}(\lambda + A)^{-1}d \lambda.$$
\end{teo}
\proof See \cite{He2} {, pg 25. \eproof}

\begin{defn}
If  $A$ is as in  Theorem \ref{142}, we define $A^\alpha$
 as the inverse of   $A^{-\alpha}$ $(\alpha >0)$ and $D(A^\alpha) = R(A^{-\alpha})$, where  $D(A^\alpha)$ is the domain of  $A^\alpha$  and $R(A^{-\alpha})$
 is the range of   $A^{-\alpha}$.
\end{defn}


\begin{defn}
If  $A$ is a sectorial operator in 
a   Banach space  $X$, we define, for  $\alpha \geq 0$
$$X^\alpha = D(A_1^\alpha)\,\,\,\mbox{with the graph norm}\,,$$
$$||\,x\,||_\alpha = ||\,A_1^\alpha x\,||\,, \, x \in X^\alpha\, ,$$
where $A_1= A + aI$ with $a$ such that  $Re \,\sigma (A_1) > 0$.
\end{defn}

\begin{teo}\label{prop_xalpha}
 If  $A$ is a sectorial operator in a Banach space  $X$, then
 $X^\alpha$ is a Banach space with norm $|| \cdot ||_\alpha$ for $\alpha \geq 0$, $X^0=X$ and, for  $\alpha \geq \beta \geq 0$, $X^\alpha$ is a dense subspace
 of  $X^\beta$ with continuous inclusion. If $A$
has compact resolvent, the inclusion  $X^\alpha \subset X^\beta$
is compact when $\alpha > \beta \geq 0$.
\end{teo}
\proof See \cite{He2} {, pg 29.}
\eproof

\begin{teo}\label{imbed_xalpha}
Let  $\Omega \subset \R^n$ be a Lipschitz domain,
$1<p<\infty$ and  $A$ a sectorial operator in 
 $X=L^p(\Omega)$ with $D(A)= X^1 \subset W^{m,p}(\Omega)$ for some
 $m \geq 1$. Then, for  $0 \leq \alpha \leq 1$ 
$$X^\alpha \subset W^{k,q}(\Omega) \,\,\,\,\mbox{when}\,\,\,\, k - \frac{n}{q} < m\alpha - \frac{n}{p},\,\,\, q \geq p ,$$
$$X^\alpha \subset C^\nu(\Omega) \,\,\,\,\mbox{when}\,\,\,\, 0 \leq \nu < m\alpha - \frac{n}{p}.$$
\end{teo}
\proof
 See \cite{He2}  {pg 39.}
\eproof
 
 \begin{lema} \label{sector}
Suppose  $A $ is a sectorial operator  with  $ \|(\lambda -
A)^{-1}\| \leq  \displaystyle\frac{M}{|\lambda -a |} $ for all $\lambda$ in the
sector $ S_{a,\phi_0} = \{ \lambda  \ | \ \phi_0 \leq
|arg(\lambda-a)| \leq \pi, \lambda \neq a \} $, for some $a\in
\mathbb{R}$ and $0\leq \phi_0<\pi/2$. Suppose also that $B$  is a
linear operator  with  $D(B ) \supset D(A)$ and  $\|Bx - A x\|
\leq \varepsilon \|Ax\| + K \|x\| $, for any  $x \in D(A) $, where
$K$  and  $ \varepsilon $  are positive constants with $
\varepsilon \leq  \displaystyle\frac{1}{4(1+LM)},\, K \leq
\displaystyle\frac{\sqrt{5}}{20M} \frac{\sqrt{2}L - 1}{L^2 - 1} $, for some
$L>1$.

Then  $B$ is also sectorial.  More precisely, if  $b =
\displaystyle\frac{L^2}{L^2 -1} a - \frac{\sqrt{2} L}{L^2 -1} |a| $, $ \phi =
\max \left\{ \phi_0, \displaystyle\frac{\pi}{4}\right\}$ and   $ M' =  2 M \sqrt{5}  $
then
\[
\|(\lambda - B)^{-1}\| \leq  \frac{M'}{|\lambda -b|},
 \]
in the sector  $ S_{b,\phi} = \{ \lambda \ | \ \phi \leq
|arg(\lambda-b)| \leq \pi, \lambda \neq b \} $.
\end{lema}

 {\proof See \cite{PP}, pg 348. \eproof}

\begin{teo} \label{cont}
  Suppose that $A$ is as in Lemma \ref{sector}, $\Lambda$ a topological
space and $\{A_\gamma\}_{\gamma \in \Lambda}$
  is a family of operators in $X$ with $A_{\gamma_0} = A$  satisfying the
following conditions:
\begin{enumerate}
  \item  $D(A_{\gamma} ) \supset D(A)$, for all $\gamma\in \Lambda$;
  \item  $\|A_{\gamma}x - A x\|  \leq \epsilon(\gamma) \|Ax\| + K(\gamma)
\|x\| $ for any $x \in D(A)$,
  where $K(\gamma)$ and $\epsilon(\gamma)$ are positive functions with
$\displaystyle\lim_{\gamma \to \gamma_0} \epsilon(\gamma)=0$
  and $\displaystyle\lim_{ \gamma \to \gamma_0} K(\gamma)  = 0 $.
\end{enumerate}
  Then, there exists a neighborhood $V$ of $\gamma_0$ such that
$A_{\gamma}$ is sectorial if $\gamma \in V$ and
 the family of (linear) semigroups $e^{-tA_{\gamma}}$  {satisfies} 
\begin{gather*}
\|  e^{-tA_{\gamma}}-  e^{-tA} \|   \leq        C(\gamma)   e^{-bt} \\
\|  A \left( e^{-tA_{\gamma}}-  e^{-tA} \right) \|   \leq
C(\gamma) \frac{1}{t} e^{-bt} \\
\| A^{\alpha} \left( e^{-tA_{\gamma}}-  e^{-tA} \right) \|  \leq
C(\gamma)
\frac{1}{ t^{\alpha}} e^{-bt}, \quad  0< \alpha < 1
\end{gather*}
for $t >0$, where $b$ is as in Lemma \ref{sector} and
$C(\gamma) \to 0$  as $\gamma \to \gamma_0$.
\end{teo} 

 {\proof See \cite{PP}, pg 349. \eproof}


 Consider now the abstract evolution equation
\begin{equation}\label{exis}
\begin{array}{lll}
\left\{
\begin{array}{lll}
\displaystyle\frac{dx}{dt} +  Ax = f(t,x) \, , \, t>t_0\,,
\\
x(t_0)=x_0 
\end{array}
\right .
\end{array}
\end{equation}
where  $A$ is a sectorial operator in $X$ and   $ f:  U \subset \R \times X^\alpha   \to X $, for some $\alpha \in [0,1)$.

\begin{defn}
  A solution of (\ref{exis}) in $[t_0, t_1)$ is a continuous function $x:[t_0,t_1) \to X^\alpha$, $x(t_0)=x_0$, such that
      $f( \cdot , x(\cdot)):[t_0,t_1) \to X$ is continuous 
 and for $t\in (t_0,t_1)$
 $\displaystyle\frac{dx}{dt}$ exists, $ {(t, x(t))} \in U$,  
 $x(t) \in D(A)$, and $x$ satisfies (\ref{exis}).
\end{defn}

\begin{teo}\label{333}
  Suppose  $A$ is a sectorial operator, $0 \leq \alpha < 1$,
  $U \subset \R \times X^\alpha$  is open and  $f: U \to X$
  is locally Lipschitz continuous in  $x$, and locally  {H$\ddot{o}$lder} continuous in $t$.
  Then, for any  $(\,t_0,x_0\,) \in U$ there exists
  $T=T(t_0,x_0) > 0$ such that  (\ref{exis}) has a unique solution
   $x(t;t_0,x_0)$  {on} $(\,t_0,t_0+T\,)$ with initial value $x(t_0)=x_0$, which is continuous in $x_0$. 
\end{teo}
\proof See \cite{He2}, Theorems 3.3.3 and 3.4.1.
\eproof

\begin{teo}\label{334}
 Suppose  $A$ and $f$  satisfy the hypotheses of Theorem \ref{333}
 and  also that   $f(B)$ is bounded in $X$, for any closed bounded set $B \subset U$.   Then, if  $x$ is  {a solution} of  (\ref{exis})  {on} $(\,t_0,t_1\,)$ and $t_1$   {is} maximal, 
then either  $t_1 = +\infty$ or there exists a sequence  $t_n \rightarrow t_{1^-}$ as  $n \rightarrow +\infty$ such that $(t_n,x(t_n)) \rightarrow \partial U$.
\end{teo}

\proof See \cite{He2}, Theorem 3.3.4.
\eproof

\begin{teo}\label{336}
  Suppose $A$ and $f$ as in Theorem \ref{333} and also assume
  that  $A$ has compact resolvent  and $f$ takes all sets $\R^+ \times B \subset U \subset \R \times X^\alpha$  with  $B$ bounded and closed
 into bounded sets of  $X$. If $x(t;t_0,x_0)$ is a solution of  (\ref{exis}) in $(\,t_0, \infty\,)$ with  $||\,x(t;t_0,x_0)\,||_\alpha$ bounded as  $t \rightarrow + \infty$, then  $\{\,x(t;t_0,x_0)\,\}_{t> t_0}$ is in a compact   {set} of  $X^\alpha$.
\end{teo}
 \proof See \cite{He2}, Theorem 3.3.6.
\eproof

\subsection{Global attractors and gradient systems} 

\begin{defn}
 A compact subset   $\mathcal{A} \subset X$ which is invariant  under the action of a
  $C^r$ semigroup  $\big\{T(t)\big\}_{t \geq 0}$
  with  $r\geq 0$  
   is called maximal  if any compact invariant subset under $T(t)$ is contained
 in   $\mathcal{A}$.  A compact invariant set  $\mathcal{A}$
 is  called a \emph{global attractor} if   $\mathcal{A}$ is invariant, maximal
 and attracts each bounded subset  $B \subset X$. 
\end{defn}

\begin{defn}
 A point $ {x \in X}$ is called an \emph{equilibrium point} of the  $C^r$ 
 semigroup
 $\big\{T(t)\big\}_{t \geq 0}$ if  $T(t)x=x$ for  $t\geq 0$.  An   {equilibrium} point $x$ is \emph{hyperbolic}
 if   $\sigma (DT(t)(x))$ does not   {intersect} the unit circle with centered at  the origin.\end{defn}

\begin{defn}
 A semigroup  $T(t): X \to X, t \geq 0$ of class $C^r$,  $r \geq 1$, is called a 
\emph{gradient system}  (or \emph{gradient flow})  if  
\begin{itemize}
\item [(i)] Each  bounded positive orbit
 is pre-compact.
\item [(ii)]  There exists a Lyapunov functional for  $T(t)$, that is,  a continuous function  $V: X \to \R$ with the following properties:
\item [(${ii}_1$)] $V(x)$ is bounded below,
\item [(${ii_2}$)] $V(x) \rightarrow \infty$ when $|\,x\,| \rightarrow \infty ,$
\item [(${ii_3}$)] $V(T(t)x)$ is decreasing in  $t$ for each  $x \in X$,
\item [(${ii_4}$)]  $T(t)x$ is defined for  $t \in \R$  and  $V(T(t)x)=V(x)$ 
for  $t \in \R$, then  $x$  is an   {equilibrium} point.
\end{itemize}
\end{defn}

\begin{teo}\label{385}
If  $\big\{T(t)\big\}_{t \geq 0}$ is a gradient system 
  {asymptotically} smooth and the set $E$ of equilibria is bounded, then there exists a  global attractor  $\mathcal{A}$ for  $T(t)$  and  $\mathcal{A}=W^u(E)=\big\{y \in X : T(-t)y \,\,\mbox{is defined for  } \,\,t\geq 0 \,\,\,\mbox{and} \,\,\,T(-t)y \rightarrow E \,\,\,\mbox{when} \,\,\,t \rightarrow \infty\,\big\}$. If
 $X$  is a Banach space, then  $\mathcal{A}$  is connected. If, in addition
 each point in  $E$  is hyperbolic, then  $E$ is finite and
   $$\mathcal{A}= \bigcup_{x \,\in\, E} W^u(x).$$
\end{teo}
\proof
 See \cite{Hal}, Theorem 3.8.5.
\eproof



 \begin{lema} \label{continabsALT}  
Suppose $Y$ is a metric space, $ \Lambda$ is  an open set in $Y$,
$\{- A_\lambda\}_{\lambda \in \Lambda}$ is a family of operators in
a Banach space $X$ satisfying the conditions of Theorem \ref{cont}
at $ \lambda = \lambda_0$,  $ U$ is an open set in $\mathbb{R}^+
\times X^{\alpha}$,  $ 0 \leq \alpha < 1$  and $f:U \times \Lambda
\to  X$ is H\"older continuous in $t$.
Suppose also that, for any bounded subset $D \subset U$, $f$ is 
 continuous in $\lambda$ at $ \lambda_0$ uniformly for $(t,x)$ in
 $D$, and there is a constant $L=L(D)$, such that 
   $ \| f(t,x,\lambda) - f(t,y,\lambda) \|
\leq L \| x-y \|_{\alpha}$
  for $(t,x)$, $(t,y)$ in $D$ and
$\lambda \in \Lambda$.
  Suppose further that  the solutions  $x(t,x_0, \lambda)$  
  of the problem
$$
\begin{gathered}
\frac{dx}{dt}   =   A_{\lambda}x + f(t,x,\lambda), \quad t> t_0  \\
  x(t_0)  =  x_0
\end{gathered}
$$
 exist and 
  remain in a bounded subset of $X^{\alpha}$  when   $x_0$ varies
in a  bounded subset of $X^{\alpha}$, 
 $\lambda $ in a neighborhood of $\lambda_0$   and $t_0 \leq t \leq T$.

 Then the function $ \lambda
\mapsto  x(t,x_0,\lambda)  \in X^{\alpha} $ is continuous   at $\lambda_0$
uniformly for $x_0$ in bounded subsets of $X^{\alpha}$
   and  $t_0 \leq t \leq T$ .
\end{lema}
{\proof See \cite{PP}, Lemma 3.7. \eproof}

 We recall now that the family of subsets  $\mathcal{A}_{\lambda}$ of
 a metric space  $(X,d)$ is
 said to be  {\em upper-semicontinuous} at $\lambda = \lambda_0$
 if $\delta(\mathcal{A}_{\lambda}, \mathcal{A}_{\lambda_0}) \to 0$
 as $ \lambda \to \lambda_0$, where
 $\delta(A,B) = \displaystyle\sup_{x\in A} d(x,B) =
\displaystyle\sup_{x\in A} \displaystyle\inf_{y \in B} d(x,y) $ 
and {\em lower-semicontinuous} if 
 $\delta(\mathcal{A}_{\lambda_0}, \mathcal{A}_{\lambda}) \to 0$
 as $ \lambda \to \lambda_0$.

\begin{teo} \label{upper}
 Suppose $Y$ is a metric space, $ \Lambda$ is  an open set in $Y$,
$\{-A_\lambda\}_{\lambda \in \Lambda}$ is a family of operators in
a Banach space $X$ satisfying the conditions of Theorem \ref{cont}
at $ \lambda = \lambda_0$ with $b>0$,  $ U$ is an open set in 
$ X^{\alpha}$,  $ 0 \leq \alpha < 1$  and $f:U \times \Lambda
\to  X$  satisfying the conditions of Lemma
 \ref{continabsALT}. Let $T(t, \lambda)(x_0)$ be the nonlinear
 semigroup in $X^{\alpha}$ given by the solutions of the problem
\begin{equation}\label{autonomous}
\begin{gathered}
\frac{dx}{dt}   =   A_{\lambda}x + f(x,\lambda), \quad t> t_0  \\
  x(t_0)  =  x_0.
\end{gathered}
\end{equation}
 Suppose also that, for each $\lambda \in \Lambda$ there
 exists a global compact attractor  $\mathcal{A}_{\lambda} $
 for  $ T(t, \lambda)(x)$, the union 
 $\displaystyle\bigcup_{\lambda \in \Lambda} \mathcal{A}_{\lambda}$  is a bounded set
  in $X$  and  $f$ maps this  union 
   into a
 bounded set of $X$. Then the family  $ \mathcal{A}_{\lambda}$
 is upper semicontinuous at $\lambda=\lambda_0$.
\end{teo} 
 {\proof See \cite{PP}, Theorem 3.9. \eproof}

Lower semicontinuity can also be proved under some additional hypotheses.
   
	\begin{teo} \label{lower}
   Suppose, in addition to the hypotheses of Theorem \ref{upper}, that
    the system generated by (\ref{autonomous})  is gradient  for all $\lambda$, 
    its equilibria are all hyperbolic and continuous in $\lambda$ and
     the local invariant manifolds of the equilibria are continuous in $\lambda$.
      Then the family  $ \mathcal{A}_{\lambda}$
   is also lower semicontinuous at $\lambda=\lambda_0$.
   \end{teo}
\proof See \cite{PP}, Theorem 4.10. \eproof
\section{Reduction to a fixed domain} \label{reduction}

 Let $\Omega$ be the unit square in $\R^2$, and 
consider the family of maps $h_{\epsilon}: \Omega \to \R^2 $, defined by  (\ref{per}).   For simplicity
 we will denote by $\Omega_{\epsilon}$  the corresponding family of ``perturbed domains'' $ \Omega_{h_\epsilon} = h_{\epsilon}(\Omega) $. 
  We first establish some basic properties of these families.

\begin{lema}\label{estimateh}
If $\epsilon> 0$ is sufficiently small, the map   $h_\epsilon$ defined by
 (\ref{per}) belongs to  $Diff^m(\Omega)$, for any $m \geq 1$ and 
 $||\,h_\epsilon - i_\Omega\,||_{\,C^1(\Omega)} \to 0 $,
 as $\epsilon \to 0 $.
\end{lema}
\proof  It is clear that $h   \in \C^m(\Omega)$ and is injective and
  $\displaystyle\frac{1}{|\,Jh_\epsilon(x)\,|}=\displaystyle\frac{1}{|\,1 + \epsilon\,sen(x_1/\epsilon^{\,\alpha})\,|}$ is bounded in $\Omega$ if $\epsilon$ is sufficiently small, proving that $h_{\epsilon} \in Diff^m(\Omega). $
Now, a simple computation shows that  

 $$||h_\epsilon - i_\Omega||_{\,C^1(\Omega)} = \max \left\{\displaystyle\sup_{(x_1,x_2)\,\in\,\Omega}\big|x_2\epsilon sen(x_1/\epsilon^\alpha)\big| , \displaystyle\sup_{(x_1,x_2)\in\Omega} \big(x_2^2\epsilon^{2 - 2\alpha}cos^2(x_1/\epsilon^\alpha) + \epsilon^2 sen^2(x_1/\epsilon^\alpha)\big)^\frac{1}{2}\right\}.$$ Therefore, if  $0 < \alpha < 1$, it follows that
  $||\,h_\epsilon - i_\Omega\,||_{C^1(\Omega)} \to 0$ as  $\epsilon \to 0$. 
\eproof

\begin{rem}
The hypothesis $\alpha < 1$ is  essential in the above proof, we do not have
 $\C^1$ convergence of $h_{\epsilon}$ if $\alpha \geq 1$. Also, convergence
 in the  $\C^2$ norm only holds if $\alpha < \frac{1}{2}$.
\end{rem}

\begin{lema} \label{isosobolev}
If  $  s> 0$
 and   $\epsilon >0$ is small enough, the map 
$$
\begin{array}{llcll}
h_\epsilon^* &:& H^s(\Omega_{\epsilon}) &\to& H^s(\Omega) \\
&& u & \longmapsto & u \circ h_\epsilon
\end{array}
$$
is an isomorphism, with inverse ${h_\epsilon^*}^{-1} = (h_\epsilon^{-1})^*$.
\end{lema}

\proof We restrict ourselves to the case $ 0 \leq s < 1$, the proof for $s \geq 1$ is similar and is left to the reader.  
 From  Lemma \ref{estimateh}   it is clear that   $h_\epsilon^*$ is  invertible with inverse   ${h_\epsilon^*}^{-1} = (h_\epsilon^{-1})^*$
if  $\epsilon >0$ is sufficiently small.
Now, since the derivatives of $h_\epsilon^{-1}$ are bounded below by a constant
 $M$, we have that
$|h_\epsilon^{-1}(r)-h_\epsilon^{-1}(w)| \geq M |r-w| $, uniformly, for $ {(r,w)}  \in {\Omega_{\epsilon}  } \times {\Omega_{\epsilon}  }$. It follows that 

\begin{eqnarray*}
||\,h_\epsilon^*u\,||_{H^s(\Omega)} &:= & \bigg\{ \int_{\Omega} |\,(u \circ h_\epsilon)(x)\,|^{\,2}dx  + \int_{\Omega} \int_{\Omega} \frac{\big|\,(u \circ h_\epsilon)(x) - (u \circ h_\epsilon)(y)\,\big|^{\,2}}{|\,x-y\,|^{\,2 + 2s}}\,dx\,dy\, \bigg\}^\frac{1}{2}  \\
&= & \bigg\{\int_{\Omega_\epsilon} |\,u(r)\,|^{\,2}|\,Jh_\epsilon^{-1}(r)\,|\,dr \,\, \\
& +&
\int_{\Omega_{\epsilon}  }  \int_{\Omega_{\epsilon}  } \frac{\big|\,u(r) - u(w)\,\big|^{\,2}}{|\,h_\epsilon^{-1}(r)-h_\epsilon^{-1}(w)\,|^{\,2 + 2s}}\,|\,Jh_\epsilon^{-1}(r)\,|\,|\,Jh_\epsilon^{-1}(w)\,|\,dr\,dw \,\bigg\}^\frac{1}{2}  \\
&= &\bigg\{\int_{\Omega_{\epsilon}  } |\,u(r)\,|^{\,2}|\,Jh_\epsilon^{-1}(r)\,|\,dr \\
 & + & \int_{\Omega_{\epsilon}  }\int_{\Omega_{\epsilon}  }
 \frac{\big|\,u(r) - u(w)\,\big|^{\,2}}{|\,r - w\,|^{\,2 + 2s}}\displaystyle\frac{|\,r - w\,|^{\,2 + 2s}}{|h_\epsilon^{-1}(r)-h_\epsilon^{-1}(w)|^{\,2 + 2s}}|Jh_\epsilon^{-1}(r)||Jh_\epsilon^{-1}(w)|drdw \bigg\}^\frac{1}{2}  \\
& \leq & \bigg\{\int_{\Omega_{\epsilon}  } |\,u(r)\,|^{\,2}|\,Jh_\epsilon^{-1}(r)\,|\,dr \\
 & + & \int_{\Omega_{\epsilon}  }\int_{\Omega_{\epsilon}  }
 \frac{\big|\,u(r) - u(w)\,\big|^{\,2}}{|\,r - w\,|^{\,2 + 2s}}
\frac{1}{M^{\,2 + 2s}}|Jh_\epsilon^{-1}(r)||Jh_\epsilon^{-1}(w)|drdw \bigg\}^\frac{1}{2}  \\
& \leq & \bigg\{{K}\int_{\Omega_{\epsilon}  } |\,u(r)\,|^{\,2}\,dr 
  +  \frac{K^2}{{(2M)^{2 + 2s}}} \int_{\Omega_{\epsilon}  }\int_{\Omega_{\epsilon}  }
 \frac{\big|\,u(r) - u(w)\,\big|^{\,2}}{|\,r - w\,|^{\,2 + 2s}}
 \, drdw \bigg\}^\frac{1}{2},  
\end{eqnarray*}
where $K$ is a bound for $|\,Jh_\epsilon^{-1}(r)\,|$ in $\Omega_{\epsilon}$.
A similar argument shows that 
${h_\epsilon^*}^{-1}= {h_\epsilon^*}^{-1}$ is also bounded.
 \eproof

 Let  $\Delta_{\Omega_{\epsilon}}$ be the Laplacian operator in the region 
 $\Omega_{\epsilon}= h_\epsilon(\Omega)$. We want to find an expression for the
  differential operator  $h_\epsilon^*\Delta_{\Omega_{\epsilon}}h_\epsilon^{*-1}$ in the fixed region  $\Omega$, in terms of $h_\epsilon$. 
  Writing
 $h_\epsilon(x)\,=\,h_\epsilon(x_1,x_2)\,=\,((h_\epsilon)_1(x),(h_\epsilon)_2(x))\,=\,(y_1,y_2)\,=\,y\,$, we obtain, for  $i=1,2$ 
\begin{equation}
\begin{split} \label{deriv}
\left(h_\epsilon^* \displaystyle\frac{\partial}{\partial y_i} h_\epsilon^{*-1}(u)\right)(x) &=  \displaystyle\frac{\partial}{\partial y_i}(u\circ h_\epsilon^{-1})(h_\epsilon(x))  \\
&= \displaystyle\frac{\partial u}{\partial x_1}(h_\epsilon^{-1}(y))
\displaystyle\frac{\partial (h_\epsilon)_1^{-1}(y)}{\partial y_i}(y) +
\displaystyle\frac{\partial u}{\partial x_2}(h_\epsilon^{-1}(y))
\displaystyle\frac{\partial (h_\epsilon)_2^{-1}(y)}{\partial y_i}(y)   \\
&= \displaystyle\sum^2_{j=1} \left[\left(\displaystyle\frac{\partial h_\epsilon}{\partial {x_j}}\right)^{-1}\right]_{j,i}(x)\frac{\partial u}{\partial x_j}(x) \\
&= \displaystyle\sum^2_{j=1} b_{ij}(x)\displaystyle\frac{\partial u}{\partial x_j}(x)\,,
\end{split}
\end{equation}

where  $b_{ij}(x)$ is the   $i,j$-entry of the
inverse transpose of the  jacobian matrix  of  $h_\epsilon$, {given by}
 
\begin{equation*}
\begin{array}{ccccl}
[{h_\epsilon}^{-1}]^{\,T}_{\,x}=
\begin{array}{ccccl}
\left(
\begin{array}{cc}
1 & \displaystyle\frac{-x_2\,\epsilon^{1-\alpha}cos(x_1/\epsilon^\alpha)}{1+ \epsilon \, sen(x_1/\epsilon^\alpha)}\\
0&
\displaystyle\frac{1}{1+ \epsilon \, sen(x_1/\epsilon^\alpha)}
\end{array}
\right)\,.
\end{array}
\end{array}
\end{equation*}

  Therefore,
\begin{equation}
 \begin{split} \label{laph}
 h_\epsilon^*\Delta_{\Omega_{\epsilon}}h_\epsilon^{*-1}(u)(x)  & =   
 \sum_{i=1}^2 \left(h_\epsilon^*\frac{\partial^2}{\partial y_i^2}h_\epsilon^{*-1}(u)\right)(x)   \\
 & =    \sum_{i=1}^2 \left(\sum_{k=1}^2 b_{i\,k}(x)\frac{\partial}{\partial x_k}\left(\sum_{j=1}^2 b_{i j}\frac{\partial u}{\partial x_j}\right) \right)(x)   \\
 &=  div \left(\displaystyle\frac{\partial u}{\partial x_1} + b_{12}\displaystyle\frac{\partial u}{\partial x_2}\,,\, \left(b_{12}^{\,2} + b_{22}^{\,2} \right) \displaystyle\frac{\partial u}{\partial x_2} + b_{12}\displaystyle\frac{\partial u}{\partial x_1}\right) (x)   \\
 & -    \displaystyle\frac{\partial b_{12}}{\partial x_2}\displaystyle\frac{\partial u}{\partial x_1} (x) - b_{12}\displaystyle\frac{\partial b_{12}}{\partial x_2}\displaystyle\frac{\partial u}{\partial x_2} (x) \,.
  \end{split}
\end{equation}

 We also need to compute the boundary condition $h_\epsilon^*\displaystyle\frac{\partial}{\partial N_{\Omega_{\epsilon}}}h_\epsilon^{*-1}u=0$  in the fixed region $\Omega$ in terms of  $h_\epsilon$.  Let  $N_{h_\epsilon(\Omega)}$ denote the {outward unit} normal to the boundary of  $h_\epsilon(\Omega):=\Omega_{\epsilon}$. From  (\ref{deriv}), we obtain 
\begin{equation}
\begin{split}
\left(h_\epsilon^*\displaystyle\frac{\partial}{\partial N_{\Omega_{\epsilon}}}h_\epsilon^{*-1}u\right)(x) &=  \sum_{i=1}^2
\left(h_\epsilon^*\displaystyle\frac{\partial}{\partial y_i}h_\epsilon^{*-1}u\right)(x)\left(N_{\Omega_{\epsilon}}\right)_i(h_\epsilon (x))   \\
&=  \sum_{i=1}^2 \displaystyle\frac{\partial}{\partial y_i}(u \circ h_\epsilon^{-1})(h_\epsilon (x))\left(N_{\Omega_{\epsilon}}\right)_i(h_\epsilon (x))   \\
&=  \sum_{i,j=1}^2 b_{ij}(x)\displaystyle\frac{\partial u}{\partial x_j}(x)\left(N_{\Omega_{\epsilon}}\right)_i(h_\epsilon(x))   \\
&=  N_{\Omega_{\epsilon}}(h_\epsilon (x))\cdot (u_{x_1}(x) + b_{12}(x)u_{x_2}(x)\,,\,b_{22}(x)u_{x_2}(x))\,. \label{normalh}
\end{split}
\end{equation}
\par Since
 $$h_\epsilon^*N_{\Omega_{\epsilon}}(x) = N_{\Omega_{\epsilon}}(h_\epsilon(x)) = \displaystyle\frac{[h_\epsilon^{-1}]_x^T N_\Omega(x)}{||\,[h_\epsilon^{-1}]_x^T N_\Omega(x)\,||} 
 = \displaystyle\frac{\big((N_\Omega(x))_1 + b_{12}(x)(N_\Omega(x))_2\, ,\, b_{22}(x)(N_\Omega(x))_2\big)}{||\,[h_\epsilon^{-1}]_x^T N_\Omega(x)\,||}$$
 (see \cite{He1}), we obtain from (\ref{normalh}) 
$$
\left(h_\epsilon^*\displaystyle\frac{\partial}{\partial N_{\Omega_{\epsilon}}}h_\epsilon^{*-1}u\right)(x) =
\displaystyle\frac{\left(N_\Omega(x)\right)_1\big[u_{x_1}(x) + b_{12}(x)u_{x_2}(x)\big] + \left(N_\Omega(x)\right)_2\left[b_{12}(x)u_{x_1}(x) + \big(b^2_{12}(x) + b_{22}^{2}(x)\big)u_{x_2}(x)\right]}{||\,[h_\epsilon^{-1}]_x^T N_\Omega(x)\,||}.
$$

 Thus, the boundary condition
$
\left(h_\epsilon^*\displaystyle\frac{\partial}{\partial N_{\Omega_{\epsilon}}}h_\epsilon^{*-1}u\right)(x) = 0\,, 
$
becomes
$$
\left(N_\Omega(x)\right)_1\left[\,-u_{x_1}(x)\,-\,b_{12}(x)u_{x_2}(x)\,\right] + \left(N_\Omega(x)\right)_2\left[\,-b_{12}(x)u_{x_1}(x)\,-\, \big(\,b^{\,2}_{12}(x) + b_{22}^{\,2}(x)\,\big)u_{x_2}(x)\,\right] = 0\,, 
$$
which can be written as 
$$
\sum_{i,j=1}^2 \left(N_\Omega(x)\right)_i(c_{ij}D_ju) = 0 \,\,\mbox{on}\,\, \partial\Omega\,,
$$
 where 
\begin{equation}\label{ces}
c_{11}=-1\,,\,\,c_{12}= c_{21} = -b_{12}\,,\,\,c_{22}=-(\,b_{12}^{\,2} + b_{22}^{\,2}\,)\,.
\end{equation}

 \vspace{3mm}

Now, observe that 
$v(.\,,t)$  is a solution  (\ref{nonlinBVP}) in  the perturbed region
 ${\Omega_{\epsilon} = h_\epsilon(\Omega)}$, if and only if   
 $u(.\,,t)={h_\epsilon^*v(.,t)}$  satisfies
\begin{equation}\label{nonlinBVP_fix}
\begin{array}{rcl}
\left\{
\begin{array}{rcl}
u_t(x,t)&=& {h_\epsilon^*\Delta_{\Omega_{\epsilon}}h_\epsilon^{*^{-1}}} u(x,t) -au(x,t) + f(u(x,t)) ,\,\, x \in \Omega \,\,\,\mbox{and}\,\,\,t>0, \\
{h_\epsilon^*\displaystyle\frac{\partial }{\partial N_{\Omega_{\epsilon}}}h_\epsilon^{*^{-1}}}u(x,t)&=&g(u(x,t)), \,\, x \in \partial\Omega \,\,\,\mbox{and}\,\,\,t>0\,,
\end{array}
\right.
\end{array}
\end{equation}

in the \emph{fixed} region $\Omega$.



 \section{Sectoriality of the perturbed operators} \label{sectorialoper}

 In this section we show that the family of  differential operators
$   - h_\epsilon^*\Delta_{\Omega_{\epsilon}}h_\epsilon^{*-1} $, appearing in 
 (\ref{nonlinBVP_fix}), generate sectorial operators in various spaces.

\subsection{Sectoriality in $L^2$}

 Consider the operator

 \begin{equation}\label{opeper}
A_{\epsilon}:= \left(\,- h_\epsilon^*\Delta_{\Omega_{\epsilon}}h_\epsilon^{*-1} +aI\,\right): L^2(\Omega) \to L^2(\Omega)\,
\end{equation}
 with domain
\begin{equation}\label{dominio}
D\left(A_{\epsilon}\right) = \left\{u \in H^2(\Omega) \,\bigg|\, 
  h_{\epsilon}^*\displaystyle\frac{\partial }{\partial N_{\Omega_{\epsilon}}}h_{\epsilon}^{*^{-1}}u = 0,\, \textrm{ on }  \,  \partial \Omega \right\}.
\end{equation}
(We will  denote simply by $A$  the unperturbed operator  
$\left(\,{-}\Delta_{\Omega}  +aI\,\right)$).

 In the case of smooth domains, it is not difficult to prove that $A_{\epsilon}$ 
 is sectorial, for $\epsilon$ sufficiently small, but for  general Lipschitz domains, we need to  address  some delicate questions of regularity. Fortunately, in our case, Theorem \ref{existconvex} can be used, since we are dealing with a convex domain. 
However, this result is applicable only  for operator in the divergence form. For this reason, using {(}\ref{laph}{)}  we   write our operator as sum 
 \begin{equation} \label{decomposition}
  A_{\epsilon}=
\left(\,C_{\epsilon} +aI  +  L_{\epsilon} \,\right)
\end{equation}
where 
\begin{equation}\label{ch}
C_{\epsilon}u=\sum_{i,j=1}^2 D_i\,(c_{ij}D_ju)
\end{equation}
 with the  $c_{ij}$ given in  (\ref{ces}),
and 
 \begin{equation*}
L_{\epsilon}u:=\displaystyle\frac{\partial b_{12}}{\partial x_2}\displaystyle\frac{\partial u}{\partial x_1} + b_{12}\displaystyle\frac{\partial b_{12}}{\partial x_2}\displaystyle\frac{\partial u}{\partial x_2}.
\end{equation*}

We now want to show that, if $\epsilon$ is small, the operator defined by
(\ref{opeper}) and (\ref{dominio}) is sectorial.
  To this end, we need some auxiliary results for the first term in the decomposition 
 (\ref{decomposition}).

\begin{lema}\label{fe}
 If $\epsilon > 0$ is sufficiently small, the differential operator 
 $C_{\epsilon}$ given by (\ref{ch}) is strongly
 elliptic.
\end{lema}
\proof If   $\epsilon >0$ is sufficiently small, we have
$c_{11}= -1$, $-\displaystyle\frac{1}{4} < c_{12}=c_{21}=-b_{12} < \displaystyle\frac{1}{4}$  
and  $c_{22}= -(\,b_{12}^{\,2} + b_{22}^{\,2}\,)< -\displaystyle\frac{1}{2}$.
Therefore, if   $\xi=(\xi_1, \xi_2) \in \R^2$, we have

\begin{equation}\label{elliptic}
\begin{split} 
\sum_{i,j=1}^2c_{ij}\xi_i\xi_j  \leq & -\displaystyle\frac{1}{2}\xi_1^{\,2} + \displaystyle\frac{1}{2}|\,\xi_1\,|\,|\,\xi_2\,| - \displaystyle\frac{1}{2}\xi_2^{\,2} \\
&= -\displaystyle\frac{1}{4}(\,\xi_1^2 + \xi_2^{\,2}\,) - \left(\,\displaystyle\frac{1}{4}\xi_1^{\,2} - \displaystyle\frac{1}{2}|\,\xi_1\,|\,|\,\xi_2\,| + \displaystyle\frac{1}{4}\xi_2^{\,2}\,\right)  \\
&=  -\displaystyle\frac{1}{4}(\,\xi_1^{\,2} + \xi_2^{\,2}\,) - \left(\,\displaystyle\frac{1}{2}|\,\xi_1\,| - \displaystyle\frac{1}{2}|\,\xi_2\,|\,\right)^2 \\
& \leq  -\displaystyle\frac{1}{4}(\,\xi_1^{\,2} + \xi_2^{\,2}\,)\,. 
\end{split}
\end{equation}

\eproof
 
   \begin{lema}\label{sim}
If $\epsilon {>0}$ is small enough, the operator $C_{\epsilon}$
in $L^2(\Omega)$ defined by ${(}\ref{ch}{)}$, with domain  
\begin{equation}\label{domicha}
D(C_{\epsilon}) = \left\{u \in H^2(\Omega) \,\bigg|\, \displaystyle\sum_{i,j=1}^2 \left(N_\Omega(x)\right)_i(c_{ij}D_ju) = 0 \,,\, x \in \partial\Omega\right\},
\end{equation} is symmetric and bounded below.
\end{lema}
\proof  Using integration by parts, we have, for any  $u, v \in D(C_{\epsilon} )$
$$
\begin{array}{lll}
\left\langle C_{\epsilon}u \,,\, v \right\rangle_{L^2(\Omega)} 
&=& \displaystyle\int_{\partial\Omega} v\left[\,\displaystyle\sum_{i,j=1}^2 \left(N_\Omega(x)\right)_i(c_{ij}D_ju)\,\right]\,d\sigma (x)
- \displaystyle\int_\Omega\sum_{i,j=1}^2 (c_{ij}D_juD_iv) \,dx 
\\
&=& - \displaystyle\int_\Omega\sum_{i,j=1}^2 (c_{ji}D_ivD_j{u}) \,dx 
 \\
&=& \left\langle u\,,\, C_{\epsilon}v 
\right\rangle_{L^2(\Omega)},
\end{array}
$$
proving that {the operator} $C_{\epsilon}$ is symmetric. 
 Now, from  (\ref{elliptic}), it follows that 

\begin{align}\label{estCh}
\left\langle C_{\epsilon}u \,,\, u \right\rangle_{L^2(\Omega)} =& \displaystyle\int_\Omega \sum_{i,j=1}^2 D_i\,(c_{ij}D_ju)\cdot u \,dx \nonumber \\
 =& - \displaystyle\int_\Omega\sum_{i,j=1}^2 (c_{ij}D_juD_iu) \,dx \nonumber \\
\geq&  \displaystyle \frac{1}{4}\int_\Omega \,|\,\nabla u\,|^{\,2} \,dx 
\end{align} 
so  $C_{\epsilon}$ is bounded below. \eproof

\begin{rem}\label{domicha_complex}
If $u$ and $v$ are complex valued
with real and imaginary parts in  $ D(C_{\epsilon} )$, we still have 
(see remark \ref{def_st_ell_complex}),
 \begin{align}\label{estCh_complex}
\left\langle C_{\epsilon}u \,,\, u \right\rangle_{L^2(\Omega)} =& \displaystyle\int_\Omega \sum_{i,j=1}^2 D_i\,(c_{ij}D_ju)\cdot \overline{u} \,dx \nonumber \\
 =& - \displaystyle\int_\Omega\sum_{i,j=1}^2 (c_{ij}D_ju \overline{D_iu}) \,dx \nonumber \\
\geq&  \displaystyle \frac{1}{4}\int_\Omega \,|\,\nabla u\,|^{\,2} \,dx .
\end{align}
\end{rem} 

\begin{lema}\label{exisol}
If  $\epsilon >0$ is sufficiently small and $a>0$, then the problem
\begin{equation}\label{ca}
\begin{array}{lll}
\left\{
\begin{array}{rcl}
C_{\epsilon}u + au  & = & f,\,\,\ x \in \Omega 
\\
\displaystyle\sum_{i,j=1}^2 \left(N_\Omega(x)\right)_i(c_{ij}D_ju) & =& 0 , \,\, x \in \partial\Omega 
\end{array}
\right.
\end{array} 
\end{equation}
 has a unique solution  $u \in H^2(\Omega)$, for any   
$f \in L^2(\Omega)$.
\end{lema}
\proof 
Since  $\Omega$ is a bounded and convex domain, we have by Theorem
 \ref{existconvex} that, for any   $f \in L^2(\Omega)$ and $a > 0$,
 there exist a unique solution  $u \in H^2(\Omega)$ such that  
\begin{equation*}
- \sum_{i,j=1}^2 \displaystyle\int_\Omega c_{ij}D_juD_iv\,dx \,+\, a\int_{\Omega} uv \,dx \,=\, \int_\Omega fv \, dx\,,
\end{equation*}
for any  $v \in H^1(\Omega)$. Integrating by parts, we obtain
\begin{equation*}
\int_{\Omega} \sum_{i,j=1}^2 D_i\,(c_{ij}D_ju)v\,dx - \displaystyle\int_{\partial\Omega}\displaystyle\sum_{i,j=1}^2 \left(N_\Omega(x)\right)_i(c_{ij}D_ju)v\,d\sigma(x) \,+\, a\int_{\Omega} uv \,dx \,=\, \int_\Omega fv \, dx\,.
\end{equation*}
Therefore, for any  $v \in H^1_0(\Omega)$, we obtain

\begin{equation*}
\int_{\Omega} \sum_{i,j=1}^2 D_i\,(c_{ij}D_ju)v\,dx \,+\, a\int_{\Omega} uv \,dx \,=\, \int_\Omega fv \, dx \,,
\end{equation*}
from which the first equation in (\ref{ca}) follows immediately.


Thus, for any $v \in H^1(\Omega)$, we obtain
\begin{equation*}
\int_{\partial \Omega} \displaystyle\sum_{i,j=1}^2 \left(N_\Omega(x)\right)_i(c_{ij}D_ju)v\,d\sigma(x) \,= \int_{\Omega} \sum_{i,j=1}^2 D_i\,(c_{ij}D_ju)v\,dx \,\, a\int_{\Omega} uv \,dx \,-\, \int_\Omega fv \, dx\ =0,
\end{equation*}
and the boundary condition in {(}\ref{ca}{)} also follows.
\eproof

  \begin{teo}\label{self}
If $\epsilon$ is small enough, the operator $C_{\epsilon}$
in $L^2(\Omega)$ defined by ($\ref{ch}$) , with domain  
given by (\ref{domicha}) is self-adjoint.
\end{teo}

  {\proof}It is clear that  $ C_{\epsilon} $ is densely defined. It is also  symmetric and lower bounded by Lemma \ref{sim}. From Lemma  
\ref{exisol} it follows that $C_{\epsilon} + a Id $
  is surjective for any $a>0$ and 
 therefore, an isomorphism from $D( C_{\epsilon})  $ to $L^2(\Omega)$. Thus
  $(C_{\epsilon} + a Id)^{-1} $ is continuous as an operator in  $L^2(\Omega)$ so  {it} has a closed graph, from which it follows that $ C_{\epsilon}$ is closed. 
  From  well-known results  (see, for instance \cite{Sc} Prop. 3.11))
  it follows that
 $C_{\epsilon}  $ is self-adjoint. 
\eproof

\vspace{3mm}

From Theorem \ref{self} it already follows that $C_{\epsilon}$ is a sectorial operator, but
 we  give  a direct proof to display the value of the constants   
appearing in Definition \ref{defset}.

\begin{teo}\label{sect}
If $\epsilon$ is small enough, the operator $C_{\epsilon}$
in $L^2(\Omega)$ defined by ($\ref{ch}$), with domain  
given by (\ref{domicha}) is sectorial and the sector in Definition \ref{defset} can be
chosen  with vertex at any
 {$b<0$} and opening angle $0< \theta <\displaystyle\frac{\pi}{2}$, and constant 
$M= cosec(\theta).$  
\end{teo}  
\proof We first observe that, as needed  when treating spectral theory, we
 work in the complexification of the relevant spaces. For simplicity however,
 we will not change notation, writing for instance,  $D(C_{\epsilon})$ for the complexification
 of the domain of the operator $C_{\epsilon}$.    Let  {$b<0$},   $0< \phi < \displaystyle\frac{\pi}{2}$
 and $\lambda = \alpha + i \beta$ a complex number in the sector
$S_{ {b}\,,\,\theta}=\left\{\,\lambda \, \big| \, \theta \leq |\,arg(\lambda -  {b} )\,| \leq \pi, \,\lambda \neq  {b}\,\right\}$.
If $u$ is  in (the complexification of)   $ D(C_{\epsilon})$, we have

\begin{align}\label{est_res}
 \| \left(C_{\epsilon} -\lambda \right) u\|_{L^2(\Omega)}  \,\| u \|_{L^2(\Omega)} & 
\geq
  \|\left\langle \left( C_{\epsilon} -\lambda \right) u  
\,,\, u \right\rangle_{L^2(\Omega)} \| \nonumber \\
 & =
\|
\left\langle \left( C_{\epsilon} -\alpha \right) u - i \beta u 
\,,\, u \right\rangle_{L^2(\Omega)} \| \nonumber  \\
 &
 = \left[ \left\langle \left( C_{\epsilon} -\alpha \right) u
  \,,\, u \right\rangle_{L^2(\Omega)}^2 +
\beta^2 
\left\langle u \,, \, u  \right\rangle_{L^2(\Omega)}^2 \right]^{\frac{1}{2}}. 
\end{align}

If $\alpha \leq  {b}$, it follows  from (\ref{est_res}) and (\ref{estCh_complex}) that
\begin{align*} 
 \| \left(C_{\epsilon} -\lambda \right) u\|_{L^2(\Omega)}  \,\| u \|_{L^2(\Omega)} & 
\geq
\left[ \left\langle \left( C_{\epsilon} - {b} \right) u
  \,,\, u \right\rangle_{L^2(\Omega)}^2 +
 \left\langle \left(  {b} -\alpha \right) u
  \,,\, u \right\rangle_{L^2(\Omega)}^2 +
\beta^2 
\left\langle u \,, \, u  \right\rangle_{L^2(\Omega)}^2 \right]^{\frac{1}{2}}
\nonumber \\   
 & \geq 
\left[\, \left(  {b} -\alpha \right)^2 +
\beta^{\,2} \, \right]^{\frac{1}{2}}
\left\langle  u
  \,,\, u \right\rangle_{L^2(\Omega)} .
\end{align*}
It follows that  
$ \| \left(C_{\epsilon} -\lambda \right) u\|_{L^2(\Omega)}  \geq
{|\lambda -  {b}|} \|u\|_{ {L^2(\Omega)}}$, for any    $ u\in D(C_{\epsilon})$ and thus 
\begin{equation}\label{estset_case1}
 {||}\left(C_{\epsilon} -\lambda \right)^{-1} u\|_{L^2(\Omega)}   \leq
\frac{1}{|\lambda -  {b}|}  \|u\|_{L^2(\Omega)}
\end{equation}
and the resolvent  inequality in  {Definition} \ref{defset} holds with $M=1$.

If $\alpha \geq  {b}$, then $\beta \geq \tan (\theta)$, and  it follows  
from (\ref{est_res}) that
\begin{align*} 
 \| \left(C_{\epsilon} -\lambda \right) u\|_{L^2(\Omega)}  \,\| u \|_{L^2(\Omega)} & 
\geq
\left[ 
\beta^2 
\left\langle u \,, \, u  \right\rangle_{L^2(\Omega)}^2 \right]^{\frac{1}{2}}
\nonumber \\
 & \geq 
  |\beta| 
\left\langle u \,, \, u  \right\rangle_{L^2(\Omega)}.
\end{align*}

It follows that  
$ \| \left(C_{\epsilon} -\lambda \right) u\|_{L^2(\Omega)}  \geq
{|\lambda -  {b}|}  \frac{|\beta|}{|\lambda -  {b}|  } \|u\|_{ {L^2(\Omega)}} 
 \geq  {|\lambda -  {b}|}  \sin{\theta  } \|u\|_{ {L^2(\Omega)}} $, for any    $ u\in D(C_{\epsilon})$ and thus 
\begin{equation}\label{estset_case2}
 {||}\left(C_{\epsilon} -\lambda \right)^{-1} u\|_{L^2(\Omega)}  \leq
\frac{cosec(\theta)}{|\lambda -  {b}|}   { \|u\|_{L^2(\Omega)} } 
\end{equation}
and the resolvent  inequality in  {Definition} \ref{defset} holds with $M= cosec(\theta)$.

From  (\ref{estset_case1})  and   {(\ref{estset_case2})}, we conclude that 
$C_{\epsilon}$ is sectorial, and the sector in   {Definition} \ref{defset} can be any sector with vertex
 in  {$b<0$}  and opening angle $0< \theta < \frac{\pi}{2}$, with constant 
$M= cosec(\theta).$   \eproof

\begin{teo}\label{sectorial} 
If  $\epsilon> 0$ is  sufficiently small and  $h_\epsilon \in \dif^1(\Omega)$,
then the operator $A_{\epsilon}= \left(\,- h_\epsilon^*\Delta_{\Omega_{\epsilon}}h_\epsilon^{*-1} +aI\,\right)$ defined by  (\ref{opeper}) and (\ref{dominio})
 is sectorial.
\end{teo}
 \proof We write  the operator as in (\ref{decomposition})  
$A_{\epsilon}=- h_\epsilon^*\Delta_{\Omega_{\epsilon}}h_\epsilon^{*-1} +aI = C_{\epsilon} + aI + L_{\epsilon}\,,$
and observe that, if $\epsilon$ is small enough, by  Theorem \ref{sect} the operator  $C_{\epsilon} + aI$ is sectorial
 with vertex in the origin  and opening angle $0< \theta < \frac{\pi}{2}$, with constant 
$M= cosec(\theta)$. For definiteness, we may  take 
 $\theta =   \displaystyle\frac{\pi}{6}$, $M=2$.    

Furthermore, by (\ref{estCh}), we obtain
 \begin{align*}
 \|\left(C_{\epsilon} + a \right)u\|_{L^2(\Omega)} \|u\|_{L^2(\Omega)} & \geq 
 \langle \left(C_{\epsilon} + a \right)u \,,\, u \rangle_{ {L^2(\Omega)}} \\
 & \geq  \displaystyle \frac{1}{4}\int_\Omega \,|\,\nabla u\,|^{\,2} \,dx
 + a \int_\Omega {  \,|\,u\,|^{\,2} \,dx }\\
& \geq  K ||u||^2_{H^1(\Omega)},
\end{align*}
where $K = \min \left\{\displaystyle\frac{1}{4}, a \right\}$ is a positive constant.
Thus $ {\|}\left(C_{\epsilon} + a \right)u\|_{L^2(\Omega)} \geq  K ||u||_{H^1(\Omega)}$.

Now, $D(C_{\epsilon} + aI + L_{\epsilon})=D(C_{\epsilon} + aI)$ and for any
$ u \in D(C_{\epsilon} + aI)$, we have

\begin{align*}
& \big|\big|\,\left(C_{\epsilon}u + au + L_{\epsilon}u\right) - \left(C_{\epsilon}u + au\right)\,\big|\big|_{L^2(\Omega)}  =  \big|\big|\,L_{\epsilon}u\,\big|\big|_{L^2(\Omega)}   \\
& \leq \bigg|\bigg|\,\displaystyle\frac{\partial b_{12}}{\partial x_2}\bigg|\bigg|_\infty\bigg|\bigg|\displaystyle\frac{\partial u}{\partial x_1} \,\bigg|\bigg|_{L^2(\Omega)} + \bigg|\bigg|\,b_{12}\displaystyle\frac{\partial b_{12}}{\partial x_2}\bigg|\bigg|_\infty\bigg|\bigg|\displaystyle\frac{\partial u}{\partial x_2}\,\bigg|\bigg|_{L^2(\Omega)}  \\
& \leq  k_1\left(\,\bigg|\bigg|\,\displaystyle\frac{\partial b_{12}}{\partial x_2}\bigg|\bigg|_\infty + \bigg|\bigg|\,b_{12}\displaystyle\frac{\partial b_{12}}{\partial x_2}\bigg|\bigg|_\infty\right)\left|\left|\,u\,\right|\right|_{H^2(\Omega)}  \\
& \leq  k_2\left(\,\bigg|\bigg|\,\displaystyle\frac{\partial b_{12}}{\partial x_2}\bigg|\bigg|_\infty + \bigg|\bigg|\,b_{12}\displaystyle\frac{\partial b_{12}}{\partial x_2}\bigg|\bigg|_\infty\right)\left|\left|\,(C_{\epsilon} + aI)u\,\right|\right|_{L^2(\Omega)} \\
& \leq {\eta}(\epsilon)
\left|\left|\,(C_{\epsilon} + aI)u\,\right|\right|_{L^2(\Omega)},
\end{align*}
where  ${\eta}(\epsilon) \to 0$ as $\epsilon \to 0.$
 From Lemma \ref{sector}, we conclude that the operator 
$$ A_{\epsilon} := C_{\epsilon} + aI + L_{\epsilon}$$
is sectorial with the sector in Definition \ref{defset}  given by    
\begin{equation}\label{def_sector}
S_{0\,,\,\frac{\pi}{6}}=\left\{\,\lambda \, \big| \, \frac{\pi}{6}  \leq |\,arg(\lambda )\,| \leq \pi, \,\lambda \neq 0\,\right\} 
\end{equation}
and  
 \begin{equation}\label{res_const}
||\,(\,\lambda - A_{\epsilon}\,)^{-1}||\leq \frac{4\sqrt{5}}{|\,\lambda \,|} \,,\,\,\mbox{ for all } \,\,\lambda \in S_{0\,,\,\frac{\pi}{6}}.
\end{equation}
\eproof

\begin{rem}
 From the proof above, we can see that 
 the sector and  constant $M$ in the resolvent inequality can be chosen
 independently of $\epsilon$, for instance, as  given in
   (\ref{def_sector}) and (\ref{res_const}).                                             
\end{rem}

\subsection{Sectoriality in $H^{-1}$}

We now want to extend the operator
  $A_{\epsilon} =
\left(-h_\epsilon^*\Delta_{\Omega_{\epsilon}}h_\epsilon^{*-1} + aI\right)$  to an operator $\widetilde{A}_{\epsilon}$
 in   $H^{-1}(\Omega)$ with  $D(\widetilde{A}_{\epsilon})=H^1(\Omega)$ and show that this extension is also sectorial for
 $\epsilon$ small enough.  

\par If   $u \in D(A_{\epsilon})=\left\{\,u \in H^2(\Omega)\,\,\bigg|\,\,h_\epsilon^*\displaystyle\frac{\partial}{\partial N_{\Omega_{\epsilon}}}h_\epsilon^{*-1}u=0\,\right\}$, $\psi \in H^1(\Omega)$, and
 $v = u \circ h_\epsilon^{-1}$, we obtain, integrating by parts

\begin{align} \label{weak_form}
\left\langle A_{\epsilon}u\,,\,\psi\right\rangle_{{-1,1}} &= - \displaystyle\int_\Omega (h_\epsilon^*\Delta_{\Omega_{\epsilon}}h_\epsilon^{*-1}u)(x)\,\psi(x)\,dx+ a\displaystyle\int_\Omega u(x)\psi(x)\,dx  \nonumber \\
&= - \displaystyle\int_\Omega \Delta_{\Omega_{\epsilon}} (u \circ h_\epsilon^{-1})(h_\epsilon(x))\,\psi(x)\,dx + a\displaystyle\int_\Omega u(x)\psi(x)\,dx  \nonumber \\
&= - \displaystyle\int_{\Omega_{\epsilon}} \Delta_{\Omega_{\epsilon}} v(y)\psi(h_\epsilon^{-1}(y)) \displaystyle\frac{1}{|Jh_\epsilon(h_\epsilon^{-1}(y))|}dy + a\displaystyle\int_{\Omega_{\epsilon}} u(h_\epsilon^{-1}(y))\psi(h_\epsilon^{-1}(y))\displaystyle\frac{1}{|Jh_\epsilon(h_\epsilon^{-1}(y))|}dy \nonumber \\
&=  -\displaystyle\int_{\partial \Omega_{\epsilon}} \displaystyle\frac{\partial v}{\partial N_{\Omega_{\epsilon}}}(y)\, \psi(h_\epsilon^{-1}(y)) \displaystyle\frac{1}{|\,Jh_\epsilon(h_\epsilon^{-1}(y))\,|}\,d\sigma (y)  \nonumber \\
& + \displaystyle\int_{\Omega_{\epsilon}} \nabla_{\Omega_{\epsilon}} v(y)\cdot \nabla_{\Omega_{\epsilon}}\psi(h_\epsilon^{-1}(y)) \displaystyle\frac{1}{|\,Jh_\epsilon(h_\epsilon^{-1}(y))\,|}\,dy  \nonumber \\
& \,+ a\displaystyle\int_{\Omega_{\epsilon}} u(h_\epsilon^{-1}(y))\psi(h_\epsilon^{-1}(y))\,\displaystyle\frac{1}{|\,Jh_\epsilon(h_\epsilon^{-1}(y))\,|}\,dy\, \nonumber \\
&=  
\displaystyle\int_{\Omega_{\epsilon}} \nabla_{\Omega_{\epsilon}} v(y)\cdot \nabla_{\Omega_{\epsilon}}\psi(h_\epsilon^{-1}(y)) \displaystyle\frac{1}{|\,Jh_\epsilon(h_\epsilon^{-1}(y))\,|}\,dy \nonumber\\
& \, + a\displaystyle\int_{\Omega_{\epsilon}} u(h_\epsilon^{-1}(y))\psi(h_\epsilon^{-1}(y))\,\displaystyle\frac{1}{|\,Jh_\epsilon(h_\epsilon^{-1}(y))\,|}\,dy\nonumber\\
&=  
 \displaystyle\int_{\Omega} \nabla_{\Omega_{\epsilon}}v(h_\epsilon(x))\cdot \nabla_{\Omega_{\epsilon}}\psi \circ h_\epsilon^{-1}\displaystyle\frac{1}{|Jh_\epsilon \circ h_\epsilon^{-1}|}(h_\epsilon(x)) |Jh_\epsilon(x)|dx \nonumber
 + a\displaystyle\int_{\Omega} u(x)\psi(x)dx\nonumber \\
&=  
\displaystyle\int_{\Omega} (h_\epsilon^*\nabla_{\Omega_{\epsilon}}h_\epsilon^{*-1}u)(x) \cdot\Bigg(h_\epsilon^* \nabla_{\Omega_{\epsilon}} h_\epsilon^{*-1} \displaystyle\frac{\psi}{|\,Jh_\epsilon\,|}\Bigg)(x)\,|\,Jh_\epsilon(x)\,|\,dx\,
+ a\displaystyle\int_{\Omega} u(x)\psi(x)\,dx .
\end{align}
 Since  (\ref{weak_form}) is well defined for  $u \in H^1(\Omega)$, we may define   an extension $\widetilde{A}_{\epsilon}$ of $A_{\epsilon}$, with values in $H^{-1}(\Omega)$ by this expression.
For simplicity, we still denote this extension by $A_{\epsilon},$ whenever there is no danger of confusion. We now show that this extension is a sectorial operator.
 
\begin{teo}\label{sectorial_weak} If $\epsilon >0$ is sufficiently small, the operator  $A_{\epsilon}$ defined by   (\ref{weak_form}), with domain
 $H^1(\Omega)$ is sectorial.
\end{teo}
\proof We  will apply  Lemma \ref{sector} again. If $\epsilon = 0$ we obtain, 
 from (\ref{weak_form}) 

\[
\left\langle A_{\epsilon}u\,,\,u\right\rangle_{ {-1,1}} = \left\langle A u\,,\,u\right\rangle_{ {-1,1}}   \geq
\displaystyle\int_{\Omega} | \nabla_{\Omega }|^2 (x)  \,dx\,
 + a\displaystyle\int_{\Omega} u^2(x) \,dx. 
\]

From the Lax-Milgram Theorem, it follows that, for any 
$\psi \in H^{-1} (\Omega)$, there exists $u \in H^{1}(\Omega)$ with 
 $Au = \psi$, so $Au$ is well defined as an operator in
 $H^{-1} (\Omega)$,  with  $ D(A) = H^{1} (\Omega)$. Using the same arguments of Theorems \ref{self} and \ref{sect}, we prove that 
 $A$ is self-adjoint and  sectorial with

$$
\|\left(A -\lambda \right)^{-1} \|_{L^2(\Omega)} \leq
\frac{cosec(\theta)}{|\lambda -  {b}|}, \textrm{ for  } \lambda \in 
S_{ {b}\,,\,\theta}=\left\{\lambda \, \big| \, \theta \leq |arg(\lambda -  {b} )| \leq \pi, \,\lambda \neq  {b}\right\}\, {and \,\theta \in \left(0,\displaystyle\frac{\pi}{2}\right)}. 
$$

 We now want to show that the family of operators
 $\big\{A_{\epsilon}\big\} $ satisfy the hypotheses of  Lemma \ref{sector}
for $\epsilon$ small.  Clearly $D\big(A_{\epsilon}\big) \supset D\big(A\big)$ for any $\epsilon \geq 0$.  We now prove that there exists a positive function    $\tau(\epsilon)$    such that
$$
\big|\big|\,\big(A_{\epsilon}  - A  \big)u\,\big|\big|_{H^{-1}(\Omega)} \leq  {\tau}(\epsilon)\big|\big|\, A  \,u\,\big|\big|_{H^{-1}(\Omega)} \,,
$$ 
for all  $u \in D\big( A  \big)$, with  $\displaystyle\lim_{\epsilon \to 0}  {\tau}(\epsilon)=0,$    which is equivalent to 
$$
\big|\,\big\langle\,\big(A_{\epsilon}  -A \big)u\,,\,\psi\,\big\rangle_{-1,1}\big| \leq \tau(\epsilon)||\,u\,||_{H^{1}(\Omega)}||\,\psi\,||_{H^{1}(\Omega)}\,, 
$$
for all  $u, \psi \in H^1(\Omega)$, with  $\displaystyle\lim_{\epsilon \to 0} \tau(\epsilon)=0.$ 
 In fact, we have, for $\epsilon >0$.

\begin{eqnarray}\label{int1}
\left|\big\langle \,\big(A_{\epsilon}  -A \,\big)u\,,\,\psi\,\big\rangle_{ {-1,1}}\right| 
&=& \bigg|\displaystyle\int_{\Omega} \big(\,h_\epsilon^*\nabla_{\Omega_{\epsilon}}h_\epsilon^{*-1}u\,\big)(x) \cdot\left(\,h_\epsilon^* \nabla_{\Omega_{\epsilon}}h_\epsilon^{*-1} \displaystyle\frac{\psi}{|\,Jh_\epsilon\,|}\,\right)(x)\,|\,Jh_\epsilon(x)\,|\,dx\nonumber\\
&&\,-\,\displaystyle\int_{\Omega}\nabla_\Omega u(x)\cdot\nabla_\Omega\psi(x)\,dx \,\bigg | \nonumber \\
&=& \bigg|\displaystyle\int_{\Omega} \left(h_\epsilon^*\nabla_{\Omega_{\epsilon}}h_\epsilon^{*-1}u - \nabla_\Omega u\right)(x) \cdot\left(h_\epsilon^* \nabla_{\Omega_{\epsilon}}h_\epsilon^{*-1} \displaystyle\frac{\psi}{|Jh_\epsilon|}\right)(x)\,|\,Jh_\epsilon(x)\,| \nonumber\\
&&\,+\, \nabla_\Omega u(x) \cdot\left(h_\epsilon^* \nabla_{\Omega_{\epsilon}}h_\epsilon^{*-1} \displaystyle\frac{\psi}{|Jh_\epsilon|}\right)(x)\,|Jh_\epsilon(x)| - \nabla_\Omega u(x)\cdot\nabla_\Omega\psi(x)\,dx \,\bigg| \nonumber\\  
&=& \bigg|\displaystyle\int_{\Omega} \left(h_\epsilon^*\nabla_{\Omega_{\epsilon}}h_\epsilon^{*-1}u - \nabla_\Omega u\right)(x) \cdot\left(h_\epsilon^* \nabla_{\Omega_{\epsilon}}h_\epsilon^{*-1} \displaystyle\frac{\psi}{|\,Jh_\epsilon\,|}\right)(x)\,|\,Jh_\epsilon(x)\,| \nonumber\\
&&\,+\,\nabla_\Omega u(x) \cdot\left[h_\epsilon^* \nabla_{\Omega_{\epsilon}}h_\epsilon^{*-1} \displaystyle\frac{\psi}{|\,Jh_\epsilon\,|}(x) - \nabla_\Omega\left(\displaystyle\frac{\psi}{|\,Jh_\epsilon\,|}\right)(x)\right]\,|\,Jh_\epsilon(x)\,| \nonumber\\
&&\,+\,\nabla_\Omega u(x) \cdot\nabla_\Omega\left(\displaystyle\frac{\psi}{|\,Jh_\epsilon\,|}\right) {(x)}|\,Jh_\epsilon(x)\,| - \nabla_\Omega u(x)\cdot\nabla_\Omega\psi(x)\,dx \,\bigg| \nonumber\\
 &\leq &  
\displaystyle\int_\Omega\left|(h_\epsilon^*\nabla_{\Omega_{\epsilon}}h_\epsilon^{*-1}u-\nabla_\Omega u)(x)\cdot
\left(h_\epsilon^*\nabla_{\Omega_{\epsilon}}h_\epsilon^{*-1}\frac{\psi}{|Jh_\epsilon|}\right)(x)\right||Jh_\epsilon {(x)}|dx \label{integral1}\\
&&\,+\, \displaystyle\int_\Omega\left|\,\nabla_\Omega u(x)\cdot\big(h_\epsilon^*\nabla_{\Omega_{\epsilon}}h_\epsilon^{*-1}-
\nabla_\Omega\big)\left(\frac{\psi}{|\,Jh_\epsilon\,|}\right)(x)\,\right|\,|Jh_\epsilon(x)|\,dx \label{integral2}\\
&&\,+\, \displaystyle\int_\Omega \psi {(x)}\nabla_\Omega u(x) \cdot\nabla_\Omega\left(\frac{1}{|\,Jh_\epsilon\,|}\right)(x)\,|\,Jh_\epsilon(x)\,|\,dx\,. \label{integral3} 
\end{eqnarray}

\par To estimate  (\ref{integral1}), (\ref{integral2}) and  (\ref{integral3}),
 we use the expression   for   $h_\epsilon^*\nabla_{\Omega_{\epsilon}}h_\epsilon^{*-1}$  in terms of the coefficients 
 of $[{h_\epsilon}^{-1}]^{\,T}_{\,x}$ given by  (\ref{deriv}). For  (\ref{integral1}), we have

$
\begin{array}{lll}
&&
\displaystyle\int_\Omega\left|\,\left(h_\epsilon^*\nabla_{\Omega_{\epsilon}}h_\epsilon^{*-1}u -\nabla_\Omega u\right)(x) \cdot
\left(h_\epsilon^*\nabla_{\Omega_{\epsilon}}h_\epsilon^{*-1}\frac{\psi}{|\,Jh_\epsilon\,|}\right)(x)\,\right||\,Jh_\epsilon(x)\,|\,dx \\
&\leq& \displaystyle\int_\Omega \left|\left(
\begin{array}{c}
\displaystyle\sum_{j=1}^2(b_{1j}(x)-\delta_{1j})\frac{\partial u}{\partial x_j}(x)\\
\displaystyle\sum_{j=1}^2(b_{2j}(x)-\delta_{2j})\frac{\partial u}{\partial x_j}(x)
\end{array}
\right) 
\left(
\begin{array}{c}
\displaystyle\sum_{j=1}^2 b_{1j}(x)\frac{\partial }{\partial x_j}\left(\frac{\psi}{|Jh_\epsilon|}\right)(x)
\\
\displaystyle\sum_{j=1}^2 b_{2j}(x)\frac{\partial }{\partial x_j}\left(\frac{\psi}{|Jh_\epsilon|}\right)(x)
\end{array}
\right)\right| \big|1 + \epsilon sen(x_1/\epsilon^\alpha)\big|dx
\\
&\leq& \displaystyle \max_{(x_1,x_2)\,\in\, ]0,1[\times ]0,1[}\big\{\,\big|\,Jh_\epsilon(x)\,\big|\,\big\}
\int_\Omega \left[\,b^{\,2}_{12}(x)\left(\frac{\partial u}{\partial x_2}\right)^{\,2}+ 
\big(\,b_{22}(x) - 1\,\big)^{\,2}\left(\frac{\partial u}{\partial x_2}\right)^{\,2}\,\right]^\frac{1}{2} 
\\
&&\,\cdot\,\displaystyle \left\{\left[\frac{\partial}{\partial x_1}\left(\frac{\psi}{|Jh_\epsilon|}\right)(x)
+ b_{12}(x)\frac{\partial}{\partial x_2}\left(\frac{\psi}{|Jh_\epsilon|}\right)(x)\right]^2
+ b^2_{22}(x)\left[\frac{\partial}{\partial x_2}\left(\frac{\psi}{|Jh_\epsilon|}\right)(x)\right]^2\right\}^\frac{1}{2}dx\\
&\leq& \displaystyle \max_{(x_1,x_2)\,\in \,]0,1[\times ]0,1[}\big\{\,\big|\,Jh_\epsilon(x)\,\big|\,\big\}\underbrace{
\left\{\,\displaystyle\int_\Omega \left[\,b^{\,2}_{12}(x) + (b_{22}(x) - 1)^{\,2}\,\right]\left(\,
\frac{\partial u}{\partial x_2}\,\right)^{\,2}dx\,\right\}^\frac{1}{2}}_{(I)}\\
&&\,\cdot\,\underbrace{\left\{\displaystyle \int_\Omega \left[\frac{\partial}{\partial x_1}\left( \frac{\psi}{|Jh_\epsilon|}\right)(x) + b_{12}(x)\frac{\partial}{\partial x_2}\left(\frac{\psi}{|Jh_\epsilon|}\right)(x)\right]^2 + 
b^2_{22}(x) \left[\frac{\partial}{\partial x_2}\left(\frac{\psi}{|Jh_\epsilon|}\right)(x)\right]^2dx \right\}^\frac{1}{2}}_{(II)}.
\end{array}
$
\par For the integral  $(I)$, we have 
$$
\begin{array}{ll}
\left\{\displaystyle\int_\Omega \left[b^{2}_{12}(x) + \big(b_{22}(x) - 1\big)^{2}\right]\left(\frac{\partial u}{\partial x_2}\right)^{2}dx\right\}^\frac{1}{2}
&=\left\{\displaystyle \int_\Omega \frac{x_2^{2}\epsilon^{2-2\alpha}cos^{2}(x_1/\epsilon^{\alpha})+ \epsilon^{2} sen^{2}(x_1/\epsilon^{\alpha})}{[1+ \epsilon sen(x_1/\epsilon^{\alpha})]^{2}}\displaystyle\left(\frac{\partial u}{\partial x_2}\right)^{2}dx\right\}^\frac{1}{2}\\
&\leq \left\{\displaystyle \int_\Omega \frac{x_2^{\,2}\epsilon^{\,2-2\alpha} + \epsilon^{\,2} }{[\,1 + \epsilon \,sen(x_1/\epsilon^{\,\alpha})\,]^{\,2}}\displaystyle\left(\frac{\partial u}{\partial x_2}\right)^{\,2}dx\,\right\}^\frac{1}{2} \\
& \leq   
\big(\,K_1(\epsilon)\,\big)^\frac{1}{2} \left\{\,\displaystyle\int_\Omega \nabla_\Omega u(x)\cdot\nabla_\Omega u(x)\,dx\,\right\}^\frac{1}{2}, 
\end{array}
$$
 where  $K_1(\epsilon):=\displaystyle\frac{x_2^{\,2}\epsilon^{\,2-2\alpha} + \epsilon^{\,2} }{[\,1 + \epsilon \,sen(x_1/\epsilon^{\,\alpha})\,]^{\,2}} \to 0$ as $h_\epsilon \to i_\Omega$  in  $C^1(\Omega)$ . 

\par To estimate the  integral $(II)$, we first observe that  

$$
\displaystyle\frac{\partial}{\partial x_1}\left(\,\displaystyle\frac{1}{|\,Jh_\epsilon(x)\,|}\,\right) =
-\displaystyle\frac{\epsilon^{\,1-\alpha}\,cos(x_1/\epsilon^{\,\alpha})}{[\,1+ \epsilon \, sen(x_1/\epsilon^{\,\alpha})\,]^{\,2}} \ \  \textrm{ and } \ \ 
\displaystyle\frac{\partial}{\partial x_2}\left(\,\displaystyle\frac{1}{|\,Jh_\epsilon(x)\,|}\,\right) = 0\,.
$$

\par Therefore, we have the following estimate for  $(II)$ : 
$$
\begin{array}{lll}
&&\left\{\,\displaystyle \int_\Omega \left[\,\frac{\partial}{\partial x_1}\left(\,\frac{\psi}{|\,Jh_\epsilon\,|}\,\right)(x)+ b_{12}(x)\frac{\partial}{\partial x_2}\left(\,\frac{\psi}{|\,Jh_\epsilon\,|}\,\right)(x)\,\right]^{\,2} + 
b^{\,2}_{22}(x) \left[\,\frac{\partial}{\partial x_2}\left(\,\frac{\psi}{|\,Jh_\epsilon\,|}\,\right)(x)\,\right]^{\,2}dx\,\right\}^\frac{1}{2}\\
&\leq& \left\{\,\displaystyle\int_\Omega 2\left[\,\frac{\partial}{\partial x_1}\left(\,\frac{\psi}{|\,Jh_\epsilon\,|}\,\right)(x)\,\right]^{\,2} + \big(\,2b^{\,2}_{12}(x) + b^{\,2}_{22}(x)\,\big)\left[\,\frac{\partial}{\partial x_2}\left(\,\frac{\psi}{|\,Jh_\epsilon\,|}\,\right)(x)
\,\right]^{\,2}dx\,\right\}^\frac{1}{2}\\
&=&\left\{\displaystyle \int_\Omega 2\left[\psi(x) \frac{\partial}{\partial x_1}
\left(\frac{1}{|Jh_\epsilon(x)|}\right) + \frac{1}{|Jh_\epsilon(x)|}\frac{\partial \psi}{\partial x_1}(x)\right]^{2} + 
\big(2b^{2}_{12}(x) + b^{2}_{22}(x)\big)\left[\frac{1}{|Jh_\epsilon(x)|}\frac{\partial \psi}{\partial x_2}(x)\right]^{2}dx\right\}^\frac{1}{2}\\
&\leq& \left\{\, \displaystyle \int_\Omega 4\psi^{\,2}\left[\,\frac{\partial}{\partial x_1}
\left(\frac{1}{|\,Jh_\epsilon(x)\,|}\right)\,\right]^{\,2} + \frac{4}{|\,Jh_\epsilon(x)\,|^{\,2}}\left(\frac{\partial \psi}{\partial x_1}(x)\right)^{\,2} \right. \\
&+& \left. \big(\,2b^{\,2}_{12}(x) + b^{\,2}_{22}(x)\,\big)\displaystyle\frac{1}{|\,Jh_\epsilon(x)\,|^{\,2}}\left(\displaystyle\frac{\partial \psi}{\partial x_2}(x)\right)^{\,2}dx\right\}^\frac{1}{2}\\
&\leq& \underbrace{\left\{\,\displaystyle\int_\Omega\left[4\left(\frac{\partial}{\partial x_1}
\left(\,\frac{1}{|\,Jh_\epsilon(x)\,|}\,\right)\right)^{\,2}\right] \psi^{\,2}(x)dx\right\}^\frac{1}{2}}_{(III)} \\
&+&
\underbrace{\left\{\,\displaystyle \int_\Omega \frac{4}{|\,Jh_\epsilon(x)\,|^{\,2}}\left(\frac{\partial \psi}{\partial x_1}(x)\right)^{\,2} + \frac{2b^{\,2}_{12}(x) + b^{\,2}_{22}(x)}{|\,Jh_\epsilon(x)\,|^{\,2}} \left(\frac{\partial \psi}{\partial x_2}(x)\right)^{\,2}dx\,\right\}^\frac{1}{2}}_{(IV)}.
\end{array}
$$

\par To estimate  $(III)$,  we first observe that 
$$
K_2(\epsilon):\,=\,4\displaystyle\left[\frac{\partial}{\partial x_1}
\left(\frac{1}{|\,Jh_\epsilon(x)\,|}\right)\right]^{\,2} 
\,=\,\displaystyle\frac{4\epsilon^{\,2-2\alpha} \,cos^{\,2}(x_1/\epsilon^{\,\alpha})}{[\,1+ \epsilon \,sen(x_1/\epsilon^{\,\alpha})\,]^{\,4}} \to 0 $$
as  $h_\epsilon \to i_\Omega$  in  $C^1(\Omega)$.  Therefore
$$
\left\{ \displaystyle \int_\Omega \left[\,4\left(\,\frac{\partial}{\partial x_1}
\left(\frac{1}{|\,Jh_\epsilon(x)\,|}\,\right)\right)^{\,2} \right] \psi^{\,2}(x)dx\right\}^\frac{1}{2}
\leq \big(K_2(\epsilon)\big)^\frac{1}{2}\left[\,\displaystyle \int_\Omega \psi^{\,2}dx\,\right]^\frac{1}{2},
$$
with  $\big(K_2(\epsilon)\big)^\frac{1}{2} \to 0$.

\par To estimate  $(IV)$, observe that 
$$\displaystyle\frac{4}{|\,Jh_\epsilon(x)\,|^{\,2}}= \frac{4}{[\,1+ \epsilon \,sen(x_1/\epsilon^{\,\alpha})\,]^{\,2}}\,\,\,\,\,\mbox{ and }\,\,\,\,\,\displaystyle\frac{ {2}b^{\,2}_{12}(x) + b^{\,2}_{22}(x)}{|\,Jh_\epsilon(x)\,|^{\,2}}=\frac{ {2}x_2^{\,2}\epsilon^{\,2-2\alpha}cos^{\,2}(x_1/\epsilon^{\,\alpha})+  {1}}{[\,1+ \epsilon \,sen(x_1/\epsilon^{\,\alpha})\,]^{\,4}}$$ are bounded for  $0 < \alpha <1$ and $\epsilon >0$ 
sufficiently small. Thus
$$
\left\{ \displaystyle \int_\Omega \frac{4}{|Jh_\epsilon(x)|^2}\left(\frac{\partial \psi}{\partial x_1}(x)\right)^2 + \frac{2b^2_{12}(x) + b^2_{22}(x)}{|Jh_\epsilon(x)|^2} \left(\frac{\partial \psi}{\partial x_2}(x)\right)^2dx\right\}^\frac{1}{2} \leq K_3\left\{\displaystyle \int_\Omega \nabla_\Omega \psi(x) \cdot \nabla_\Omega \psi(x)\,dx\right\}^\frac{1}{2}, 
$$
where  $K_3$ is a positive constant. We then have the following estimate for
 (\ref{integral1}):

$$
\begin{array}{lll}
&&\displaystyle\int_\Omega\left|\,\big(h_\epsilon^*\nabla_{\Omega_{\epsilon}}h_\epsilon^{*-1}u -\nabla_\Omega u\big)(x)\cdot \left(\
h_\epsilon^*\nabla_{\Omega_{\epsilon}}h_\epsilon^{*-1}\frac{\psi}{|\,Jh_\epsilon\,|}\right)(x)\,\right|\,\big|\,Jh_\epsilon(x)\,\big|\,dx \\
&\leq& \displaystyle \max_{(x_1,x_2)\,\in\, ]0,1[\times ]0,1[}\big\{\,\big|\,Jh_\epsilon(x)\,\big|\,\big\}
\,[\,K_1(\epsilon)K_2(\epsilon)\,]^\frac{1}{2}\left\{\displaystyle\int_\Omega \psi^{\,2}\,dx\right\}^\frac{1}{2}\left\{\,\displaystyle\int_\Omega \nabla_\Omega u \cdot \nabla_\Omega u\,dx\,\right\}^\frac{1}{2} \\
&+& \displaystyle \max_{(x_1,x_2)\,\in\, ]0,1[\times ]0,1[}\big\{\,\big|\,Jh_\epsilon(x)\,\big|\,\big\}
[\,K_1(\epsilon)\,]^\frac{1}{2}K_3\,\left\{\,\displaystyle\int_\Omega \nabla_\Omega \psi \cdot \nabla_\Omega \psi\,dx\,\right\}^{\,\frac{1}{2}}\left\{\displaystyle\int_\Omega \nabla_\Omega u\cdot 
\nabla_\Omega u \, dx\right\}^{\,\frac{1}{2}}.
\end{array}
$$

Taking
$$C_0(\epsilon):=\displaystyle \max_{(x_1,x_2)\,\in\,]0,1[\times ]0,1[}\big\{\big|Jh_\epsilon(x)\big|\big\}
[K_1(\epsilon)K_2(\epsilon)]^\frac{1}{2}\,\,\mbox{ and }\,\,C^{\,'}_0(\epsilon):=\displaystyle \max_{(x_1,x_2)\,\in\,]0,1[\times ]0,1[}\big\{\big|Jh_\epsilon(x)\big|\big\}[K_1(\epsilon)]^\frac{1}{2}K_3\,,$$
we have that  $C_0(\epsilon)$ and $C^{\,'}_0(\epsilon) \to 0$
 when  $h_\epsilon \to i_\Omega$  in  $C^1(\Omega)$.
\par  We estimate (\ref{integral2}) in a similar way:
$$
\begin{array}{ll} 
&
\displaystyle\int_\Omega\left|\,\nabla_\Omega u(x)\cdot\big(h_\epsilon^*\nabla_{\Omega_{\epsilon}}h_\epsilon^{*-1}-
\nabla_\Omega\big)\left(\frac{\psi}{|\,Jh_\epsilon\,|}(x)\right)\,\right|\,|\,Jh_\epsilon(x)\,|\,dx \\
&\leq \displaystyle \max_{(x_1,x_2)\,\in\, ]0,1[\times ]0,1[}\big\{\big|Jh_\epsilon(x)\big|\big\}
\displaystyle \left\{\int_\Omega \big|\nabla_\Omega u(x)\big|^2dx\right\}^\frac{1}{2}
\underbrace{
\displaystyle \left\{\int_\Omega \sum^2_{i=1} \left[\sum^2_{j=1} \big(b_{ij}(x) - \delta_{ij}\big) \frac{\partial}{\partial x_j}\left(\frac{\psi}{|Jh_\epsilon|}(x)\right)\right]^{2}dx\right\}^\frac{1}{2}}_{(V)}.
\end{array}
$$
 For the integrand in $(V)$, we have: 
$$
\begin{array}{ll}
&\displaystyle\sum^2_{i=1}\left[\sum^2_{j=1}\big(b_{ij}(x) - \delta_{ij}\big)\frac{\partial}{\partial x_j}\left(\,\frac{\psi}{|\,Jh_\epsilon\,|}(x)\,\right)\right]^2 \\
=& \displaystyle\sum^2_{i=1} \left[\,\big(b_{i1}(x) - \delta_{i1}\big)\frac{\partial}{\partial x_1}\left(\,\frac{\psi}{|\,Jh_\epsilon\,|}(x)\,\right) + \big(\,b_{i2}(x) - \delta_{i2}\,\big)\frac{\partial}{\partial x_2}\left(\frac{\psi}{|\,Jh_\epsilon\,|}(x)\right)\,\right]^2\\ 
\leq& \displaystyle\sum^2_{i=1} \left\{\,2\big(b_{i1}(x)-\delta_{i1}\big)^{\,2} \left[\,\frac{\partial}{\partial x_1}\left(\,\frac{\psi}{|\,Jh_\epsilon\,|}\,\right)(x)\,\right]^{\,2}
+  2\,\big(b_{i2}(x) - \delta_{i2}\big)^{\,2}\left[\,\displaystyle\frac{\partial}{\partial x_2}\left(\frac{\psi}{|\,Jh_\epsilon\,|}\right)(x)\,\right]^{\,2}\,\right\} \\ 
=& \big[\,2b_{12}^{\,2}(x) +  {2}\big(\,b_{22}(x)-1\,\big)^{\,2}\,\big] \displaystyle \left[\,\frac{\partial}{\partial x_2}\left(\frac{\psi}{|\,Jh_\epsilon\,|}\right)(x)\,\right]^{\,2} \\ 
=& \big[\,2b_{12}^{\,2}(x) +  {2}\big(\,b_{22}(x) - 1\,\big)^{\,2}\,\big]\displaystyle\frac{1}{|\,Jh_\epsilon(x)\,|^2}\left(\,\frac{\partial \psi}{\partial x_2}(x)\,\right)^{\,2}.
\end{array}
$$
We then have, for the  integral $(V)$
$$
\displaystyle \left\{\int_\Omega \sum^2_{i=1} \left[\sum^2_{j=1} \big(b_{ij}(x) - \delta_{ij}\big)\frac{\partial}{\partial x_j}\left(\frac{\psi}{|Jh_\epsilon|}(x)\right)\right]^2dx\right\}^\frac{1}{2}\leq\displaystyle \left\{\int_\Omega \frac{2b_{12}^2(x) +  {2}\big(b_{22}(x) - 1\big)^2}{|Jh_\epsilon(x)|^2}\left(\frac{\partial \psi}{\partial x_2}(x)\right)^2dx\right\}^\frac{1}{2}.
$$
and 
$$
K_4(\epsilon):= \displaystyle\frac{2b_{12}^{\,2}(x) +  {2}\big(b_{22}(x) - 1\big)^{\,2}}{|\,Jh_\epsilon(x)\,|^{\,2}}
\,=\,\displaystyle\frac{2x_2^{\,2}\,\epsilon^{\,2-2\alpha}cos^{\,2}(x_1/\epsilon^{\,\alpha})+  {2}\epsilon^{\,2}sen^{\,2}(x_1/\epsilon^{\,\alpha})
}{|\,Jh_\epsilon(x)\,|^{\,2}}\  \to 0 
$$
 when $h_\epsilon \to i_\Omega$ in $C^1(\Omega)$. Thus, taking  $C_1(\epsilon):= \displaystyle \max_{(x_1,x_2)\,\in\, ]0,1[\times ]0,1[}\big\{|\,Jh_\epsilon(x)\,|\big\}
[K_4(\epsilon)]^\frac{1}{2}$, we have the following estimate  for (\ref{integral2}):
$$ 
\begin{array}{lll}
&&
\displaystyle\int_\Omega\left|\,\nabla_\Omega u(x)\cdot\big(h_\epsilon^*\nabla_{\Omega_{\epsilon}}h_\epsilon^{*-1}-
\nabla_\Omega\big)\left(\frac{\psi}{|\,Jh_\epsilon\,|}(x)\right)\,\right|
|\,Jh_\epsilon(x)\,|\,dx \\
&\leq& \displaystyle C_1(\epsilon) \left\{\,\displaystyle\int_\Omega \nabla_\Omega \psi \cdot \nabla_\Omega \psi \,dx\,\right\}^\frac{1}{2}
\left\{\,\displaystyle\int_\Omega \nabla_\Omega u\cdot \nabla_\Omega u\, dx\,\right\}^\frac{1}{2},
\end{array} 
$$
  with  $C_1(\epsilon) \to 0$  as $h_\epsilon \to i_\Omega$  in  $C^1(\Omega)$.

\par  Finally we have for the integral (\ref{integral3}):
$$
\begin{array}{lll} 
&&
\displaystyle\int_\Omega\left|\,\psi(x) \nabla_\Omega u(x)\cdot\nabla_\Omega\left(\frac{1}{\,\big|\,Jh_\epsilon(x)\,\big|\,}\right)\,\right|\,\big|\,Jh_\epsilon(x)\,\big|\,dx \\
&\leq& \displaystyle \max_{(x_1,x_2)\,\in\,]0,1[ \times ]0,1[}\big\{\,\big|\,Jh_\epsilon(x)\,\big|\,\big\} \displaystyle\int_\Omega\psi {(x)}\left[\,\sum^2_{j=1}\left(\frac{\partial u}{\partial x_j}\right)^{\,2}\,\right]^\frac{1}{2}\left[\sum^2_{j=1}\left(\frac{\partial}{\partial x_j}\left(\frac{1}{|\,Jh_\epsilon(x)\,|}\right)\right)^{\,2}\right]^\frac{1}{2}dx \\ 
&\leq&\displaystyle \max_{(x_1,x_2)\,\in\,]0,1[\times ]0,1[}\big\{\big|Jh_\epsilon(x)\big|\big\}\left\{\displaystyle\int_\Omega\psi^{2} {(x)} dx\right\}^\frac{1}{2} 
\left\{\displaystyle\int_\Omega \sum^2_{j=1}\left(\frac{\partial u}{\partial x_j}\right)^2\left[\sum^2_{j=1}\left(\frac{\partial}{\partial x_j}\left(\frac{1}{|Jh_\epsilon(x)|}\right)\right)^2\right]dx\right\}^\frac{1}{2}
\\ 
&=&
\displaystyle \max_{(x_1,x_2)\,\in\, ]\,0,1\,[\times ]\,0,1\,[}\big\{\,\big|\,Jh_\epsilon {(x)}\,\big|\,\big\} \left\{\displaystyle\int_\Omega\psi^{\,2}  {(x)} dx\right\}^\frac{1}{2}
\left\{\displaystyle\int_\Omega \sum^2_{j=1}\left(\frac{\partial u}{\partial x_j}\right)^{\,2}\left[\frac{\epsilon^{\,2-2\alpha}
cos^{\,2}(x_1/\epsilon^{\,\alpha})}{|\,Jh_\epsilon(x)\,|^{\,2}}\right]\,dx\right\}^\frac{1}{2}.
\end{array}
$$
\par Since $K_5(\epsilon):=\displaystyle\frac{\epsilon^{\,2-2\alpha}
cos^{\,2}(x_1/\epsilon^\alpha)}{|\,Jh_\epsilon(x)\,|^{\,2}} \to 0 $ as  $h_\epsilon \to i_\Omega$ in  $C^1(\Omega)$. We have the following estimate
 for the integral (\ref{integral3}):  
$$\displaystyle\int_\Omega\left|\,\psi(x) \nabla_\Omega u(x)\cdot
\nabla_\Omega\left(\frac{1}{|\,Jh_\epsilon(x)\,|}\right)\,\right|\,\big|\,Jh_\epsilon(x)\,\big|\,dx \,\leq\, C_2(\epsilon)\left\{\displaystyle\int_\Omega\psi^{\,2} \,dx\right\}^\frac{1}{2}\left\{\displaystyle\int_\Omega \nabla_\Omega u \cdot \nabla_\Omega u\,dx\right\}^\frac{1}{2},
$$
with
$C_2(\epsilon):= \displaystyle \max_{(x_1,x_2)\,\in\,]\,0,1\,[\times ]\,0,1\,[}\big\{\,\big|\,Jh_\epsilon(x)\,\big|\,\big\}\,[K_5(\epsilon)]^\frac{1}{2}
\to 0 $ as $h_\epsilon \to i_\Omega$  in  $C^1(\Omega)$ .

\par We conclude that  
$$
\big|\,\big\langle\,\big(\,A_{\epsilon} - A\,\big)u\,,\,\psi\,\big\rangle_{-1,1}\,\big| \leq C(\epsilon)||\,u\,||_{H^{1}(\Omega)}||\,\psi\,||_{H^{1}(\Omega)}\, 
$$
 with  $\displaystyle\lim_{\epsilon \to 0+ } C(\epsilon)=0$ (independently of $u$)
and therefore

\begin{align} \label{dif_estimate}
\| \left({A}_{\epsilon} - {A}  \right)u \|_{ {H^{-1}(\Omega)}}
& \leq  C(\epsilon)||\,u\,||_{H^{1}(\Omega)} \nonumber \\
& \leq \tau(\epsilon)||\,{A} u\,||_{ {H^{-1}(\Omega)}}
\end{align}
 with  $ \displaystyle\lim_{\epsilon \to 0+ } \tau(\epsilon)=0$,
 (and $\tau(\epsilon)$ does not depend on $u$).
 Therefore, the result follows from Lemma \ref{sector}. \eproof  

\begin{rem}\label{rem_sector_weak}
From the above proof, it also follows that the sector and the constant $M$ in the resolvent inequality are the same as the ones for ${A}_{\epsilon}$,
 and can be chosen as 
 in  (\ref{def_sector}) and (\ref{res_const}).                                   \end{rem}


\section{The abstract problem in a scale of Banach spaces} 
\label{abstract}
 Our goal in this section is to pose the  problem  (\ref{nonlinBVP})  in a convenient abstract setting.
We proved in Theorem \ref{sectorial} 
  that, if $\epsilon$ is small,  the operator $A_{\epsilon}$  in $L^2(\Omega)$ defined by
 (\ref{opeper}) with domain given in (\ref{dominio}) is sectorial as well as its extension 
 $\widetilde{A}_{\epsilon}$  to $H^{-1}(\Omega)$.  
It is then well-known that the domains $X_{\epsilon}^{\alpha}$ (resp. $\widetilde{X}_\epsilon^{\alpha})$,
 $ \alpha \geq 0$  of the fractional powers of  $A_{\epsilon}$  (resp. $\widetilde{A}_{\epsilon}$) 
are Banach spaces,  $X_{\epsilon}^0 = L^2(\Omega)$, (resp. $\widetilde{X}_\epsilon^0 =  H^{-1} (\Omega)$),
 $X_{\epsilon}^1 = D(A_{\epsilon})$,  (resp. $\widetilde{X}_\epsilon^1 =  D(\widetilde{A}_{\epsilon})  $),
 $X_{\epsilon}^\alpha$,  ( $\widetilde{X}_\epsilon^\alpha$)   is compactly embedded in
  $ X_{\epsilon}^\beta$,   ($  \widetilde{X}_\epsilon^\beta$)  when  $\alpha > \beta \geq 0$,
and    $X_{\epsilon}^{\alpha} = H^{2\alpha}$, when $2 \alpha$ is an integer number.


 Since $ X_{\epsilon}^{\frac{1}{2}}= \widetilde{X}_{\epsilon}^1 $, it follows easily that
   $ X_{\epsilon}^{\alpha - \frac{1}{2}  }= \widetilde{X}_{\epsilon}^\alpha  $, for $ \frac{1}{2} 
  \leq \alpha \leq 1 $  and, by an abuse of notation,
 we will still write  $X_{\epsilon}^{\alpha - \frac{1}{2} } $ instead of
  $\widetilde{X}_\epsilon^{\alpha}$,   for  $  0  \leq 
\alpha \leq    \frac{1}{2} $  so we may denote by 
  $ \{ X_{\epsilon}^{\alpha}, \, \,     -\frac{1}{2} \leq  \alpha \leq 1 \} =
 \{ X_{\epsilon}^{\alpha}, \, \,     0 \leq   \alpha \leq 1 \} \cup
 \{ \widetilde{X}_\epsilon^{\alpha}, \, \,   0 \leq   \alpha \leq 1 \},$
 the whole family of fractional power spaces. 
We will denote  simply by  $X^{\alpha}$ the fractional power spaces associated to the unperturbed operator $A$.

 For any  $  \displaystyle -\frac{1}{2} \leq \beta \leq 0$, we may now
define an operator in  these spaces as the restriction of
$\widetilde{A}_{\epsilon}$. We then have the following result
 \begin{teo}\label{sec_scale} For any   $-\frac{1}{2} \leq \beta \leq 0$ and $\epsilon$ sufficiently small, the operator
 $(A_{\epsilon})_{\beta}$ in $ X_{\epsilon}^{\beta}$, obtained by restricting
 $ \widetilde{A}_{\epsilon}$, with domain  $ X_{\epsilon}^{\beta+1}$ is a sectorial operator.
\end{teo}
\proof
 Writing  $\beta =  -\frac{1}{2} + \delta $, for some $ 0 \leq \delta \leq
\frac{1}{2}$, we have $(A_{\epsilon})_{\beta} =
   \widetilde{A}_{\epsilon}^{-\delta} \widetilde{A}_{\epsilon}  \widetilde{A}_{\epsilon}^{\delta}. $ Since $\widetilde{A}_{\epsilon}^{\delta}$ is an isometry from
$X_{\epsilon}^{\beta}$ to $X_{\epsilon}^{ {-\frac{1}{2}}} = H^{-1} (\Omega),$ the result follows easily.  
 {\eproof}

We now show that  the   scale
  of Banach spaces $ \{X_{\epsilon}^{\alpha}, \ -\frac{1}{2}
\leq \alpha \leq \frac{1}{2} \} = \{ \widetilde{X}_{\epsilon}^{\alpha}, \ 0
\leq \alpha \leq 1\} $, do
 not change when   $\epsilon >0$ varies.  More precisely

\begin{teo}\label{uniform_scale}
For any $ -\frac{1}{2} \leq \alpha \leq \frac{1}{2}$, let $\|u \|_{\epsilon,\alpha}  $ denote the norm in $X_{\epsilon}^{\alpha}=  \widetilde{X}_{\epsilon}^{\alpha + \frac{1}{2}} $ and 
 $\|u \|_{\alpha}$  the norm in $X^{\alpha }=  \widetilde{X}^{\alpha + \frac{1}{2}} $. Then, we have 
 $\|u \|_{\epsilon,\alpha}  \leq K_1(\epsilon) \|u \|_{\alpha}
 \leq K_2(\epsilon)  \|u \|_{\epsilon,\alpha},  $
 with   $K_1(\epsilon), K_2(\epsilon) \to 1 $  as $\epsilon \to 0$, uniformly
 in $\alpha$. In particular, 
$X_{\epsilon}^{\alpha}= X^{\alpha} $,  with equivalent norms, uniformly 
 in $\epsilon$.  
\end{teo} 
\proof
 As observed above, using the same arguments of Theorems \ref{self} and \ref{sect}, we prove that  $\widetilde{A}_{\epsilon}$
  is self-adjoint with respect to the usual inner product in $H^1({\Omega})$.
  Similarly, one can prove that  $\widetilde{A}_{\epsilon}$ is self-adjoint with respect to the inner product  with ``weight'' $|Jh_{\epsilon}|.$ It then follows, from well-known results,  that the fractional powers of order $\alpha$ of these operators coincide isometricaly  with  the interpolation spaces
    $[D(\widetilde{A}), H^{-1}({\Omega})]_{\alpha}$,
  $[D(\widetilde{A}_{\epsilon}), H^{-1}({\Omega})]_{\alpha}$ respectively (see, for instance, Theorem 16.1 of \cite{Yagi}).  Let $I: X \to Y$ denote the inclusion operator from $X\subset Y$ to $Y$ and
  %
 by  $\|I \|_{\mathcal{L}(X,Y)}$ its norm as a  linear operator.
   From Theorem  {1.15} of \cite{Yagi} we then have, for any
   $u \in D(\widetilde{A}_{\epsilon})$

\begin{eqnarray} \label{estimate_L}
  \|I \|_{\mathcal{L}([D(\widetilde{A}), H^{-1}({\Omega})]_{\alpha} \,,\,
  [D(\widetilde{A}_{\epsilon}), H^{-1}({\Omega})]_{\alpha} ) }
  & \leq &
  \|I \|_{\mathcal{L}(H^{-1}(\Omega)\,,\,H^{-1}(\Omega)) }^{(1-\alpha)} 
  \|I \|_{\mathcal{L}(D(\widetilde{A}) \,,\,  D(\widetilde{A}_{\epsilon})) }^\alpha  \nonumber \\
  & \leq & \|I \|_{\mathcal{L}(D(\widetilde{A})\,,\,   D(\widetilde{A}_{\epsilon}  ) )}^\alpha ,
  \end{eqnarray}

Now, from \eqref{dif_estimate}  it follows that
   $\|I u \|_{   D(\widetilde{A}_{\epsilon})  }=
 \|u\|_{D(\widetilde{A}_{\epsilon})  } = \| \widetilde{A}_{\epsilon} u  \|_{H^{-1}(\Omega)} 
 \leq  (1+ \tau(\epsilon)) 
 \|   \widetilde{A} u \|_{H^{-1}(\Omega)} = 
 (1+ \tau(\epsilon))  \|u\|_{D(\widetilde{A})  }$, where
 $ \tau(\epsilon) \to 0$, uniformly for
  $u \in  D(\widetilde{A})$.
Thus  $\|I \|_{\mathcal{L}(D(\widetilde{A})\,,\,D(\widetilde{A}_{\epsilon})   )} \leq  (1+ \tau(\epsilon))$ as $\epsilon \to 0$ and it follows from \eqref{estimate_L} that
  \[
  \|u\|_{[D(\widetilde{A}_{\epsilon})\,,\, H^{-1}({\Omega})]_{\alpha} } 
 \leq   (1+ \tau(\epsilon))^{\alpha}   \|u\|_{[D(\widetilde{A})\,,\, H^{-1}({\Omega})]_{\alpha} },  \] 
 as $\epsilon \to 0$. The reverse inequality follows similarly.
   \eproof

 \vspace{3mm}


Using the results of Theorem \ref{sec_scale} and \ref{uniform_scale}, we
 can now  pose the  problem (\ref{nonlinBVP_fix})  as an abstract problem in a (fixed)  scale of Banach spaces $ \{X^{\beta},
 \, -\frac{1}{2} \leq \beta \leq 0 \}  $.
 \begin{equation}\label{abstract_scale}
\begin{array}{lll}
\left\{
\begin{array}{lll}
u_t + (A_{\epsilon})_{\beta}u = (H_{\epsilon})_{\beta}u \, , \, t>t_0\,;
\\
u(t_0)=u_0 \in X^{\eta}\, ,
\end{array}
\right .
\end{array}
\end{equation}
where 
\begin{equation} \label{defH}
  (H_{\epsilon})_{\beta} = H(\cdot,\epsilon):=(F_{\epsilon})_{\beta}
  + (G_{\epsilon})_{\beta} :X ^{\eta} \to X ^{\beta},
  \ \ \epsilon >0 \textrm{ and } 0 \leq \eta \leq \beta +1,
  \end{equation}
\begin{itemize}
\item[(i)] $(F_{\epsilon})_{\beta} = F(\cdot,\epsilon) :X ^{\eta} \to X ^{\beta}$ is given by
\begin{eqnarray}\label{Fh}
\left\langle F(u,\epsilon)\,,\,\Phi \right\rangle_{\beta \,,\, - \beta} =\displaystyle\int_{\Omega} f(u)\,\Phi\,dx,  \ \ \textrm{ for any } \Phi \in X^{{-} \beta},
\end{eqnarray}
\item[(ii)] $(G_{\epsilon})_{\beta} = G(\cdot, {\epsilon}) :
 X ^{\eta} \to X ^{\beta} $ is given by 
\begin{eqnarray}\label{Gh}
\left\langle G(u,\epsilon)\,,\,\Phi \right\rangle_{\beta \,,\,  - \beta } =\displaystyle\int_{\partial \Omega}g(\gamma(u))\,\gamma(\Phi)\left|\frac{J_{\partial \Omega}h_\epsilon}{Jh_\epsilon}\right|\,d\sigma(x)\,, \ \ \textrm{ for any } \Phi \in X^{ {-} \beta},
\end{eqnarray}
where $\gamma$ is the trace map  {and $J_{\partial\Omega}h_\epsilon$ is the determinant of the Jacobian matrix of the diffeomorphism $h_\epsilon: \partial\Omega \longrightarrow \partial h_\epsilon(\Omega)$}.
\end{itemize}

\section{Local well-posedness}
\label{wellposed}
  In order to prove local  well-posedness for the abstract problem 
  (\ref{abstract_scale}), we will need the
 following  growth conditions for the functions $f,g: \R \to \R$:

\begin{enumerate}
\item  $f$ is in  $C^1(\R,\R)$ and there exist real numbers  $\lambda_1 >0$ e $L_1 >0$ such that
\begin{equation} \label{growthf}
|\,f(u_1)-f(u_2)\,| \leq L_1\,(1 + |\,u_1\,|^{\,\lambda_1} + |\,u_2\,|^{\,\lambda_1})|\,u_1-u_2\,|,  \ \textrm{for all} \   u_1,u_2 \in \R.
\end{equation}
\item  $g$ is in $C^2(\R,\R)$ and there exists a 
 real number  $\lambda_2 >0$ e $L_2 >0$ such that
\begin{equation} \label{growthg}
|\,g(u_1)-g(u_2)\,| \leq L_2\,(1 + |\,u_1\,|^{\,\lambda_2} + |\,u_2\,|^{\,\lambda_2})|\,u_1-u_2\,|   \ \textrm{for all} \   u_1,u_2 \in \R.
\end{equation}
\end{enumerate}

\begin{lema}\label{Fbem}
Suppose that $f$ satisfies the growth condition (\ref{growthf}) and $  {\eta >} \frac{1}{2} - \frac{1}{2 (\lambda_{ {1} } +1)}$. 
Then, the operator 
 $ {(F_{\epsilon})_\beta}= F  :X^{\eta}   {\rightarrow}   X^{\beta}$ given by 
 (\ref{Fh}) is well defined, and bounded in bounded sets.
\end{lema}
\proof
 If
 $u \in X^{\eta} $  and
 $\Phi \in X^{-\beta}$, 
$$
\begin{array}{lll}
\big|\left\langle F(u\,,\,\epsilon)\,,\,\Phi \right\rangle_{ {\beta, -\beta}}\big| &\leq& \displaystyle\int_\Omega |\,f(u)\,|\,|\,\Phi\,|\,dx \\
&\leq& L_1\displaystyle\int_\Omega |\,u\,|\,|\,\Phi\,|\,dx + L_1\displaystyle\int_\Omega |\,u\,|^{\,\lambda_1 + 1}\,|\,\Phi\,|\,dx + \displaystyle\int_\Omega |\,f(0)\,|\,|\,\Phi\,|\,dx \\
&\leq& L_1||u||_{L^2(\Omega)}||\Phi||_{L^2(\Omega)} + L_1||u^{\lambda_1 + 1}||_{L^2(\Omega)}||\Phi||_{L^2(\Omega)} +
||f(0)||_{L^2(\Omega)}||\Phi||_{L^2(\Omega)} \\
&=& L_1||u||_{L^2(\Omega)}||\Phi||_{L^2(\Omega)} + L_1||u||^{\lambda_1 + 1}_{L^{2(\lambda_1 + 1)} (\Omega)}||\Phi||_{L^2(\Omega)} +
||f(0)||_{L^2(\Omega)}||\Phi||_{L^2(\Omega)}\,.\\
\end{array}
$$
\par By  {Theorem} \ref{imbed_xalpha} we have   $X^{-\beta} 
 \subset L^2(\Omega)$,
 $X^{\eta} \subset L^2(\Omega)$
 and $X^{\eta} \subset L^{2(\lambda_1 + 1)}(\Omega)$,
with embedding constants $K_1$, $K_2$ and $K_3$, respectively.
 Thus
$$
\begin{array}{lll}
\big|\left\langle F(u,\epsilon)\,,\,\Phi \right\rangle_{ { \beta \,,\,-\beta}}\big| 
&\leq& L_1K_1K_2||\,u\,||_{X^\eta}||\,\Phi\,||_{X^{-\beta}} + L_1 {K_1}K_3^{ {\lambda_1 + 1} }||\,u\,||^{\,\lambda_1 + 1}_{X^{\eta}}||\,\Phi\,||_{X^{-\beta}} \\
&& +\,
 {K_1}||\,f(0)\,||_{L^2(\Omega)}||\,\Phi\,||_{X^{-\beta}}\,.
\end{array}
$$
Therefore, if  $u \in X^{\eta}$ and $\Phi \in X^{-\beta}$,
we have
$$
\big|\big|\,F(u,\epsilon)\,\big|\big|_{X^{\beta}}
\leq L_1K_1K_2||\,u\,||_{X^{\eta}} + L_1 {K_1}K_3^{ {\lambda_1 + 1}}||\,u\,||^{\,\lambda_1 + 1}_{X^{\eta}} +
 {K_1}||\,f(0)\,||_{L^2(\Omega)}\, ,
$$
which concludes the proof.
\eproof

\begin{lema}\label{Flip}
 Suppose that $f$ satisfies the growth condition (\ref{growthf}) and let  $p$ and $q$ be conjugated exponents, with $ \frac{1}{\lambda_1} <
p < \infty$.
Then, if $\eta>  \max\left\{\frac{1}{2} - \frac{1}{2p \lambda_1},
\frac{1}{2} - \frac{1}{2q} \right\}$, the
   operator 
 $F(u,\epsilon)= F(u) :X^{\eta}\times \R   {\rightarrow} X^{\beta}$ given by (\ref{Fh}) 
is locally  Lipschitz continuous in $u$.
\end{lema}
\proof
If  $u_1, u_2 \in X^{ {\eta}}$, we have
$$
\begin{array}{lll}
\left|\left\langle F(u_1,\epsilon) - F(u_2,\epsilon)\,,\,\Phi \right\rangle_{\beta, -\beta}\right |&=&
\left|\displaystyle\int_{\Omega} [f(u_1)-f(u_2)]\,\Phi\,dx\,\right| \\
&\leq& \displaystyle\int_{\Omega} L_1(1 + |\,u_1\,|^{\,\lambda_1} + |\,u_2\,|^{\,\lambda_1})\,|\,u_1-u_2\,|\,|\,\Phi\,|\,dx  \\
&\leq& L_1\Bigg\{\displaystyle\int_{\Omega} (1 + |u_1|^{\,\lambda_1} + |u_2|^{\,\lambda_1})^{\,2}|u_1-u_2|^{\,2}dx \Bigg\}^\frac{1}{2} 
\Bigg\{\displaystyle\int_{\Omega} |\,\Phi\,|^{\,2}dx \Bigg\}^\frac{1}{2} \\
&\leq& L_1\Bigg\{\displaystyle\int_{\Omega} (1 + |\,u_1\,|^{\,\lambda_1} + |\,u_2\,|^{\,\lambda_1})^{\,2p}\Bigg\}^\frac{1}{2p}\Bigg\{\displaystyle\int_{\Omega} |\,u_1-u_2\,|^{\,2q}dx \Bigg\}^\frac{1}{2q}\\
&& \cdot \, \Bigg\{\displaystyle\int_{\Omega} |\,\Phi\,|^{\,2}dx \Bigg\}^\frac{1}{2}\\
&=& L_1||\,1 + |\,u_1\,|^{\,\lambda_1} + |\,u_2\,|^{\,\lambda_1}\,||_{L^{2p}(\Omega)}\,||\,u_1-u_2\,||_{L^{2q}(\Omega)}
||\,\Phi\,||_{L^2(\Omega)} \\
&\leq& L_1\big( |\Omega|^{ {\frac{1}{2p}}} + ||u_1^{\lambda_1}||_{L^{2p}(\Omega)} + ||u_2^{\,\lambda_1}||_{L^{2p}(\Omega)}\big)\,||u_1-u_2||_{L^{2q}(\Omega)}
||\,\Phi\,||_{L^2(\Omega)} \\
&=& L_1\big( |\Omega|^{ {\frac{1}{2p}}} + ||u_1||^{\lambda_1}_{L^{2p\lambda_1}(\Omega)} + ||u_2||^{\lambda_1}_{L^{2p\lambda_1}(\Omega)}\big)||u_1-u_2||_{L^{2q}(\Omega)}||\Phi||_{L^2(\Omega)} \\
&\leq& L_1\big(\,|\,\Omega\,|^{ {\frac{1}{2p}}} + K_4^{  {\lambda_1}}||\,u_1\,||^{\lambda_1}_{X^{\eta}} + K_4^{  {\lambda_1}}||\,u_2\,||^{\,\lambda_1}_{X^{\eta}}\,\big)\\
&& \cdot \,K_5||\,u_1-u_2\,||_{X^{\eta}}
  {K_1}||\,\Phi\,||_{X^{-\beta}}\,, 
\end{array} 
$$
where $|\,\Omega\,|$ is the measure of  $\Omega$. 

\par By  {Theorem} \ref{imbed_xalpha} we have $X^{-\beta} 
 \subset L^2(\Omega)$, 
 $X^{\eta} \subset L^{2p \lambda_1}(\Omega)$,
 and $X^{\eta} \subset L^{2q}(\Omega)$, 
with embedding constants  {$K_1$}, $K_4$ and $K_5$, respectively.
 Thus

$$
\begin{array}{lll}
\left|\,\left\langle F(u_1,\epsilon) - F(u_2,\epsilon)\,,\,\Phi \right\rangle_{\beta\,,\, -\beta}\,\right|
&\leq& L_1\big(\,|\,\Omega\,|^{  {\frac{1}{2p}}} + K_4^{  {\lambda_1}}||\,u_1\,||^{\,\lambda_1}_{X^{\eta}} + K_4^{  {\lambda_1}}||\,u_2\,||^{\,\lambda_1}_{X^{\eta}} \big) \\
&&\cdot \,K_5||\,u_1-u_2\,||_{X^{\eta}}
  {K_1}||\,\Phi\,||_{X^{-\beta}}\,. 
\end{array} 
$$
 If  $U$ is a bounded subset of  $X^{\eta}$ and $u_1$, $u_2 \in U$,
we have
$$
||\,F(u_1\,,\,\epsilon)- F(u_2\,,\,\epsilon)\,||_{X^{\beta}}
\leq K_{\lambda_1\,,\,U}\,||\,u_1-u_2\,||_{X^{\eta}},
$$
where $K_{\lambda_1\,,\,U}$ is a positive constant depending
 on  $\lambda_1$ and  $U$. This concludes the proof.
\eproof 

\par For  the regularity properties of $(G_{\epsilon})_{\beta} $ we will need to compute the  function
$$\theta(x, \epsilon) := 
\left|\,\displaystyle\frac{J_{\partial\Omega}h_\epsilon}{Jh_\epsilon}(x) \,\right|,\,x =(x_1,x_2) \in \R^2,$$
  in  each of the four  segments of $\partial \Omega$.
 For $I_1:= \left\{  (  {x_1},1) \, | \, 0\leq   {x_1} \leq 1 \right\}$, we have
 $h_\epsilon(  {x_1}, 1) = (  {x_1},   1 + \epsilon \,sen(  {x_1}/ \epsilon^\alpha) )$, so
$ \frac{d}{d   {x_1}} h_\epsilon(  {x_1}, 1) = Dh_\epsilon(  {x_1},1) \cdot (1,0) = 
  (1, \epsilon^{1- \alpha} \cos ( {x_1}/ \epsilon^\alpha))$. Therefore the restriction
 of $h_\epsilon$ to $I_1$ takes the unit vector $(1,0)$ tangent to 
 $ I_1$ to the vector $(1, \epsilon^{1- \alpha} \cos ( {x_1}/ \epsilon^\alpha))$.
 It follows that $ \|D {h_\epsilon}_ {\restriction I_1}( {x_1},1)\| =
   |J_{\partial\Omega}h_\epsilon ( {x_1},1) | =    
   \sqrt{1 +  {\epsilon}^{\,2-2\alpha}{cos}^{\,2}( {x_1}/  {\epsilon}^{\,\alpha})} $.
  
Thus
 \begin{eqnarray} \label{theta_1}
\theta(x, \epsilon)& := &  \frac{\sqrt{1 +  {\epsilon}^{\,2-2\alpha}{cos}^{\,2}( {x_1}/  {\epsilon}^{\,\alpha})}  }{1+ \epsilon \sin( {x_1}/\epsilon^{\alpha})}, 
  \textrm{ for }  x \in I_1.
\end{eqnarray}

 Similar computations  give
\begin{eqnarray}\label{theta_2}
 \theta(x, \epsilon)& := & \frac{1+ \epsilon \sin( {1}/\epsilon^{\alpha}) }{1+ \epsilon \sin( {1}/\epsilon^{\alpha})} = 1  \textrm{ for }  x \in I_2:= \left\{  (1, {x_2}) \, | \, 0
\leq  {x_2} \leq 1 \right\}, \nonumber \\
\theta(x, \epsilon) &  := &  \frac{1 }{1+ \epsilon \sin( {x_1}/\epsilon^{\alpha}) }  \textrm{ for }  x \in I_3:= \left\{  ( {x_1},0) \, | \, 0
\leq  {x_1} \leq 1 \right\}, \nonumber \\
\theta(x, \epsilon) &  := & \frac{1 }{1} = 1  \textrm{ for }  x \in I_4:= \left\{  (0, {x_2}) \, 
  | \,  {0} \leq  {x_2} \leq 1 \right\}. 
\end{eqnarray}

\begin{lema}\label{Gbem}
Suppose that $g$ satisfies the growth condition (\ref{growthg}),
  $\eta > \frac{1}{2}- \frac{1}{4(\lambda_2+1)}$,  $\beta < -\frac{1}{4}$ and
 $\epsilon_0 < 1$. Then, the operator $(G_{\epsilon})_{\beta}  {=G} :X^{\eta}   {\rightarrow} X^{\beta} $
given by (\ref{Gh})  is well defined, for $0 \leq \epsilon <\epsilon_0$
and  bounded in bounded sets, uniformly in $\epsilon$.
\end{lema}
\proof  If $u \in X^{\eta}$ and $\Phi \in X^{-\beta}$, we have
$$
\begin{array}{lll}
\big|\left\langle G(u,\epsilon)\,,\,\Phi \right\rangle_{\beta \,,\, -\beta }\big| &\leq&\displaystyle\int_{\partial \Omega} |\,g(\gamma(u))\,|\,|\,\gamma(\Phi)\,|\left|\frac{J_{\partial \Omega}h_\epsilon}{Jh_\epsilon}\right|\,d\sigma(x) \\
&\leq& ||\theta||_\infty \displaystyle\int_{\partial \Omega} L_2\big[|\gamma(u)| + |\gamma(u)|^{\lambda_2 + 1}\big]|\gamma(\Phi)| + |g(\gamma(0))|\,|\gamma(\Phi)| \, d\sigma(x) \\
&\leq& ||\,\theta\,||_\infty  L_2\,\displaystyle\bigg(\int_{\partial \Omega}\,|\,\gamma(u)\,|^{\,2}d\sigma(x)\bigg)^\frac{1}{2}\displaystyle\bigg(\int_{\partial \Omega}|\,\gamma(\Phi)\,|^{\,2}d\sigma(x)\bigg)^\frac{1}{2} \\
&&+\, ||\,\theta\,||_\infty  L_2\,\displaystyle\bigg(\int_{\partial \Omega}\,|\,\gamma(u)\,|^{\,2(\lambda_2 + 1)}d\sigma(x)\bigg)^\frac{1}{2}\displaystyle\bigg(\int_{\partial \Omega}|\,\gamma(\Phi)\,|^{\,2}d\sigma(x)\bigg)^\frac{1}{2} \\
&&+\, ||\,\theta\,||_\infty||\,g(\gamma (0))\,||_{L^2(\partial \Omega)}\,||\,\gamma(\Phi)\,||_{L^2(\partial \Omega)} \\
&=& ||\,\theta\,||_\infty L_2\,\big[\,||\,\gamma(u)\,||_{L^2(\partial \Omega)} + ||\,\gamma(u)\,||_{L^{2(\lambda_2 + 1)}(\partial \Omega)}^{\,\lambda_2 + 1}\,\big]\,||\,\gamma(\Phi)\,||_{L^2(\partial \Omega)} \\
&&+\, ||\,\theta\,||_\infty||\,g(\gamma (0))\,||_{L^2(\partial \Omega)}\,||\,\gamma(\Phi)\,||_{L^2(\partial \Omega)}\,,
\end{array}
$$
where  $\|\theta\|_\infty = \sup \left\{ |\theta(x, \epsilon)| \, | \, 
 x\in \partial \Omega, \,  0 \leq \epsilon
 \leq  \epsilon_0 \right\}$ is finite by 
(\ref{theta_1}) and (\ref{theta_2}). If $s= 2\eta$,   then $X^{\eta} \subset H^{s}(\Omega)$, by Theorem
\ref{imbed_xalpha} and   $\gamma : H^{s}(\Omega) \mapsto
 L^{2 (\lambda_2+1)} (\partial\Omega)$, by Theorem \ref{trace_hk}.
 Thus, we obtain
$||\gamma(\Phi)||_{L^2(\partial \Omega)} \leq \overline{K}_1\,||\Phi||_{X^{-\beta}}\,,\
||\gamma(u)||_{L^{2(\lambda_2 + 1)}(\partial \Omega)} \leq \overline{K}_2\,||u||_{X^{\eta}}\,, \
||\gamma(u)||_{L^2(\partial \Omega)} \leq \overline{K}_3\,||u||_{X^{\eta}}\, ,
$
where $\overline{K}_1$, $\overline{K}_2$ and $\overline{K}_3$
 are embedding constants. It follows that 
$$
\begin{array}{lll}
\big|\big|\,G(u,\epsilon)\,\big|\big|_{X^{\beta}} 
&\leq& L_2 \overline{K}_3\overline{K}_1||\,\theta\,||_\infty ||\,u\,||_{X^{\eta}} + L_2 \overline{K}_2^{ {\lambda_2 + 1}}\overline{K}_1||\,\theta\,||_\infty ||\,u\,||_{X {^\eta}} \\
& &+\, \overline{K}_1||\,\theta||_\infty||\,g(\gamma (0))\,||_{L^2(\partial \Omega)}\,,
\end{array}
$$
proving that  $(G_{\epsilon})_{\beta} $ is well defined.
\eproof

\begin{lema}\label{Glip}
Suppose that $g$ satisfies the growth condition (\ref{growthg}) and let $p$ and $q$ be conjugated exponents, with $ \frac{1}{2\lambda_2} <
p < \infty$  and $\epsilon_0 < 1$.
Then, if $ \eta >  \max\left\{\frac{1}{2} - \frac{1}{4p \lambda_2},
\frac{1}{2} - \frac{1}{4q} \right\}$  and $ \beta <- \frac{1}{4}$,
 the operator
 $G(u,\epsilon)   {=G(u)} :X^{\eta} \times [0, \epsilon_0] \to X^{\beta}$ given
 by  (\ref{Gh})
 is uniformly continuous in $\epsilon$, for $u$ in bounded sets of $X^{\eta}$
   and
  locally Lipschitz continuous in $u$, uniformly in  
  $\epsilon$. 
\end{lema}
\proof
 We first show that $(G_{\epsilon})_{\beta} $  is locally Lipschitz continuous in  $u \in X^{\eta}$. Let $u_1, u_2 \in X^{\eta}$, $\Phi \in X^{-\beta}$  and $ \epsilon \in [0, \epsilon_0]$. Then
$$
\begin{array}{lll}
\left|\left\langle G(u_1,\epsilon) - G(u_2,\epsilon),\Phi \right\rangle_{\beta,-\beta}\right |
&\leq& \displaystyle\int_{\partial\Omega} |\,g(\gamma(u_1))-g(\gamma(u_2))\,|\,\big|\,\gamma(\Phi)\,\big|\,\left|\frac{J_{\partial\Omega} h_\epsilon}{Jh_\epsilon}\right|\,d \sigma (x) \\
&\leq& \displaystyle\int_{\partial\Omega} L_2\big(1 + |\gamma(u_1)|^{\,\lambda_2} + |\gamma(u_2)|^{\,\lambda_2}\big)|\gamma(u_1)-\gamma(u_2)|\left|\gamma(\Phi)\right|\left|\frac{J_{\partial\Omega} h_\epsilon}{Jh_\epsilon}\right|d\sigma (x) \\
&\leq &
L_2||\theta||_\infty \displaystyle\int_{\partial\Omega} \big(1 + |\gamma(u_1)|^{\,\lambda_2} + |\gamma(u_2)|^{\,\lambda_2}\big)|\gamma(u_1)-\gamma(u_2)|\left|\gamma(\Phi)\right|\,d \sigma (x) \\
&\leq&
L_2||\theta||_\infty \Bigg\{\displaystyle\int_{\partial\Omega} \big(1 + |\gamma(u_1)|^{\,\lambda_2} + |\gamma(u_2)|^{\,\lambda_2}\big)^{\,2}\,|\gamma(u_1)-\gamma(u_2)|^{\,2}\,d\sigma(x)\Bigg\}^\frac{1}{2} \\
&\cdot&\bigg\{\displaystyle\int_{\partial\Omega}\left|\,\gamma(\Phi)\,\right|^{\,2}\,d \sigma (x)\bigg\}^\frac{1}{2} \\
&\leq &
L_2||\,\theta\,||_\infty \Bigg\{\displaystyle\int_{\partial\Omega} \big(1 + |\,\gamma(u_1)\,|^{\,\lambda_2} + |\,\gamma(u_2)\,|^{\,\lambda_2}\big)^{\,2p}\,d\sigma(x)\Bigg\}^\frac{1}{2p} \\
&\cdot& \Bigg\{\displaystyle\int_{\partial\Omega} |\,\gamma(u_1)-\gamma(u_2)\,|^{\,2q}d\sigma(x)]\Bigg\}^\frac{1}{2q}||\,\gamma(\Phi)\,||_{L^2(\partial \Omega)} \\
&\leq&
L_2||\,\theta\,||_\infty \big(\,|\,\partial\Omega\,| {^\frac{1}{2p}} + ||\,\gamma(u_1)\,||^{\,\lambda_2}_{L^{2p\lambda_2}(\partial \Omega)} + ||\,\gamma(u_2)\,||^{\lambda_2}_{L^{2p\lambda_2}(\partial \Omega)}\,\big) \\
&\cdot& ||\,\gamma(u_1)-\gamma(u_2)\,||_{L^{2q}(\partial \Omega)}||\,\gamma(\Phi)\,||_{L^2(\partial \Omega)}\,,
\end{array}
$$
where $|\,\partial\Omega\,|$ is the measure of $\partial\Omega$,  $1 < p,q < \infty$ 
 are conjugated exponents with  $p > \displaystyle\frac{1}{2\lambda_2}$ and 
 $\|\theta\|_\infty = \sup \left\{ |\theta(x, \epsilon)| \, | \,
 x\in \Omega, \,  0 \leq \epsilon
 \leq  \epsilon_0 \right\}$   .  Reasoning as  in Lemma
 \ref{Gbem}, we obtain  the following estimates:
$$
\begin{array}{lll}
||\,\gamma(\Phi)\,||_{L^2(\partial \Omega)} &\leq& \overline{K}_1\,||\,\Phi\,||_{X^{-\beta}}\,,\\
||\,\gamma(u)\,||_{L^{2p\lambda_2}(\partial \Omega)} &\leq& \overline{K}_4\,||\,u\,||_{X^{\eta}}\,, \\
||\,\gamma(u)\,||_{L^{2q}(\partial \Omega)} &\leq& \overline{K}_5\,||\,u\,||_{X^{\eta}}\,, \\
\end{array}\,
$$
 where $\overline{K}_1$, $\overline{K}_4$ e $\overline{K}_5$ are embedding constants.
 Thus, we obtain:
$$
\begin{array}{lll}
\left|\,\left\langle G(u_1,\epsilon) - G(u_2,\epsilon)\,,\,\Phi \right\rangle_{\beta\,,\,-\beta}\,\right| 
&\leq& L_2||\,\theta\,||_\infty \big(\,|\,\partial\Omega\,| {^\frac{1}{2p}} + \overline{K}_4 ^{\lambda_2}||\,u_1\,||^{\,\lambda_2}_{X^{\eta}} + \overline{K}_4^{\lambda_2}||\,u_2\,||^{\,\lambda_2}_{X^{\eta}}\,\big)\\
&&\cdot\, \overline{K}_5||\,u_1 - u_2\,||_{X^{\eta}}\overline{K}_1||\,\Phi\,||_{X^{-\beta}}\,.
\end{array}
$$
 If  $U$  is a bounded subset of  $X^{\eta}$,  $u_1$, $u_2 \in U$,  we have
$$||\,G(u_1,\epsilon) - G(u_2,\epsilon)\,||_{X^{\beta}}\leq \overline{K}_{\lambda_2\,,\,U}\,||\,u_1-u_2\,||_{X^{\eta}}\,,$$
 where $\overline{K}_{\lambda_2\,,\,U}$ is a positive constant
 depending on  $\lambda_2$ and $U$. Therefore,  $(G_{\epsilon})_{\beta} $ is locally Lipschitz in $u$.

\par Now, if  $u \in X^{\eta}$, $\Phi \in X^{-\beta}$  and
 $ \epsilon_1, \epsilon_2 \in [0, \epsilon_0]$, we have
 \begin{eqnarray*}
\big|\langle G(u,{\epsilon_1})-G(u,{\epsilon_2}),\Phi\rangle_{\beta,-\beta}\big| &\leq&
\displaystyle\int_{\partial\Omega}|\,\gamma(g(u))\,|\,|\,\gamma(\Phi)\,|\left|\,\left(\,\left|\frac{J_{\partial \Omega}h_{\epsilon_1}}{Jh_{\epsilon_1}}\right|-\left|\frac{J_{\partial \Omega}h_{\epsilon_2}}{Jh_{\epsilon_ 2}}\right|\,\right)\,\right|\,d\sigma(x)  \\
&\leq&
\|\theta_{\epsilon_1} - \theta_{\epsilon_2} \|_{\infty}
\displaystyle\int_{\partial\Omega} |\,g(\gamma(u))\,|\,|\,\gamma(\Phi)\,|\, d\sigma(x),
\end{eqnarray*}
 and $\|\theta_{\epsilon_1}- \theta_{\epsilon_2} \|_\infty =
 \sup \left\{ |\theta(x, \epsilon_1  )  - \theta(x, \epsilon_1  ) | \, | \,
 x\in \Omega, \,  \right\} \to 0 $ as $ |\epsilon_1 - \epsilon_2| \to 0 $,
by (\ref{theta_1}) and (\ref{theta_2}). \eproof

\begin{teo}\label{loc_exist}
  Suppose  $f$ and $g$ satisfy the growth conditions 
 (\ref{growthf}) and (\ref{growthg}), respectively. Suppose also, that
 $\beta$ and $\eta$ satisfy the  {hypotheses} of Lemmas \ref{Fbem}, \ref{Flip},
\ref{Gbem},  \ref{Glip}, $\eta < 1+ \beta$ and $ \epsilon> 0$ is sufficiently small.  
 Then, for any $(t_0, u_0) \in \R \times X^\eta$, the  problem (\ref{abstract_scale}) 
 has a unique solution $u(t,t_0, u_0,\epsilon)$  with initial value $u(t_0) = u_0$. The map $\epsilon \mapsto u(t,t_0, u_0,\epsilon) \in X^{ \eta}$ is continuous at
 $\epsilon = \epsilon_0$,  {uniformly} for $u_0$ in bounded sets of $X^{\eta}$
  and $t_0 \leq t \leq T < \infty$.
\end{teo}
\proof
 From Theorem \ref{sec_scale}  it follows that $ {(}A_{ {\epsilon}} {)_\beta}$ is a sectorial operator  in
 $X^{\beta}_{ {\epsilon}}$, with domain  $X^{1+\beta}_{ {\epsilon}}$, if $\epsilon $ is small enough. From Lemmas  \ref{Fbem}, \ref{Flip},
\ref{Gbem} and  \ref{Glip} it follows  that $ {(}H_{ {\epsilon}} {)_\beta}$ is well defined and locally Lipschitz continuous in $X^\eta$, and bounded in bounded sets of  $X^\eta$. The result follows then from Theorems \ref{333} and \ref{continabsALT}. 
\eproof

\section{ {Lyapunov} functionals  and global existence}
\label{lyapunov}

 We now want to show that the solutions given by Theorem \ref{loc_exist} are globally defined, if a additional (dissipative)  hypotheses on $f$ and $g$ is assumed.
  We start by stating these hypotheses:

  There exist constants  $c_0$,  {$d_0$} and $d_0'$ such that
\begin{equation} \label{dissipative}
\displaystyle\limsup_{|\,u\,| \to \infty }\frac{f(u)}{u} \leq c_0\,, \quad
\displaystyle\limsup_{|\,u\,| \to \infty }\frac{g(u)}{u} \leq d_0'
\end{equation}

and, if  $d_0 > d_0'$, the first eigenvalue $\mu_1$ of the problem
\begin{equation} \label{compet}
\left\{
\begin{array}{lll}
-\Delta u + (a - c_0)u = \mu u \,\,\mbox{em} \,\,\Omega   \\
\displaystyle\frac{\partial u}{\partial N_\Omega}=d_0\,u \,\, \mbox{em} \,\, \partial\Omega
\end{array}
\right.
\end{equation}
is positive.

\begin{rem}
 Observe that the hypotheses (\ref{dissipative}) and (\ref{compet}) still hold
 for the perturbed operator $ h_{\epsilon}^{*}  \Delta_{ {\Omega_\epsilon}} {h_{\epsilon}^{*}}^{-1} $, with perturbed
 boundary conditions
 $ h_{\epsilon}^{*} \displaystyle\frac{\partial u}{\partial N_{\Omega_{ {\epsilon}}}}{h_{\epsilon}^{*}}^{-1}$, if $\epsilon$ is small enough,  since  the eigenvalues change continuously with $\epsilon$ by (\ref{dif_estimate}).
\end{rem}

In order to prove global existence, we  work first in the natural ``energy space''  $H^{1}(\Omega)$, that is we choose $\eta= \frac{1}{2}$
 (and $\beta < -\frac{1}{4}$). It is also convenient to work first in the
 \emph{perturbed domain} $\Omega_{\epsilon}$.  
 More precisely, we consider initially  the following
 abstract version of problem
  (\ref{nonlinBVP})  in the  Banach space $Y^{\beta}$, where
 $Y^{\alpha}, \  -\frac{1}{2} \leq \alpha \leq 1$,  now denote the fractional powers of
 the  operator $- \Delta_{\Omega_\epsilon} + aI $  in the perturbed domain
 $\Omega_\epsilon$,  

 \begin{equation}\label{abstract_perturbed}
\begin{array}{lll}
\left\{
\begin{array}{lll}
v_t + A_{\beta}v = H_{\beta}v \, , \, t>t_0\,;
\\
v(t_0)=v_0 \in H^1(\Omega_\epsilon) \, ,
\end{array}
\right .
\end{array}
\end{equation}
where 
 $ H_{\beta} :=F_{\beta}
  +  G_{\beta} :  H^1(\Omega_\epsilon) \to X ^{\beta},
 \ \ \epsilon >0$,  
\begin{itemize}
\item[(i)] $ F_{\beta}  :  H^1(\Omega_\epsilon) \to X ^{\beta}$ is given by
\begin{eqnarray*}
\left\langle F(v)\,,\,\Psi \right\rangle_{\beta \,,\, - \beta} =
\displaystyle\int_{\Omega_{\epsilon}} f(v)\,\Psi\,dy,  \ \ \textrm{ for any } \Psi
 \in Y^{ {-}\beta},
\end{eqnarray*}
\item[(ii)] $ G_{\beta}  :  H^1(\Omega_\epsilon)
  \to X ^{\beta} $ is given by 
\begin{eqnarray*}
\left\langle G(v)\,,\,\Psi \right\rangle_{\beta \,,\,  - \beta } =\displaystyle\int_{\partial \Omega {_\epsilon}}g(\gamma(v))\,\gamma(\Psi)  {\left|\frac{J_{\partial\Omega}h_\epsilon}{Jh_\epsilon}\right|}\,d\sigma(y)\,, \ \ \textrm{ for any } \Psi \in Y^{ {-}\beta},
\end{eqnarray*}
where $\gamma$ is the trace map.
\end{itemize}

 It is not difficult to show that  $v$ is a solution of  
(\ref{abstract_perturbed}) if and only if $u = v \circ h {_\epsilon}$ is a solution
 of (\ref{abstract_scale}), with $\eta = \frac{1}{2}$. Therefore, the local
 well-posedness of (\ref{abstract_perturbed})   follows immediately
 from  {Theorem} \ref{loc_exist}. 

 We now  prove the existence of a Lyapunov functional for the dynamical system generated by 
(\ref{abstract_perturbed}). 
 
 \begin{lema}\label{funly}
   Suppose that  the hypotheses  of
   Theorem \ref{loc_exist} are satisfied, with $\eta=\frac{1}{2}$
   and, additionally, that $f$ and $g$ satisfy the
   dissipative conditions (\ref{dissipative}) and (\ref{compet}) and
 consider the map 
$$
\begin{array}{llll}
W_\epsilon :& H^1(\Omega_{\epsilon}) & {\longrightarrow} &\R \\
& \,\,\,\,\, v  & \longmapsto & W_\epsilon(v)= \displaystyle\frac{1}{2}\displaystyle\int_{\Omega_{\epsilon}}|\,\nabla v\,|^2dx + \displaystyle\frac{a}{2} 
\displaystyle\int_{\Omega_{\epsilon}}|\,v\,|^2dx - \displaystyle\int_{\Omega_{\epsilon}}F(v)\,dx - \displaystyle\int_{\partial\Omega_{\epsilon}}G(\gamma(v))\,dS
\end{array} 
$$
where $a$ is a positive number and  $F,G : \R \to \R$ 
 are primitives of $f$ and $g$ respectively. 
 Then, if $\epsilon>0$ is sufficiently small,  $W_\epsilon$
 is a Lyapunov functional for the problem (\ref{abstract_perturbed}) and
  there constants $K_1(\epsilon)$ and $K_2 (\epsilon)$ such that 
  $W_{\epsilon}(v)  \leq K_1(\epsilon) \|v\|_{H^1(\Omega_\epsilon)} +K_2 (\epsilon)$,
  for any $v\in H^1 {(}{\Omega_\epsilon} {)}$.
\end{lema} 
\proof  If $v$ is a solution of  (\ref{abstract_perturbed}) in $H^2(\Omega {)}$, we have

$$
\begin{array}{lll}
\displaystyle\frac{d}{dt}W_{\epsilon}(v(t)) 
 &=& \displaystyle\frac{d}{dt} \left(\displaystyle\frac{1}{2}\displaystyle\int_{\Omega_{\epsilon}}|\,\nabla v\,|^{\,2}dx + \displaystyle\frac{a}{2} 
\displaystyle\int_{\Omega_{\epsilon}}|\,v\,|^{\,2}dx - \displaystyle\int_{\Omega_{\epsilon}}F(v)\,dx - \displaystyle\int_{\partial\Omega_{\epsilon}}G(\gamma(v))\,dS\right) \\
&=& \displaystyle\int_{\Omega_{\epsilon}}\left(\displaystyle\frac{d}{dx}v_t\right)\nabla v\,dx + a\displaystyle\int_{\Omega_{\epsilon}}v_tv\,dx - \displaystyle\int_{\Omega_{\epsilon}}v_tf(v)\,dx - \displaystyle\int_{\partial\Omega_{\epsilon}}v_tg(\gamma(v))\,dS \\
&=& -\displaystyle\int_{\Omega_{\epsilon}}v_t\Delta v\,dx + \displaystyle\int_{\partial\Omega_{\epsilon}}v_t\frac{\partial v}{\partial N}dS +a\displaystyle\int_{\Omega_{\epsilon}}v_tv\,dx - \displaystyle\int_{\Omega_{\epsilon}}v_tf(v)\,dx - \displaystyle\int_{\partial\Omega_{\epsilon}}v_tg(\gamma(v))\,dS \\
&=& -\left(\displaystyle\int_{\Omega_{\epsilon}}v_t\Delta v\,dx - a\displaystyle\int_{\Omega_{\epsilon}}v_tv\,dx + \displaystyle\int_{\Omega_{\epsilon}}v_tf(v)\,dx \right) \\
&=& -\displaystyle\int_{\Omega_{\epsilon}}|\,v_t\,|^{\,2}dx \\
&=& -\, ||\,v_t\,||_{L^2(\Omega_{\epsilon})}^{\,2} .  
\end{array}
$$

The equality
$\displaystyle{\frac{d}{dt} {W}_{\epsilon}(v(t)) =  {-}
  ||   \,v_t\, ||_{L^2(\Omega_{\epsilon})}^{\,2}} $ is established, 
 supposing that $ v$ is a solution in $H^2(\Omega)$. But, since both sides are well defined and continuous functions of $v \in X^{\frac{1}{2}}=H^1 (\Omega_{\epsilon})$, it remains true for any solution. Therefore $W_{\epsilon}$ is decreasing along
 the solutions of (\ref{abstract_perturbed}). It is clear from its formula, that $W_\epsilon$ is continuous.  We now want to obtain an estimate for $W_\epsilon(v)$ in terms of the norm of $v \in  H^1 (\Omega_{\epsilon}) $.
 From (\ref{dissipative}), there exist $\epsilon_f >0$ and  $M(\epsilon_f)>0$ such that
 $\displaystyle\frac{f(s)}{s} - c_0 \leq \epsilon_f$ for  $|\,s\,| >M(\epsilon_f)$. Therefore, if $s >0$ we have
  $$\displaystyle\int_{\Omega_{\epsilon}}F(v)\,dx=\displaystyle\int_{\Omega_{\epsilon}}\left(\int_0^v f(s)\,ds \right)dx\leq \displaystyle\int_{\Omega_{\epsilon}}\left(\int_0^v (c_0s + \epsilon_fs)\,ds\right)dx \leq
\frac{c_0}{2}\displaystyle\int_{\Omega_{\epsilon}}|\,v\,|^{\,2}dx + k_0\,,
$$
where $k_0$ depends on  $\epsilon_f$, $M(\epsilon_f)$  and the first eigenvalue of (\ref{compet}).
  A similar argument gives the same estimate for  $s<0$.

  Now, from   (\ref{dissipative}) there exist  $\epsilon_g >0$  and $N(\epsilon_g)>0$ such that
  $\displaystyle\frac{g(s)}{s} - d_0'\leq \epsilon_g$ for  $|\,s\,| >N(\epsilon_g)$. 
If $s >0$ we have  $g(s)\leq d_0's + \epsilon_gs$\,. Therefore 
$$\displaystyle\int_{\partial\Omega_{\epsilon}}G(\gamma(v))dS=\displaystyle\int_{\partial\Omega_{\epsilon}}\left(\int_0^v g(s)ds \right)dS\leq \displaystyle\int_{\partial\Omega_{\epsilon}}\left(\int_0^v (d_0's + \epsilon_gs)ds\right)dS \leq
\frac{d_0'}{2}\displaystyle\int_{\partial\Omega_{\epsilon}}|\gamma(v)|^2dS + k_0'\,,
$$
 where  $k_0'$ depends on  $\epsilon_g$, $N(\epsilon_g)$  and the first eigenvalue of
 (\ref{compet}). We obtain a similar estimate for  $s<0$.
 
 Using  these estimates, we obtain

 $$
\begin{array}{lll}
W_{\epsilon}(v) &\leq& \displaystyle\frac{1}{2}\displaystyle\int_{\Omega_{\epsilon}}|\,\nabla v\,|^{\,2}dx + \displaystyle\frac{a}{2} 
\displaystyle\int_{\Omega_{\epsilon}}|\,v\,|^{\,2}dx + \frac{c_0}{2}\displaystyle\int_{\Omega_{\epsilon}}|\,v\,|^{\,2}dx + k_0 + \frac{d_0'}{2}\displaystyle\int_{\partial\Omega_{\epsilon}}|\,\gamma(v)\,|^{\,2}dS + k_0'\\
& \leq & K_1  {(\epsilon)}\|v \|_{H^1(\Omega_{\epsilon})} + K_2(\epsilon).
\end{array}
$$
On the other hand, we also have

$$
\begin{array}{lll}
W_{\epsilon}(v) &\geq& \displaystyle\frac{1}{2}\displaystyle\int_{\Omega_{\epsilon}}|\,\nabla v\,|^{\,2}dx + \displaystyle\frac{a}{2} 
\displaystyle\int_{\Omega_{\epsilon}}|\,v\,|^{\,2}dx - \frac{c_0}{2}\displaystyle\int_{\Omega_{\epsilon}}|\,v\,|^{\,2}dx - k_0 - \frac{d_0'}{2}\displaystyle\int_{\partial\Omega_{\epsilon}}|\,\gamma(v)\,|^{\,2}dS - k_0' \\
&=& \displaystyle\frac{1}{2}\displaystyle\int_{\Omega_{\epsilon}}|\,\nabla v\,|^{\,2}dx + \displaystyle\frac{(a-c_0)}{2} 
\displaystyle\int_{\Omega_{\epsilon}}|\,v\,|^{\,2}dx  - \frac{d_0'}{2}\displaystyle\int_{\partial\Omega_{\epsilon}}|\,\gamma(v)\,|^{\,2}dS - (k_0 + k_0'),
\end{array}
$$
 and, since  $d_0 > d_0'$, we may choose  $d_0''\neq 0$ such that  
$$d_0 > d_0'' >d_0' \,\,\, \mbox{ and } \,\,\, \left(\frac{d_0}{d_0''}-1\right)\frac{a-c_0}{2} + \frac{\lambda_0}{2} >0,$$ 
and thus
$$
\begin{array}{lll}
W_{\epsilon}(v) &\geq& \displaystyle\frac{1}{2}\displaystyle\int_{\Omega_{\epsilon}}|\,\nabla v\,|^{\,2}dx + \displaystyle\frac{(a-c_0)}{2} 
\displaystyle\int_{\Omega_{\epsilon}}|\,v\,|^{\,2}dx  - \frac{d_0''}{2}\displaystyle\int_{\partial\Omega_{\epsilon}}|\,\gamma(v)\,|^{\,2}dS - (k_0 + k_0') \\
&=& \displaystyle\frac{d_0''}{d_0}\left[\displaystyle\frac{d_0}{d_0''}\displaystyle\frac{1}{2}\displaystyle\int_{\Omega_{\epsilon}}|\,\nabla v\,|^{\,2}dx + \displaystyle\frac{d_0}{d_0''}\displaystyle\frac{(a-c_0)}{2} 
\displaystyle\int_{\Omega_{\epsilon}}|\,v\,|^{\,2}dx  - \frac{d_0}{2}\displaystyle\int_{\partial\Omega_{\epsilon}}|\,\gamma(v)\,|^{\,2}dS -  \displaystyle\frac{d_0}{d_0''}(k_0 + k_0') \right] \\
&=& \displaystyle\frac{d_0''}{d_0}\left[\left(\displaystyle\frac{d_0}{d_0''}-1\right)\displaystyle\frac{1}{2}\displaystyle\int_{\Omega_{\epsilon}}|\,\nabla v\,|^{\,2}dx + \displaystyle\frac{1}{2}\displaystyle\int_{\Omega_{\epsilon}}|\,\nabla v\,|^{\,2}dx + \left(\displaystyle\frac{d_0}{d_0''}-1\right)\displaystyle\frac{(a-c_0)}{2} 
\displaystyle\int_{\Omega_{\epsilon}}|\,v\,|^{\,2}dx \right. \\
&&+\, \left.\displaystyle\frac{(a-c_0)}{2} 
\displaystyle\int_{\Omega_{\epsilon}}|\,v\,|^{\,2}dx - \frac{d_0}{2}\displaystyle\int_{\partial\Omega_{\epsilon}}|\,\gamma(v)\,|^{\,2}dS -  \displaystyle\frac{d_0''}{d_0}(k_0 + k_0') \right] \\
&\geq& \displaystyle\frac{d_0''}{d_0}\left\{\left(\displaystyle\frac{d_0}{d_0''}-1\right)\displaystyle\frac{1}{2}\displaystyle\int_{\Omega_{\epsilon}}|\nabla v|^2dx + \left[\left(\displaystyle\frac{d_0}{d_0''}-1\right)\displaystyle\frac{(a-c_0)}{2} + \displaystyle\frac{\lambda_0}{2}\right]\displaystyle\int_{\Omega_{\epsilon}}|v|^2dx - \displaystyle\frac{d_0''}{d_0}(k_0 + k_0') \right\}.
\end{array}
$$

 From this estimate, we see that  $W_{\epsilon}(v)$ is bounded below by a constant times $||\,v\,||_{H^1(\Omega_{\epsilon})}$  plus a constant depending on the first eigenvalue of (\ref{compet})
 and the nonlinearities.  Thus  $|\,W_{\epsilon}(v)\,| \to \infty$ as  $||\,v\,||_{H^1(\Omega_{\epsilon})} \to \infty$. 
 Finally, if  $v(t)$ is a solution defined for all $ {t \in}\, \R$ and
 $W_{\epsilon}(v(t))=W_{\epsilon}(v_0)$,  with  $v_0 \in H^1(\Omega_{\epsilon})$, we have 
$$\frac{d}{dt}W_{\epsilon}(v(t))=\frac{d}{dt}W_{\epsilon}(v_0)\,  {\Longrightarrow} -||\,v_t\,||_{L^2(\Omega_{\epsilon})}^{\,2}=0\,,$$ so $v$ must be an equilibrium of
 (\ref{abstract_perturbed}).

Therefore $W_{\epsilon}$ is a Lyapunov function for the flow generated by 
 (\ref{abstract_perturbed}), as claimed.
\eproof

The Lyapunov functionals $W_{\epsilon}$ in the  perturbed regions $\Omega_\epsilon$
 approach the functional $W = W_0$, in $\Omega_0:= \Omega$ when $\epsilon \to 0$ in the following sense: 

\begin{lema}\label{Wrel} If $W_\epsilon: H^1(\Omega_\epsilon) \to \R$ is as in Lemma \ref{funly}, then
  $$K_1(h {_\epsilon})|\,W(u)\,| \leq |\,W_\epsilon(u \circ h_\epsilon^{-1})\,|
  \leq K_2(h {_\epsilon})|\,W(u)\,|\,\,\, \forall \,u \in H^1(\Omega),$$
      with  $K_1(h {_\epsilon}), K_2(h {_\epsilon}) \to 1$  when
       $\epsilon \to 0$.
\end{lema}
\proof The result was proved in \cite{OPP} in the case where
$h_{\epsilon} \to i_{\Omega}$  in $\mathcal{C}^2$.
  However, since the inequalities involve only first order derivatives
   the extension 
$\mathcal{C}^1$ is immediate. \eproof

 We now define a  functional  ${V}_{\epsilon} : H^1(\Omega) \to \R $ by 
\begin{equation}\label{lyapfunctilde}
V_{\epsilon}(u) = W_{\epsilon}(u \circ h_{\epsilon}^{-1})\,.
\end{equation}

 \begin{lema}\label{lyapfunctildelemma}
 The functional   $V_{\epsilon} $ defined by (\ref{lyapfunctilde}) is
 a Lyapunov functional for (\ref{abstract_scale}), with $\eta = \frac{1}{2}$ and $-\frac{1}{2}<\beta < -\frac{1}{4}$, and the following estimates hold:
 \begin{enumerate}
 \item $\displaystyle K_1 
    \|u\|_{H^1(\Omega)} - K_2  \leq  V_{\epsilon}(u) \leq K_1
    \|u\|_{H^1(\Omega)} + K_2 $, for some constants $K_1$ and $K_2$,
   \item  $\displaystyle K_1(h{ {_\epsilon}})|\,V(u)\,| \leq |\, V_\epsilon(u )\,|
  \leq K_2(h{ {_\epsilon}})|\,V(u)\,|\,\,\, \forall \,u \in H^1(\Omega),$
      with  $K_1(h{ {_\epsilon}}), K_2(h{ {_\epsilon}}) \to 1$  when
      $\epsilon \to 0$.
      \end{enumerate}
\end{lema}
\proof
The required properties for  $V_{\epsilon} $ follow easily from the properties
of $W_{\epsilon}$ and the fact that  $ {h_{\epsilon}^{*}}^{-1} : H^{1}(\Omega) \to  H^{1}(\Omega_{\epsilon})$ is an  isomorphism and takes solutions of   (\ref{abstract_scale}) into solutions
of  (\ref{abstract_perturbed}).
\eproof

 Using the properties of the Lyapunov functional $V_\epsilon$, we now prove the
 following result of global existence for (\ref{abstract_scale}):

 \begin{teo}\label{global_exist}
   Suppose that $\beta$, $\eta$, $f$ and $g$ satisfy the conditions of
   Theorem \ref{loc_exist} and, additionally, that $f$ and $g$ satisfy the
   dissipative conditions (\ref{dissipative}) and (\ref{compet}). Then
    if $\epsilon> 0$ is sufficiently small the solutions of
    (\ref{abstract_scale}) are globally defined.
    \end{teo}
 \proof
 We first consider the case $\eta=\frac{1}{2}$. From Theorem \ref{loc_exist},
 for each $(t_0, u_0) \in \R \times H^{1}(\Omega)$, there exists
 $T=T(t_0,u_0)>0$ such that the problem (\ref{abstract_scale})
 has a unique solution  $u$ in $(t_0, t_0 + T)$, with  $u(t_0)=u_0$.
 From  {Lemmas} \ref{Flip} and  \ref{Glip} it follows  that $ {(}H {_\epsilon)_\beta}$ is locally Lipschitz, so it
 takes bounded sets  in $X^\eta  {=} H^1(\Omega  {)}$, into bounded sets
 of  in $X^{\beta}$.  Suppose that  $T < \infty$.  Then, by Theorem \ref{334}
 there exists a sequence  $t_n \to T^{-}$ such that $||\,u(t_n)\,||_{H^1(\Omega)} \to \infty$, and thus  $|V_{\epsilon}(u(t_n))| \to \infty$, which is a contradiction with the fact that  $V_{\epsilon}$ is decreasing along orbits.

 Now, suppose that $\eta= \eta_0 <\frac{1}{2}$, and let $u(t,u_0)$, be the solution
 with initial value $u(t_0)=u_0 \in X^{\eta_0}$, defined for $t_0 < t < T$.
 Then
 $u(t, u_0) \in X^{1+ \beta} \subset H^{1}(\Omega)$, for $t_0 < t < T$. Let
 $t_1 \in (t_0,T)$, $u_1 = u(t_1, u_0)$. Then the solution $v(t, u_1)$,
 $ v(t_1)=u_1$ of (\ref{abstract_scale}), with $\eta=\frac{1}{ {2}}$ is defined for $t \in (t_1, \infty)$.
 Since  $v(t, u_1)$ is also a solution of  (\ref{abstract_scale})
 with $\eta = \eta_0$, it must coincide with $u(t, u_0)$, for
 $t \in (t_1, T)$.  Define $ {\widetilde{u}}$ in $0 < t < \infty$, by
 $\widetilde{u}(t)= u(t, u_0)$, if  $t_0 < t < T$ and
 $\widetilde{u}(t)= v(t, u_0)$, if  $t_1< t < \infty$. Then
 $ {\widetilde {u}}$ is a solution of  (\ref{abstract_scale})
 with $\eta = \eta_0$ in  $(t_0, \infty)$.

  Finally, suppose that $\eta= \eta_0 >\frac{1}{2}$, and let $u(t,u_0)$, be the solution
  with initial value $u(t_0)=u_0 \in X^{\eta_0}$, defined for $t_0 < t < T$ and
$v(t, u_0)$,
  the solution  of (\ref{abstract_scale}), with $\eta=1$ and $v(t_0) = u_0$  defined for
  $t \in (t_0, \infty)$. Since $u(t,u_0)$ and $v(t,u_0)$ are both solutions of
  (\ref{abstract_scale}), with $\eta=1$, they must coincide in $(t_0, T)$.
  Since $v(t,u_0) \in X^{\eta_0}$, for $t>0$, it is a global solution
   of (\ref{abstract_scale})
 with $\eta = \eta_0$ in  $(t_0, \infty)$. \eproof

 From now on, we denote the flow generated by (\ref{abstract_scale})
 by $T_{\eta,\beta}(t,t_0,u)$  or $T_{\eta,\beta}(t,u)$ and also sometimes do not include the parameters
  $\eta$ and $\beta$ to simplify the notation.

\section{Existence of  global attractors} \label{existglob}

In this section, we prove that the flow $T_{\epsilon, \eta,\beta}(t,t_0,u)$  admits a global attractor. As in the previous section, it is convenient to start with
the especial case $\eta=\frac{1}{2}$.
\begin{lema}\label{gradsys}
  If  $\epsilon > 0$ is sufficiently small, the nonlinear semigroup
  $T_\epsilon(t)$ generated by  (\ref{abstract_scale}), with
   $\eta=\frac{1}{2}$ and $ -\frac{1}{2} < \beta < -\frac{1}{4}$ in
  $X^\frac{1}{2}=H^1(\Omega)$ is a gradient flow.
\end{lema}
\proof
We know, from Lemma \ref{lyapfunctildelemma} that the map $V_{\epsilon}$
defined  by  (\ref{lyapfunctilde}) is a Lyapunov functional for the flow.
Since $(A_{\epsilon})_\beta$ has compact resolvent and $ (H {_\epsilon)_\beta}$ takes bounded subsets of
$H^1(\Omega)$, into bounded subsets of $X^{\beta}$, it follows from Theorem
\ref{336} that bounded positive orbits are precompact, so $T_\epsilon(t)$
is gradient.
\eproof

We now want to show that the set of equilibria of (\ref{abstract_scale})
is bounded. As in the previous sections, it is convenient to prove this first for the problem in the perturbed domain.

\begin{lema}\label{equil_perturbed}
   Suppose that  $f$ and $g$ satisfy the conditions of
   Theorem \ref{loc_exist} and the
   dissipative conditions (\ref{dissipative}) and (\ref{compet}). Then
    if $\epsilon> 0$ is sufficiently small, the set
    $E_{\epsilon}$ of equilibria of the system  generated by
    (\ref{abstract_perturbed}) is uniformly bounded in
     $H^1(\Omega_{\epsilon})$.
\end{lema}
\proof The equilibria of (\ref{abstract_perturbed}) are the solutions of the
 problem:
$$
\begin{array}{rcl}
\left\{
\begin{array}{rcl}
\Delta u(x) -au(x) + f(u(x)) &=& 0 \,,\,\,\ x \in \Omega_{\epsilon} 
\\
\displaystyle\frac{\partial u}{\partial N}(x)&=&g(u(x))\,, \,\, x \in \partial\Omega_{\epsilon} \,.
\end{array}
\right.
\end{array}
$$
\par Multiplying by $u$ and integrating, we obtain

\begin{eqnarray} \label{eqest1}
  0& =& \displaystyle\int_{\Omega_{\epsilon}}u\Delta u \,dx - a\displaystyle\int_{\Omega_{\epsilon}} |\,u\,|^{\,2}dx + \displaystyle\int_{\Omega_{\epsilon}} uf(u)\,dx \nonumber \\
  &=&  - \displaystyle\int_{\Omega_{\epsilon}} |\,\nabla u\,|^{\,2}dx +\displaystyle\int_{\partial\Omega_{\epsilon}}u\frac{\partial u}{\partial N}\,dS - a\displaystyle\int_{\Omega_{\epsilon}} |\,u\,|^{\,2}dx + \displaystyle\int_{\Omega_{\epsilon}} uf(u)\,dx \nonumber \\
  &=&  - \displaystyle\int_{\Omega_{\epsilon}} |\,\nabla u\,|^{\,2}dx +
 \displaystyle\int_{\partial\Omega_{\epsilon}}ug(u)\,dS\,
  - a\displaystyle\int_{\Omega_{\epsilon}} |\,u\,|^{\,2}dx + \displaystyle\int_{\Omega_{\epsilon}} uf(u)\,dx\,. 
\end{eqnarray}

From (\ref{dissipative}), for any $\delta >0$ there exist   constants $K_\delta$ such that

\[  f(u)u \leq (c_0 + \delta)u^{\,2}  + K_\delta\  \textrm{and}  \
\ g(u)u \leq (d_0' + \delta)u^{\,2} + K_\delta  \  \textrm{for} \  u \in \R\,.
\]

From  (\ref{eqest1}), we then obtain
\begin{equation} \label{eqest2}
  \displaystyle\int_{\Omega_{\epsilon}} |\,\nabla u\,|^{\,2}dx
  \leq  -  a\displaystyle\int_{\Omega_{\epsilon}} |\,u\,|^{\,2}dx
 + 
 (c_0 +\delta) \displaystyle\int_{\Omega_{\epsilon}}  u^2\,dx 
  +(d_0'+\delta)\displaystyle\int_{\partial\Omega_{\epsilon}}u^2\,dS\,
  +K_\delta ( |\Omega| + |\partial \Omega|).
\end{equation}

On the other hand, since the first eigenvalue 
$\lambda_0(\epsilon)$ of the problem (\ref{compet}) is positive,
it follows that 
$ \displaystyle{0< \lambda_0(\epsilon) = \inf_{u\, \in \,H^1(\Omega_{\epsilon})}\frac{\left\langle -\Delta u + (a - c_0)u\,,\,u\right\rangle }{||\,u\,||_{L^2(\Omega_{\epsilon})}}}\,.$
 Thus, for all $u \in H^1(\Omega_\epsilon)$, we have
 \begin{eqnarray*}
   \lambda_0(\epsilon)\displaystyle\int_{\Omega_{\epsilon}} |\,u\,|^{\,2}dx\,
& \leq & \displaystyle\int_{\Omega_{\epsilon}} |\,\nabla u\,|^{\,2}\,dx - d_0\displaystyle\int_{\partial\Omega_{\epsilon}} |\,u\,|^{\,2}\,dS\, + (a - c_0)\displaystyle\int_{\Omega_{\epsilon}} |\,u\,|^{\,2}dx\,.
 \end{eqnarray*}

 From  {(\ref{eqest2})}, we then obtain
\begin{eqnarray} 
 \lambda_0(\epsilon)\displaystyle\int_{\Omega_{\epsilon}} |u|^2dx\, &\leq &
 -  a\displaystyle\int_{\Omega_{\epsilon}} |u|^2dx
 + 
 (c_0 +\delta) \displaystyle\int_{\Omega_{\epsilon}}  u^2\,dx 
  +(d_0'+\delta)\displaystyle\int_{\partial\Omega_{\epsilon}}u^2\,dS\,
  +K_\delta ( |\Omega| + |\partial \Omega|) \nonumber  \\ 
  & &  - d_0\displaystyle\int_{\partial\Omega_{\epsilon}} |\,u\,|^{\,2}\,dS\, + (a - c_0)\displaystyle\int_{\Omega_{\epsilon}} |\,u\,|^{\,2}dx \nonumber \\
  & \leq &
 \delta \displaystyle\int_{\Omega_{\epsilon}}  u^2\,dx 
  +(d_0'- d_0+\delta)\displaystyle\int_{\partial\Omega_{\epsilon}}u^2\,dS\,
  +K_\delta ( |\Omega| + |\partial \Omega|)  .
     \end{eqnarray}
Choosing $\delta <\min \{\lambda_0(\epsilon), (d_o-d_0')   \}$, we obtain
\begin{equation} \label{eqest3}
\displaystyle\int_{\Omega_{\epsilon}}  u^2\,dx 
  + \displaystyle\int_{\partial\Omega_{\epsilon}}u^2\,dS\, \leq
   \frac{K_\delta ( |\Omega| + |\partial \Omega|)}{l_\delta}, 
  \end{equation}
where $ l_{\delta}= \min \{\lambda_0(\epsilon)- \delta, (d_o-d_0')-\delta   \}$.

Finally, using (\ref{eqest2}) once again, we obtain

\begin{equation} \label{eqest4}
  \displaystyle\int_{\Omega_{\epsilon}} |\,\nabla u\,|^{\,2}dx
  \leq  
 (c_0 +d_0'+2\delta)  \frac{K_\delta ( |\Omega| + |\partial \Omega|)}{l_\delta}
  +K_\delta ( |\Omega| + |\partial \Omega|).
\end{equation}

The claim now follows immediately from  (\ref{eqest3}) and
 (\ref{eqest4}), observing that the constants in those estimates can be chosen uniformly in $\epsilon$ for $\epsilon$ close to $0$.
\eproof

\begin{cor} \label{equil_H1}
  Suppose that  $f$ and $g$ satisfy the conditions of Lemma
  \ref{equil_perturbed}. Then
    if $\epsilon> 0$ is sufficiently small, the set
    $E_{_\epsilon}$ of equilibria of the system  generated by
    (\ref{abstract_scale}) is uniformly bounded in
     $H^1(\Omega)$.
  \end{cor}
\proof
Since $u$ is an equilibria of (\ref{abstract_scale}) if and only if
$u= h {_\epsilon}^{* {-1}}v$, the result follows from Lemmas \ref{isosobolev} and \ref{equil_perturbed}. 
\eproof

We are now in a position to prove the existence of attractors in the particular
 case $ X^{\eta} = X^{\frac{1}{2}}= H^{1}(\Omega)$.
\begin{teo}\label{global_attract_H1}
Suppose that  $f$ and $g$ satisfy the conditions of
   Theorem \ref{loc_exist} with $\eta= \frac{1}{2}$, and the
   dissipative conditions (\ref{dissipative}) and (\ref{compet}). Then
   if $\epsilon > 0$ is sufficiently small, the flow $T_\epsilon(t,u)$
   generated by (\ref{abstract_scale}) has a global attractor $\mathcal{A}_{\epsilon}$ in
     $X^{\frac{1}{2}}= H^{1}(\Omega)$ for $\epsilon_0$ sufficiently small.
  \end{teo}
\proof
We apply Theorem \ref{385}. We know, from Lemma   {\ref{gradsys}}, that
the flow is gradient and from  Corollary \ref{equil_H1}, its set of equilibria is bounded, so it remains to be proved that it is asymptotically smooth. The proof follows from the regularizing properties of the flow and is by now very standard, so we omit it (see, for example   {Theorem \ref{336}}). 
\eproof

We now extend the previous result to other phase spaces.

\begin{teo}\label{global_attract}
Suppose that  $f$ and $g$ satisfy the conditions of
   Theorem \ref{loc_exist}, and the
   dissipative conditions (\ref{dissipative}) and (\ref{compet}). Then
   if $\epsilon > 0$ is sufficiently small, the flow  $T_{\epsilon, \eta, \beta}(u)$
   generated by
(\ref{abstract_scale}) has a global attractor $\mathcal{A}_\epsilon$, for
   in
   $X^{\eta}$, for $\epsilon_0$ sufficiently small.
   Furthermore $\mathcal{A}_\epsilon$
    does not depend
  on either  $\eta$ or $\beta$.
  \end{teo}
\proof Suppose first that $\eta = \eta_0< \frac{1}{2}$.

Let $B$ be a bounded set
in $X^{\eta}$. We may suppose that $B=B_R$, the ball of radius $R$ in $X^{\eta}$,
centered at the origin.  By Lemmas \ref{Fbem} and \ref{Gbem}, there is a constant $N$, such that $\|(H_{\epsilon})_{\beta}(u)\|_{ X^{\beta} } <N$, for any $u$, with
$ \|u\|_{X^{\eta}} < 2R $. Now, let  $  T(t)u_0=
 u(t;t_o,u_0)$ be the solution
of (\ref{abstract_scale}) with $\eta= \eta_0$ and initial condition  $u_0 \in B_R$.  We may suppose also that
Re $\sigma \big((A_{\epsilon})_{\beta}\big) > w > 0$. Then, while
$\|u(t;t_0, u_0)\|_{X^{\eta}}\leq 2R$, we have, by the variation of constants formula:
\begin{eqnarray}\label{bound_Xalpha}
  ||\,u(t;t_0,u_0)\,||_{X^{\eta}} &\leq&
    ||\,e^{-(A_{\epsilon})_{\beta}(t-t_0)}u_0\,||_{X^{\eta}} + \displaystyle\int_{t_0}^t ||\,e^{-(A_{\epsilon})_{\beta}(t-s)}(H_\epsilon)_{\beta} \,u(s)\,||_{X^{\eta}}\,ds \nonumber \\
&\leq&  {||\,(A_\epsilon)_\beta^{\eta - \beta} e^{-(A_{\epsilon})_{\beta}(t-t_0)}u_0\,||_{X^\beta}} 
+\, \displaystyle\int_{t_0}^t ||\,(A_{\epsilon})_{\beta}^{\,\eta-\beta}e^{-(A_{\epsilon})_{\beta}(t-s)}(H_\epsilon)_{\beta}\,u(s)\,||_{X^{\beta}}\,ds  \nonumber \\
&\leq & {\overline{C}_{\beta, \eta}(t - t_0)^{\beta - \eta}}e^{-w(t-t_0)}||u_0||_{X^\beta} + NC_{\beta, \eta}\displaystyle\int_{t_0}^t(t-s)^{-(\eta-\beta)}e^{-w(t-s)}ds.
\end{eqnarray}

Let $T = \sup\left\{t \,\geq\, t_0 \, | \, u(s;t_0,u_0) \in B_{2R}, \textrm{ for all }
s\leq t \right\}$, and $\delta $ such that
$$ {\overline{C}_{\beta , \eta}(t - t_0)^{\beta - \eta}}e^{-w(t-t_0)}||\,u_0\,||_{X^\beta} + NC_{\beta, \eta}\displaystyle\int_{t_0}^t (t-s)^{-(\eta-\beta)}e^{ {-w(t-s)}}ds < 2R.$$ 

From (\ref{bound_Xalpha}), it follows
that $T\geq \delta$,
so the solutions with initial conditions in the ball of radius $R$ in $X^{\eta}$
remain in
 the ball of radius $2R$ in $X^{\eta}$ for $0 \leq t \leq T$.
  Now, using the variation of constants formula  again, we obtain, for any  $0 \leq t \leq T$.

\begin{eqnarray}
  ||\,u(t;t_0,u_0)\,||_{X^{\frac{1}{2}}} &\leq&
    ||\,e^{-(A_{\epsilon})_{\beta}(t-t_0)}u_0\,||_{X^{\frac{1}{2}}} + \displaystyle\int_{t_0}^t ||\,e^{-(A_{\epsilon})_{\beta}(t-s)}(H_\epsilon)_{\beta} \,u(s)\,||_{X^{\frac{1}{2}}}\,ds \nonumber \\
&\leq& ||(A_{\epsilon})_{\beta}^{\,\frac{1}{2}-\eta} e^{-(A_{\epsilon})_{\beta}(t-t_0)}u_0||_{X^\eta} 
+\, \displaystyle\int_{t_0}^t 
 ||\,(A_{\epsilon})_{\beta}^{\,\frac{1}{2}-\beta}e^{-(A_{\epsilon})_{\beta}(t-s)}(H_\epsilon)_{\beta}\,u(s)\,||_{X^{\beta}}\,ds  \nonumber \\
&\leq &  {\overline{C}}_{\beta} (t-t_0)^{-(\frac{1}{2}-\eta)} e^{-w(t-t_0)}||u_0||_{X^\beta} + N {C_\beta}\displaystyle\int_{t_0}^t (t-s)^{-(\frac{1}{2}-\beta)}e^{-w(t-s)}ds \,. \nonumber
\end{eqnarray}

Therefore  $T(t) B_R$ is in a bounded set of 
 $X^{\frac{1}{2}} =  H^{1} {(\Omega)}$ if $ 0 < t \leq T$. 
By Theorem \ref{global_attract_H1} the global attractor
 $\mathcal{A} {_\epsilon}$ of (\ref{abstract_scale}) with $\eta= \frac{1}{2}$ atracts $T(t)B$  in the norm of $X^{\frac{1}{2}} =
 H^1{ {(}\Omega {)}}$. Thus $\mathcal{A} {_\epsilon}$ also attracts $B$ in the 
norm of $X^{\eta}$. Since $\mathcal{A} {_\epsilon}$ is invariant for the flow $T(t)$, it must be the attractor (\ref{abstract_scale}) with $\eta= \eta_0$.

 Suppose now that $\frac{1}{2} <\eta = \eta_0 <
 \beta+1$.

 If $B$ be a bounded set
in $X^{\eta}$, it is also a bounded set in $X^{\frac{1}{2}}$, so it
is attracted by the 
 the global attractor
  $\mathcal{A} {_\epsilon}$ of (\ref{abstract_scale}) with $\eta= \frac{1}{2}$ in the norm of $X^{\frac{1}{2}} =
 H^1{ {(}\Omega {)}}$,
under the flow $T_{\frac{1}{2}}(t)$ of
 (\ref{abstract_scale}) with $\eta= \frac{1}{2}$
 (which coincides with the flow of
 (\ref{abstract_scale}) with $\eta= \eta_0 $ in $X^{\eta}$).
 Now we prove that  $T_{\frac{1}{2}}(t)$ is continuous as a map from $X^{\frac{1}{2}}$ into $X^{\eta}$, for $t> 0$.
 If $u_1, u_2 \in X^{\frac{1}{2}}$ and 
  $ u_i(t;t_0,u_i) = T_{\frac{1}{2}}(t)u_i$,$i=1,2$, we obtain
 from the variation of constants formula

 \begin{eqnarray*}
 && ||\,u(t;t_0,u_1) -u(t;t_0,u_1) \,||_{X^{\eta}} \nonumber \\ 
& \leq & 
     ||\,e^{-(A_{\epsilon})_{\beta}(t-t_0)}(u_1 -u_2) \,  ||_{X^{\eta}} +
 \displaystyle\int_{t_0}^t ||\,e^{-(A_{\epsilon})_{\beta}(t-s)}((H_\epsilon)_{\beta} \,(u_1(s)) -( H_\epsilon)_{\beta} \,(u_2(s))) \,||_{X^{\eta}}\,ds \nonumber \\
&\leq& ||(A_{\epsilon})_{\beta}^{\eta- \frac{1}{2}}
 e^{-(A_{\epsilon})_{\beta}(t-t_0)}(u_1-u_2)||_{X^\frac{1}{2}} 
+\, \displaystyle\int_{t_0}^t 
 ||(A_{\epsilon})_{\beta}{ {^{\eta-\frac{1}{2}}}}e^{-(A_{\epsilon})_{\beta}(t-s)}((H_\epsilon)_{\beta}\,u_1(s) - (H_\epsilon)_{\beta}\,u_1(s) ) \,||_{X{ {^\frac{1}{2}}}}\,ds  \nonumber \\
&\leq & C_{\eta} (t-t_0){ {^{\frac{1}{2} - \eta}}} e^{-w(t-t_0)}
||\,u_1 -u_2\,||_{X^\frac{1}{2}} + 
L_{\beta, \eta}C_{\beta, \eta}
\displaystyle\int_{t_0}^t (t-s)^{ {\frac{1}{2} - \eta}}e^{-w(t-s)} 
\|u_1(s) - u_2(s) \|_{ {X^\frac{1}{2}}} 
ds \,,
\end{eqnarray*}
where $L_{\beta, \eta}$ is a  (local) Lipschitz constant of 
$(H_\epsilon)_{\beta}$, which exists by Lemmas \ref{Flip} and
 \ref{Glip}.
 The claimed continuity follows then from the Gronwall's inequality.
 Therefore, if $V_{\delta}$ is a small neighborhood of the attractor
 $\mathcal{A} {_\epsilon}$ in $X^{\frac{1}{2}}$ which contains $T_{\frac{1}{2}}(t)B= T_{\eta}(t)B $, then $T_{\frac{1}{2}}(1)U$ is a small neighborhood
of $\mathcal{A} {_\epsilon}$ in $X^{\eta}$ which contains 
$T_{\frac{1}{2}}(t+1) B= T_{\eta}(t+1)B $.
Since $\mathcal{A} {_\epsilon} = T_{\frac{1}{2}}(1) \mathcal{A} {_\epsilon} \in X^{\eta}$
is invariant, it must be the attractor of $T_{\eta}$.  
 \eproof

 \section{Continuity of the attractors} \label{contattractors}

 In this section we prove that
the family of attractors of  (\ref{abstract_scale}) is continuous in $X^{\eta}$ 
at $\epsilon = 0$, for $  \eta < \frac{1}{2}$.

\subsection{Uppersemicontinuity}

 We start by  proving  that the family of attractors  $\mathcal{A}_{\epsilon}$
is uniformly bounded  in $H^1(\Omega)$. More precisely,

\begin{teo} Suppose the hypotheses of Theorem \ref{global_attract} hold.
\label{bound_attract} Then, the  family of attractors
    $\left\{ \mathcal{A}_{\epsilon}, \ \epsilon < \epsilon_0
    \right\}$ 
    of (\ref{abstract_scale}) is uniformly  bounded in 
    $H^1(\Omega)$  for some $\epsilon_0> 0$ and
     $\mathcal{A}_{\epsilon}$ is bounded in $L^{\infty}$ for each $\epsilon$.
\end{teo}
  \proof
   Denote by $T_\epsilon(t)$, the flow generated by
   (\ref{abstract_scale}), and $T= T_0$. By Corollary \ref{equil_H1}, there is an $\epsilon_0 >0$, such
  that  the set $E_{\epsilon}$ of equilibria of
  (\ref{abstract_scale}) is in the open ball of radius $r$,  $B_r$ of $H^1(\Omega)$
  for $\epsilon < \epsilon_0$.   If $u \in \mathcal{A_\epsilon}$, by  Theorem
  \ref{385}, there is a $t_u$ such that  $u = T (t_u) u_0$ for some $u_0 \in B_r$. 
Let $V_\epsilon$ be the Lyapunov functional of (\ref{abstract_scale}) given by
 (\ref{lyapfunctilde}).
  We have
  $$\displaystyle V_{\epsilon}(u_0) \leq K_1
    \|u_0\|_{H^1(\Omega)} + K_2 \leq   {K_1}r + K_2. $$
  It follows that  $V_\epsilon(u_0 ) \leq R$, for some constant $R$ depending only on $r$. Thus
     
      \[
        V_\epsilon(u) \leq  V_\epsilon(T(t_u)u_0) \\
         \leq  V_\epsilon(u_0) \\
         \leq   R.
      \]

      From  {Lemma} \ref{lyapfunctildelemma}, it follows that
        $$\|u\|_{H^1(\Omega)} \leq \frac{1}{K_1}\left(   V_{\epsilon}(u) +K_2\right),
      $$
      for any $u \in \displaystyle\bigcup_{\epsilon\leq \epsilon_0} \mathcal{A_\epsilon}$,
      if $\epsilon_0$ is small enough.
      The $L^{\infty}$ estimate follows immediately from the fact that the attractors do not depend on $\eta$, and $X^{\eta}$ is continuously embedded in $L^{\infty}$, for $ \eta > \frac{1}{2}$, by Theorem \ref{c}.
   \eproof

   We are now  ready prove the upper semicontinuity property of the family.

   \begin{teo}\label{upperapp} Suppose the hypotheses of Theorem \ref{global_attract} hold with  $\beta = -\frac{1}{2}$.
   Then the family of attractors
     \{ $\mathcal{A}_{\epsilon } \, {,}  $ \ $0\leq \epsilon \leq \epsilon_0$\} 
     of  the flow  $T_{\epsilon,\eta}(t,u)$, generated by
     (\ref{abstract_scale}), whose existence is guaranteed by
     Theorem \ref{global_attract} is uppersemicontinuous in   $X^\eta$.
     (We observe that the conditions on  $\eta$ hold if
     $\frac{1}{2} -\delta \eta< \frac{1}{2}$, with $\delta$ sufficiently small).
\end{teo}
\proof
We have to check that the conditions of Theorem \ref{upper} hold.
In fact we have
\begin{enumerate}
\item We showed, during the proof of Theorem \ref{sectorial_weak}
  (see estimate (\ref{dif_estimate})), that the family of operators
  $A_\epsilon=(A_\epsilon)_{-\frac{1}{2}} $
  defined in $X^{-\frac{1}{2}}= H^{-1}(\Omega)$ with domain  $X^{\frac{1}{2}}= H^{1}(\Omega)$ satisfy the conditions of
  Theorem \ref{cont}.
\item From  {Lemmas} \ref{Fbem} and \ref{Gbem}, it follows that
  $H_\epsilon= F_{\epsilon}+ G_{\epsilon}$ takes bounded sets in $X^{\eta}$ into bounded sets of $X^{-\frac{1}{2}}= H^{-1}(\Omega)$.
\item From  {Lemmas} \ref{Flip} and \ref{Glip}, $H_\epsilon$ is continuous in $\epsilon$, uniformly for $u$ in bounded sets of $X^{\eta}$ and locally Lipschitz in $u$,
   uniformly in $\epsilon$.
\item By Theorem \ref{bound_attract}, the  set
  $\displaystyle\bigcup_{\epsilon\leq \epsilon_0} {\mathcal{A}_{\epsilon }}$  is bounded in  $X^{\frac{1}{2}}= H^{1}(\Omega)$ and thus by (2)  $H_\epsilon$ takes bounded sets in $X^{\frac{1}{2}}$ into bounded sets of $X^{-\frac{1}{2}}= H^{-1}(\Omega)$.
\end{enumerate}
Therefore, the upper semicontinuity follows from Theorem \ref{upper}.
\eproof

From the semicontinuity of attractors, we can easily prove the corresponding
 property for the equilibria.
\begin{cor}\label{upperequil}
 Suppose the hypotheses of Theorem \ref{global_attract} hold with  $\beta = -\frac{1}{2}$.
   Then the family of sets of equilibria
   $ {\{} E_{\epsilon } \,   {|} \, 0\leq \epsilon \leq \epsilon_0  {\}}$, of
     the problem  (\ref{abstract_scale}) is uppersemicontinuous in $X^\eta$.
\end{cor}
\proof
The result is well-known, but we sketch a proof here for completeness.
Suppose $u_{n}$, with $\displaystyle\lim_{n\to \infty}\epsilon_n= 0$. We choose an arbitrary subsequence and still call it   $(u_{n})$, for simplicity. It is enough to show that, there exists a subsequence   $(u_{n_k})$, which converges to a point $u_0 \in E_0$. Since $(u_{n}) \to \mathcal{A} {_\epsilon}$ and $\mathcal{A} {_\epsilon}$ is compact,  there
exists  a subsequence   $(u_{n_k})$, which converges to a point $u_0 \in
\mathcal{A} {_\epsilon}$. Now, since the flow $T_\epsilon(t)$ is continuous in $\epsilon$
we have, for
any $t>0$
\[ u_{n_k} \to u_0 \Leftrightarrow  T_{\epsilon_{n_k}}(t) u_{n_k} \to
  T_{0}(t) u_0 \Leftrightarrow   u_{n_k} \to
  T_{0}(t) u_0, \]
  so by uniqueness of the limit $T_{0}(t) u_0 =u_0$, for any $t> 0$, so
 $u {_0} \in E_0. $ 
\eproof

\subsection{Lowersemicontinuity}

For the lower semicontinuity  we will need to 
assume additional properties for the  nonlinearities.

\begin{equation} \label{boundfg}
 f \textrm{ and } g  \textrm{ are in }  C^1(\R,\R)  \textrm{ with bounded derivatives }.  
\end{equation}



\begin{lema}
 \label{FGateaux}
 If $f$ satisfies (\ref{boundfg}) and $\eta>0$  then
the operator 
 $F :X^{\eta}\times  \R  {\rightarrow} X^{\beta}$ given by 
 (\ref{Fh}) is Gateaux differentiable with respect to
 $u$, with Gateaux differential
 $\displaystyle{\frac{\partial F}{\partial u}(u,\epsilon)w}$ given by 
 \begin{equation}\label{FGateaux_form}
\left\langle \frac{\partial F}{\partial u}(u,\epsilon)w\,,\,
\Phi\right\rangle_{\beta,-\beta} = \displaystyle\int_\Omega f^{\,'}(u)w\,\Phi\,dx\,,
\end{equation}
for all  $w \in X^\eta$ and $\Phi \in X^{-\beta}$.
\end{lema}
\proof
 Observe first that $F(u, {\epsilon})$ is well-defined, since the conditions of Lemma \ref{Fbem} are met, with $\lambda_1=0$.

It is clear that $\displaystyle\frac{\partial F}{\partial u}(u,\epsilon)$ is linear. We now show that it is bounded. In fact we have,
 for all 
 $u,w \in X^{\eta} $  and
 $\Phi \in X^{-\beta}$
\begin{eqnarray*}
\left|\left\langle \frac{\partial F}{\partial u}(u,\epsilon)w\,,\,
\Phi\right\rangle_{\beta,-\beta} \right|  &\leq & \displaystyle\int_\Omega 
|\,f^{\,'}(u)\,|\, |\,w\,| \,|\,\Phi\,|\, dx\ \\
&\leq &  \|\,f'\,\|_{\infty}\displaystyle\int_\Omega |\,w\,| \, |\,\Phi\,|\,dx \\
&\leq &  \|\,f'\,\|_{\infty}\displaystyle \|\,w\,\|_{L^2(\Omega)} \,
 \|\,\Phi\,\|_{L^2(\Omega)}\,dx \\
&\leq &  \|\,f'\,\|_{\infty}\displaystyle K_1 K_3\|w\|_{ X^{\eta}} \,
 \|\,\Phi\,\|_{ X^{-\beta}}\,dx\,, 
\end{eqnarray*}
where $\|f'\|_{\infty}= \sup \{f'(x)  \,|\, x \in \R \} $ and $K_1, K_3$ are 
embedding constants. This proves boundeness. 

Now, we  have,  for all 
 $u,w \in X^{\eta} $  and
 $\Phi \in X^{-\beta}$
\begin{eqnarray*}
&& \left|\frac{1}{t}\left\langle F(u + tw,\epsilon) - F(u,\epsilon) - t  \frac{\partial F}{\partial u}(u,\epsilon) w,\Phi\right\rangle_{\beta,-\beta}\right| \\
&\leq& \frac{1}{|t|}\displaystyle \int_\Omega \big|\,[\,f(u + tw) - f(u)
  - tf^{\,'}(u)w\,]\,\Phi\,\big|\,dx \\
&\leq & K_1 \frac{1}{|t|} \left(\displaystyle\int_\Omega \big|f(u + tw) - f(u) - tf^{\,'}(u) w\big|^{2}dx\right)^\frac{1}{2}||\Phi||_{X^{-\beta}} \\
& \leq & 
 K_1 \left(\,\displaystyle \underbrace{\int_\Omega \big|\,
 \left(f'(u + \bar{t}w ) -  f^{\,'}(u)\right) w\,\big|^{\,2} \,dx\,
}_{(I)}\right)^\frac{1}{2}||\,\Phi\,||_{X^{-\beta}}, 
\end{eqnarray*}
where $K_1$ is the embedding constant of $X^{-\beta}$ in
 $L^2(\Omega)$, and $  0 \leq \bar{t} \leq t$.  
Since $f'$ is bounded, the integrand of $(I)$
  is bounded by an integrable function and goes to $0$ as $t \to 0$.
 Thus, the integral $(I)$ goes to $ 0$  as $t \to 0$, from Lebesgue's Dominated
Convergence Theorem. It follows that 
$ \displaystyle{\lim_{t \to 0} \frac{ F(u + tw {,\epsilon}) - F(u {,\epsilon})}{t}  = \frac{\partial F}{\partial u}(u,\epsilon) w \ \textrm{ in}  \  X^{\beta},}$
 for all  $u,w \in X^{\eta} $; so $F$ is Gateaux differentiable with 
 Gateaux differential given by (\ref{FGateaux_form}).
\eproof

We now want to prove that the Gateaux differential of $F(u, \epsilon)$
is continuous in $u$. Let us denote by
 $\mathcal{B}(X, Y)$  
the space of
  linear bounded operators from $X$ to $Y$.
 We will need the following result, whose simple proof is omitted.

 \begin{lema}\label{strong_uniform_operators} 
Suppose  $X,Y$ are Banach spaces and $ T_n : X \to Y$ is a sequence
 of linear operators converging strongly to the linear operator
 $T:X \to Y$. Suppose also that 
 $X_1 \subset X$  is a Banach space, the inclusion 
 $i: X_1 \hookrightarrow X$
 is compact and let $ \widetilde{T}_n = T_n \circ i$ 
 and  $ \widetilde{T} = T \circ i$.   Then
  $\widetilde{T}_n \to  \widetilde{T} $ uniformly for $x$ in a 
 bounded subset of $X_1$ (that is, in the or norm  of 
 $\mathcal{B}(X_1,Y )$).
\end{lema}

\begin{lema}\label{FGateaux_cont}
If $f$ satisfies (\ref{boundfg})  and $\eta>0$  then
the Gateaux differential of $F(u,\epsilon)$, with respect to $u$ is 
 continuous in $u$, that is, the map 
  $ u \mapsto \displaystyle\frac{\partial F}{\partial u}( u,\epsilon)
 \in \mathcal{B}(X {^\eta}, X^{\beta})$
 is continuous. 
\end{lema}
\proof Let  $u_n$ be a sequence converging to
 $u$ em $X^{\eta}$, and choose $0<\widetilde{\eta}< \eta$.
 Then, we have for any
 $\Phi \in X^{-\beta}$ and  $w \in X^{\widetilde{\eta}}$:
\begin{eqnarray*} 
\bigg|\left\langle \left(  \frac{\partial F}{\partial u}(u_n, {\epsilon}) -  \frac{\partial F}{\partial u}(u, {\epsilon})\right)w\,,\,\Phi \right\rangle_{\beta\,,\, -\beta}\bigg|
&\leq& \displaystyle\int_\Omega \bigg|\,\big(\,f^{\,'}(u) - f^{\,'}(u_n)\,\big)w\,\Phi\,\bigg|\,dx \nonumber \\
&\leq& \bigg(\displaystyle\int_\Omega \big|\big(f^{\,'}(u) - f^{\,'}(u_n)\big)w\big|^{2}dx\bigg)^\frac{1}{2}\bigg(\displaystyle\int_\Omega |\Phi\big|^{2}dx\bigg)^\frac{1}{2} \nonumber \\
&\leq& K_1\,\bigg(\displaystyle\underbrace{\int_\Omega \big|\big(f^{\,'}(u) - f^{\,'}(u_n)\big)w\big|^{2}dx}_{(I)}\bigg)^\frac{1}{2}\,||\Phi||_{ {X^{-\beta}}}\,, 
\end{eqnarray*}
where $K_1$ is the embedding constant of    $X^{-\beta}$ in
 $L^2(\Omega)$.

\par Now, the integrand in $(I)$ is bounded by the
 integrable function $\,||\,f^{\,'}||_\infty^{\,2}\,w^{\,2}$ and
 goes to $0$ a.e. as   $u_n \to u$ in $X^\eta$.
 Therefore the sequence of operators
$\displaystyle\frac{\partial F}{\partial u}( u_n,\epsilon)$ converges strongly
in the space $\mathcal{B}(X^{\widetilde{\eta}}, X^{\beta})$  to the operator  $ \displaystyle\frac{\partial F}{\partial u}( u,\epsilon)$.
 From Lemma \ref{strong_uniform_operators} 
 the convergence holds in the norm of  $\mathcal{B}(X^{\eta}, X^{\beta})$,
 since  $X^{\eta}$ is compactly embedded in $X^{\widetilde{\eta}}$, by Theorem \ref{prop_xalpha}.
\eproof

\begin{lema}
 \label{GGateaux}
 If $g$ satisfies  (\ref{boundfg}) and  $\eta>\frac{1}{4}$  then
the operator 
 $G :X^{\eta}\times  \R  {\rightarrow} X^{\beta}$ given by 
 (\ref{Gh}) is Gateaux differentiable with respect to
$u$, with Gateaux differential
\begin{equation}\label{GGateaux_form}
\left\langle \frac{\partial G}{\partial u}(u, {\epsilon})w\,,\,\Phi\right\rangle_{ {\beta, - \beta}} = \displaystyle\int_{\partial\Omega} g^{\,'}(\gamma(u))\gamma(w)\,\gamma(\Phi)\,\left|\displaystyle\frac{J_{\partial\Omega}h_\epsilon}{Jh_\epsilon}\right|\,d\sigma(x)\,,
\end{equation}
for all  $w \in X^\eta$ and $\Phi \in X^{-\beta}$.
\end{lema}
\proof
Observe first that $G(u, {\epsilon})$ is well-defined, since the conditions of Lemma
\ref{Gbem} are met, with $\lambda_2=0$.

It is clear that $\displaystyle\frac{\partial G}{\partial u}(u,\epsilon)$ is linear. We now show that it is bounded. In fact we have,
 for all 
 $u,w \in X^{\eta} $  and
 $\Phi \in X^{-\beta}$
 \begin{eqnarray*}
   \left| \left\langle \frac{\partial G}{\partial u}(u, {\epsilon})w\,,\,\Phi\right\rangle_{\beta, -\beta} \right| &  = &  \left| \displaystyle\int_{\partial\Omega} g^{\,'}(\gamma(u))\gamma(w)\,\gamma(\Phi)\,\left|\displaystyle\frac{J_{\partial\Omega}h_\epsilon}{Jh_\epsilon}\right|\,d\sigma(x)\, \right| \\
   &\leq &
   \|\theta \|_{\infty} \, \|g'\|_{\infty}
   \displaystyle\int_{\partial\Omega} | \gamma(w)|\, |\gamma(\Phi) |\,
   \,d\sigma(x)\,
   \\
   &\leq &  \|\theta \|_{\infty} \, \|g'\|_{\infty}
   \displaystyle \|\gamma(w)\|_{L^2(\partial \Omega)} \,
 \|\,\gamma(\Phi)\,\|_{L^2(\partial \Omega)}\\
 &\leq &
  K_1 K_2  \|\theta \|_{\infty} \, \|g'\|_{\infty}
   \displaystyle \|w\|_{X^{\eta}}} \,
 \|\,\Phi\,\|_{X^{-\beta} \,,
\end{eqnarray*}
 where $\|g'\|_{\infty}= \sup \left\{g'(x)  \,|\, x \in \R \right\} $,
$\|\theta \|_{\infty} =   \sup \left\{ |\theta(x, \epsilon)| \, | \, 
 x\in \partial \Omega \right\} =  \sup \left\{
 \displaystyle\frac{J_{\partial\Omega}h_\epsilon}{Jh_\epsilon} (x) \, \bigg| \, 
 x\in \partial \Omega \right\}$
 and $K_1$, $K_2$ are embedding constants given by Theorems \ref{imbed_xalpha} and
  \ref{trace_hk}. This proves boundeness. 

Now, we  have,  for all 
 $u,w \in X^{\eta} $  and
 $\Phi \in X^{-\beta}$
\begin{eqnarray*}
 & & \left| \frac{1}{t} \left\langle G(u + tw {,\epsilon}) - G(u {,\epsilon}) - t  \frac{\partial G}{\partial u}(u,\epsilon) w\,,\,\Phi\right\rangle_{\beta,-\beta}\right|  \\ 
 &\leq&  \frac{1}{|t|}\displaystyle \int_{\partial  {\Omega}} \left|\,\left[\,g(\gamma(u + tw)) -g(\gamma(u))
 - tg'(\gamma(u))\right]\gamma(w)\,\right|\,  {\left|\gamma(\Phi)\right|}\,\left| \frac{J_{\partial\Omega}h_\epsilon}{Jh_\epsilon} \right|\, {d\sigma(x)}   \\
  & \leq& K_1  \|\theta \|_{\infty} 
 \frac{1}{|t|}\displaystyle \left\{ \int_{\partial  {\Omega}} \left|\,\left[\,g(\gamma(u + tw)) -
 g(\gamma(u))
  - tg'(\gamma(u))\,\right]\gamma(w)\right|\, {d\sigma(x)} \right\}^{\frac{1}{2}} 
  \| \Phi \|_{X^{-\beta}}    \\  
 &\leq& K_1 \|\theta \|_{\infty}  \displaystyle \left\{ \underbrace{\int_{\partial  {\Omega}}
 \left|\,\left[\,g'(\gamma(u + \bar{t}w)) 
  - g'(\gamma(u))\,\right]\gamma(w)\,\right|\, {d\sigma(x)}}_{(I)} \right\}^{\frac{1}{2}} 
  \| \Phi \|_{X^{-\beta}} \,,  
\end{eqnarray*}
where  $K_1$  is
the  embedding constant given by Theorems \ref{imbed_xalpha} and
  \ref{trace_hk},  and $  0 \leq \bar{t} \leq t$.  
Since $g'$ is bounded, the integrand of $(I)$
  is bounded by an integrable function and goes to $0$ as $t \to 0$.
 Thus, the integral $(I)$ goes to $ 0$  as $t \to 0$, from Lebesgue's Dominated
Convergence Theorem. It follows that 
$ \displaystyle{\lim_{t \to 0} \frac{ G(u + tw {,\epsilon}) - G(u {,\epsilon})}{t}  =
 \frac{\partial G}{\partial u}(u,\epsilon) w \ \textrm{ in}  \  X^{\beta},}$
 for all  $u,w \in X^{\eta} $; so $G$  is Gateaux differentiable with 
 Gateaux differential given by (\ref{GGateaux_form}).
\eproof

\begin{lema}\label{GGateaux_cont}
If $g$ satisfies (\ref{boundfg}) and $\eta>\frac{1}{4}$  then
the Gateaux differential of $G(u,\epsilon)$, with respect to $u$ is 
 continuous in $u$ (that is, the map 
  $ u \mapsto \displaystyle\frac{\partial G}{\partial u}( u,\epsilon)
 \in \mathcal{B}(X^{\eta}, X^{\beta})$
 is continuous) and uniformly continuous in $\epsilon$ for $u$ in bounded sets
 of $X^{\eta}$ and
 $0\leq \epsilon \leq \epsilon_0 <1$. 
\end{lema}

\proof Let $0\leq \epsilon \leq \epsilon_0$,  $u_n$ be a sequence converging to
 $u$ em $X^{\eta}$, and choose $\frac{1}{4} < \widetilde{\eta} < \eta$.
 Then, we have for any
 $\Phi \in X {^{-\beta}}$ and  $w \in X^{\widetilde{\eta}}$:
\begin{eqnarray*}
& &\bigg|\,\left\langle \,\left(\,  \frac{\partial G}{\partial u}(u_n,\epsilon) -  \frac{\partial G}{\partial u}(u,\epsilon)\,\right)w\,,\,\Phi \,\right\rangle_{\beta\,,\, -\beta}\,\bigg| \\
& \leq & \displaystyle\int_{\partial\Omega}
\left| \left(g'(\gamma(u)) - g'(\gamma(u_n))\right)
\gamma(w)\,\gamma(\Phi)\right|\,\left|\displaystyle\frac{J_{\partial\Omega}h_\epsilon}{Jh_\epsilon}\right|\,d\sigma(x)
\nonumber \\
& \leq & \|\theta_{\epsilon} \|_{\infty} \left\{ \displaystyle\int_{\partial\Omega}
\left| (g'(\gamma(u)) - g'(\gamma(u_n))
\gamma(w)\,\,\right|^2\,d\sigma(x)\, \right\}^{\frac{1}{2}}
  \left\{ \displaystyle\int_{\partial\Omega}
\left| \,\gamma(\Phi)\,\right|^2\,d\sigma(x)\, \right\}^{\frac{1}{2}} 
\nonumber \\
 &\leq &K \|\theta_{\epsilon}\|_{\infty} \left\{ \displaystyle \underbrace{\int_{\partial\Omega}
\left| (g'(\gamma(u)) - g'(\gamma(u_n))
\gamma(w)\,\,\right|^2\,d\sigma(x)}_{(I)}\, \right\}^{\frac{1}{2}}
 \|\Phi\|_{ {X^{-\beta} }}\,,
\nonumber 
\end{eqnarray*}
where $K$ is the  constant due to continuity of the trace map from
$X^{-\beta}$ into
 $L^2(\partial \Omega)$.

\par Now, the integrand in $(I)$ is bounded by the
 integrable function $4 ||\,g'\,||_\infty^2 \left| 
\gamma(w)\right|^2$ and
goes to $0$ a.e. as   $u_n \to u$  {in} $X^\eta$,
 by Theorems \ref{imbed_xalpha} and \ref{trace_hk}.
 Therefore the sequence of operators
 $ \displaystyle\frac{\partial G}{\partial u}( u_n,\epsilon)$ converges strongly
in the space $\mathcal{B}(X^{\widetilde{\eta}}, X^{\beta})$  to the operator  $ \displaystyle\frac{\partial G}{\partial u}( u,\epsilon)$.
 From Lemma \ref{strong_uniform_operators} 
 the convergence holds in the norm of  $\mathcal{B}(X^{\eta}, X^{\beta})$,
 since  $X^{\eta}$ is compactly embedded in $X^{\widetilde{\eta}}$, by Theorem \ref{prop_xalpha}.

 On the other hand, if $0\leq \epsilon_1 \leq \epsilon_2 <\epsilon_0$, we have
 for any
 $\Phi \in  {X^{-\beta}}$ and  $w \in X^{\eta}$:
\begin{eqnarray*} 
& &\bigg|\,\left\langle \,\left(\,  \frac{\partial G}{\partial u}(u,\epsilon_1) -  \frac{\partial G}{\partial u}(u,\epsilon_2)\,\right)w\,,\,\Phi \,\right\rangle_{\beta\,,\, -\beta}\,\bigg| \\
&  \leq & \displaystyle\int_{\partial\Omega}
\left|\,g'(\gamma(u))
\gamma(w)\,\gamma(\Phi)\,\right|\,\left|\displaystyle \theta_{\epsilon_1} -
\theta_{\epsilon_2}\right|\,d\sigma(x)\,,
\nonumber \\
& \leq&  \|\theta_{\epsilon_1} -
\theta_{\epsilon_2}  \|_{\infty} \left\{ \displaystyle\int_{\partial\Omega}
\left| g'(\gamma(u))
\gamma(w)\,\,\right|^2\,d\sigma(x)\, \right\}^{\frac{1}{2}}
  \left\{ \displaystyle\int_{\partial\Omega}
\left| \,\gamma(\Phi)\,\right|^2\,d\sigma(x)\, \right\}^{\frac{1}{2}} 
\nonumber \\
 & \leq & K K'\|\,g'\,\|_{\infty}\|\|\,w\, \|_{ {X^\eta}}
 \|\Phi\|_{ {X^{-\beta}}}
  \|\theta_{\epsilon_1}-\theta_{\epsilon_2}\|_{\infty},
 \nonumber
\end{eqnarray*}
where $K'$ is the  constant due to continuity of the trace map from
$X^{\eta}$ into
 $L^2(\partial \Omega)$.
 This proves  uniform continuity in $\epsilon$.
\eproof

 \begin{lema}\label{Hfrechet}
   If  {$f$} and $g$   satisfy the condition (\ref{boundfg})
   and $\eta>\frac{1}{4}$, then 
  the map
 $ {(}H {_\epsilon)_\beta}=  {(}F {_\epsilon)_\beta} +  {(}G {_\epsilon)_\beta}  :X^{\eta}\times  \R \mapsto X^{\beta}$ given by 
 (\ref{defH}) is continuously Fr\'echet differentiable with respect to
  $u$ and the derivative $\displaystyle\frac{\partial G}{\partial u}$ is uniformly continuous with respect to $\epsilon$, for $u$ in bounded sets
  of $X^{\eta}$ and $0\leq \epsilon \leq\epsilon_0 < 1$.
\end{lema}
\proof
The proof follows from Lemmas \ref{FGateaux_cont}, \ref{GGateaux_cont}
and Proposition 2.8 in \cite{Rall}.
\eproof 

We now prove lower semicontinuity for the equilibria.

\begin{teo}\label{equicon}
   {If $f$} and $g$   satisfy the conditions of Theorem  \ref{global_attract}
  and
  also (\ref{boundfg}), the equilibria
  of  (\ref{abstract_scale}) with $\epsilon = 0$ are all hyperbolic  
  and $\frac{1}{4}<\eta< \frac{1}{2}$, then the  family of sets of equilibria
  $\{ E_{\epsilon} \, | \, 0 \leq \epsilon <\epsilon_0 \}$ of
  (\ref{abstract_scale}) is
  lower semicontinuous in $X^{\eta}$ at $\epsilon = 0$.
\end{teo}
\proof A point $e \in X^{\eta} $ is an equilibrium of (\ref{abstract_scale})
if and only if it is a root of the map
$$
\begin{array}{rlc}
Z: H^1(\Omega) \times  {\R}& \longrightarrow &X^{-\frac{1}{2}} \, \\
(u\,,\,\epsilon)& \longmapsto & (A_{ {\epsilon}})_{-\frac{1}{2}}(u) + (H_{\epsilon})_{-\frac{1}{2}}(u)\,,
\end{array}
$$

By Lemma \ref{Hfrechet} the map $ {(}H_\epsilon {)_{-\frac{1}{2}}}: X^{\eta} \to X^{-\frac{1}{2}}$ is continuously Fr\'echet differentiable with
respect to $u$ and by Lemmas \ref{Glip} and \ref{Flip} it is also continuous in $\epsilon$
if $\eta= \frac{1}2 - \delta$, with $\delta>0$ is sufficiently small.
Therefore, the same holds if $\eta = \frac{1}{2}$.

The map  $A_\epsilon= -h_\epsilon^{*} \Delta_{\Omega_\epsilon} h_\epsilon^{*} \,+\, aI$
is a bounded  linear operator from  $H^1(\Omega)$ to $X^{-\frac{1}{2}}$.
It is also
 continuous in $\epsilon$ since it is analytic as a function of
$h {_\epsilon} \in Diff {^1}(\Omega)$ and
 $ h_\epsilon$ is continuous in  $\epsilon$.
 
 Thus, the map $Z$ is continuously differentiable in $u$ and continuous in
  $\epsilon$. 
 The derivative of $\displaystyle\frac{\partial Z}{\partial u}(e, 0)$
 is an isomorphism by hypotheses.
 Therefore,  the Implicit Function Theorem apply, implying that the
 zeroes of $Z(\cdot, \epsilon)$ are given by a continuous function
 $ e(\epsilon)$. This proves the claim. \eproof

 For to prove the lower semi continuity, we also need to prove the continuity of local unstable manifolds at equilibria, more precisely 

\begin{teo}\label{manifcont}
   Suppose  $f$ and $g$   satisfy the conditions of Theorem \ref{equicon},
   $u_0$ is an equilibrium of   (\ref{abstract_scale}) with $\epsilon = 0$,
   and for each $\epsilon>0$ sufficiently small, let $u {_\epsilon}$
   be the unique equilibrium of  (\ref{abstract_scale}), whose existence
   is asserted by Corollary \ref{upperequil} and Theorem \ref{equicon}.
   Then, for $\epsilon$ and $\delta$ sufficiently small, there exists a
   local unstable manifold
$
W_{\rm loc}^u(u_{\epsilon}) 
$ of $u_{\epsilon}$, and if we denote
$ W_{\delta}^u(u_{\epsilon}) =\{ w \in  W_{\rm loc}^u(u_{\epsilon})  \ | \
\|w-u_{\epsilon} \|_{X^{\eta}} < \delta  \}, then$
\[
  \beta \Big(W_{\delta}^u(u_{\epsilon}),W_{\delta}^u( u_0) \Big) \quad \textrm{and} \quad
   \beta \Big(W_{\delta}^u(u_{0}),W_{\delta}^u( u_{\epsilon}) \Big)
\]
approach zero as $\epsilon \to 0$, where
  $\beta(O,Q)=\displaystyle\sup_{o \in O} \inf_{q \in Q}
\|q-o\|_{ {X^{\eta}}}$ for $O$, $Q\subset X^{\eta}$.
\end{teo}
 \proof
 Let $H_{\epsilon}(u)=H(u,\epsilon)$ be the map defined by (\ref{defH}), with $\beta = -\frac{1}{2}$ and $u_{\epsilon}$ a hyperbolic equilibrium of (\ref{abstract_scale}). Since $H(u,\epsilon)$ is differentiable by Lemma \ref{Hfrechet},
 it follows that  $H_{\epsilon}(u_{\epsilon}+w , \epsilon)=
 H_{\epsilon}(u_{\epsilon},\epsilon)
   + H_u(u_{\epsilon} , \epsilon)w + r(w,\epsilon)= A_{\epsilon}u_{\epsilon} + H_u(u_{\epsilon} , \epsilon)w + r(w , \epsilon)$,
 with $r(w,\epsilon)=o(\|w\|_{X^\eta})$, as $\|w\|_{X^\eta} \to 0$.
 The claimed result was proved in \cite{PP}, assuming the following properties of $H_{\epsilon}$:
 
 \begin{itemize}
\item[a)]  $||\,r(w,0)-r(w,\epsilon)\,||_{X^{\beta}} \leq C({\epsilon})$,
  with  $C({\epsilon}) \to 0 \textrm{ when } \epsilon \to 0$, uniformly for $w$ in a neighborhood of $0$ in $X^{\eta}$.
\item[b)] $||\,r(w_1,\epsilon)-r(w_2,\epsilon)\,||_{X^{\beta}} \leq k(\rho) ||\,w_1-w_2\,||_{X^{\eta}}$ , for $||\,w_1\,||_{X^{\eta}}\leq \rho$, $||\,w_2\,||_{X^{\eta}}\leq \rho$, with $k(\rho) \to 0$ when $\rho \to  0^+$  and $k(*)$ is non decreasing.
\end{itemize}

 Property a) follows from easily from the fact that both $H(u,\epsilon)$ and
 $H_u(u,\epsilon)$ are
 uniformly continuous in $\epsilon$ for $u$ in bounded sets
 of $X^{\eta}$, by Lemmas \ref{Glip},  {\ref{Flip}}  and  \ref{Hfrechet}.
 It remains to prove property b).
 
 If  $w_1,w_2 \in X^{\eta}$ and  $\epsilon \in [0, \epsilon_0]$,  with
  $0 < \epsilon_0 <1 $  small enough, we have 
\begin{eqnarray}
||\,r(w_1\,,\,\epsilon)-r(w_2\,,\,\epsilon)\,||_{X^{\beta}} &=& ||\,H(u_{\epsilon} + w_1\,,\,\epsilon) - H(u_{\epsilon}\,,\,\epsilon) - H_u(u_{\epsilon}\,,\,\epsilon)w_1 \nonumber\\ &&-\,
H(u_{\epsilon} + w_2\,,\,\epsilon) + H_{\epsilon}(u_{\epsilon}\,,\,\epsilon) + H_u(u_{\epsilon}\,,\,\epsilon)w_2\,||_{X^{\beta}} 
\nonumber\\ 
& \leq& ||\,F(u_{\epsilon} + w_1\,,\,\epsilon) - F(u_{\epsilon}\,,\,\epsilon) - F_u(u_{\epsilon}\,,\,\epsilon)w_1 \label{7}\\ 
&&-\, F(u_{\epsilon} + w_2\,,\,\epsilon) + F(u_{\epsilon}\,,\,\epsilon) + F_u(u_{\epsilon}\,,\,\epsilon)w_2\,||_{X^{\beta}}\nonumber \\
&&+\, ||\,G(u_{\epsilon} + w_1\,,\,\epsilon) -G(u_{\epsilon}\,,\,\epsilon) - G_u(u_{\epsilon}\,,\,\epsilon)w_1 \label{8}\\ 
&&-\, G(u_{\epsilon} + w_2\,,\,\epsilon) + G(u_{\epsilon}\,,\,\epsilon) + G_u(u_{\epsilon}\,,\,\epsilon)w_2\,||_{X^{\beta}} \,.\nonumber
\end{eqnarray}
 We first estimate (\ref{7}). Since $f'$ is  bounded by (\ref{boundfg}), we have

\begin{eqnarray*}
&& \bigg|\,\left\langle \,F(u_{\epsilon} + w_1\,,\,\epsilon) - F(u_{\epsilon}\,,\,\epsilon) - F_u(u_{\epsilon}\,,\,\epsilon)w_1 \right.\\
  &&\left. -F(u_{\epsilon} + w_2\,,\,\epsilon) + F(u_{\epsilon}\,,\,\epsilon) + F_u(u_{\epsilon}\,,\,\epsilon)w_2 \,,\,\Phi\,\right\rangle_{ {\beta, -\beta}}
  \,\bigg| \\
&\leq&
\displaystyle\int_{\Omega} \left|\,[\,f(u_{\epsilon}+w_1)-f(u_{\epsilon})- f'(u_{\epsilon})w_1
  -f(u_{\epsilon}+w_2)+ f(u_{\epsilon})+ f'(u_{\epsilon})w_2\,]\,\Phi\, \right|\,dx\,
\nonumber\\
& =&
\displaystyle\int_{\Omega} \left|\,[\,f'(u_{\epsilon}+ \xi_x)-f'(u_{\epsilon})\,](w_1(x)-w_2(x))\,\Phi \,\right|\,dx \nonumber\\
& \leq&
K_2  \displaystyle \left\{\int_{\Omega} \left|\,[\,f'(u_{\epsilon}+ \xi_x)-f'(u_{\epsilon})\,]^2
(w_1(x)-w_2(x))^2\,\right|\,dx\right\}^{\frac{1}{2}} \|\Phi \|_{H^1(\Omega)}\,,
\nonumber\\
\end{eqnarray*}
where   $K_2$ is the embedding constant 
of   $X^{\frac{1}{2}}=H^1(\Omega)$  {into} $L^2(\Omega)$, and
$ w_1(x) \leq \xi_x \leq w_2(x)$ or $ w_2(x) \leq \xi_x \leq w_1(x)$. 

Now choosing  $q = \displaystyle\frac{1}{1- 2 \eta}$, we have by Theorems \ref{imbed_wk}
and \ref{imbed_xalpha} that $X^{\eta} $ is continuously embedded in $L^{2q}(\Omega)$
Therefore, if $p$ is the conjugate exponent of $q$ we have, by
H\"{o}lder's inequality
$$
\begin{array}{ll}
&\left\{\displaystyle\int_{\Omega} [\,f'(u_{\epsilon}+ \xi_x)-f'(u_{\epsilon})\,]^{\,2}(w_1(x)-w_2(x))^{\,2}dx \right\}^\frac{1}{2} 
\\ \leq&
\left\{ \displaystyle\int_{\Omega} [\,f'(u_{\epsilon}+ \xi_x)-f'(u_{\epsilon})\,]^{\,2p}dx \right\}^\frac{1}{2p}\left\{ \displaystyle\int_{\Omega}(w_1(x)-w_2(x))^{\,2q}dx \right\}^\frac{1}{2q}
\\ \leq& 
\left\{ \displaystyle\int_{\Omega} [\,f'(u_{\epsilon}+ \xi_x)-f'(u_{\epsilon})\,]^{\,2p}dx \right\}^\frac{1}{2p}
||\,w_1-w_2\,||_{L^{2q}(\Omega)}\\
\leq& K_3\left\{ \displaystyle\int_{\Omega} [\,f'(u_{\epsilon}+ \xi_x)-f'(u_{\epsilon})\,]^{\,2p}dx \right\}^\frac{1}{2p}
||\,w_1-w_2\,||_{X^{\eta}}\,,
\end{array}
$$
where $K_3$ is the embedding constant of
the embedding $X^{\eta}$  {into} $L^{2q}(\Omega$).
Therefore, we have
$$
\begin{array}{lll}
&&||\,F(u_{\epsilon} + w_1,\epsilon) - F(u_{\epsilon},\epsilon) - F_u(u_{\epsilon},\epsilon)w_1 -
F(u_{\epsilon} + w_2,\epsilon) + F(u_{\epsilon},\epsilon) + F_u(u_{\epsilon},\epsilon)w_2||_{X^{\beta}} \\
&\leq& K_2K_3\left\{ \displaystyle\int_{\Omega} [\,f'(u_{\epsilon}+ \xi_x)-f'(u_{\epsilon})\,]^{\,2p}dx \right\}^\frac{1}{2p}
||\,w_1-w_2\,||_{X^{\eta}}\,.
\end{array}
$$
Now the integrand above is bounded by  $4^{\,p}||f'||_\infty^{\,2p}$
 and goes a.e. to $0$ as $\rho \to 0$, since   $||\,w_1\,||_{X^{\eta}}\leq \rho$,  $||\,w_2\,||_{X^{\eta}}\leq \rho$  and  $w_1(x)\leq \xi_x\leq w_2(x)$. 
Thus, the integral goes to $0$ by Lebesgue's bounded convergence Theorem.
 
\quad We now estimate (\ref{8}):
$$
\begin{array}{lll}
&\bigg|\,\left\langle G(u_{\epsilon} + w_1\,,\,\epsilon) - G(u_{\epsilon}\,,\,\epsilon) - G_u(u_{\epsilon}\,,\,\epsilon)w_1 \right.\\
& \left. -G(u_{\epsilon} + w_2\,,\,\epsilon) + G(u_{\epsilon}\,,\,\epsilon) + G_u(u_{\epsilon}\,,\,\epsilon)w_2 \,,\,\Phi\,\right\rangle_{ {\beta, -\beta}}\,\bigg| \\
\leq&
\displaystyle\int_{\partial\Omega} \bigg|\,[\,g(\gamma(u_{\epsilon}+w_1)) - g(\gamma(u_{\epsilon}))- g'(\gamma(u_{\epsilon}))w_1 \\
&-g(\gamma(u_{\epsilon}+w_2)) + g(\gamma(u_{\epsilon}))+ g'(\gamma(u_{\epsilon}))w_2\,]
\,\gamma(\Phi)\gamma\left(\,\left|\,\displaystyle\frac{J_{\partial\Omega} {h_\epsilon}}{J {h_\epsilon}}\,\right|\,\right)\,\bigg|\,d\sigma(x)\\
=&
\displaystyle\int_{\partial\Omega} \left|\,[\,g'(\gamma(u_{\epsilon}+ \xi_x)) - g'(\gamma(u_{\epsilon}))\,]\gamma(\,w_1(x)-w_2(x)\,)\,\gamma(\Phi)\gamma\left(\,\left|\,\displaystyle\frac{J_{\partial\Omega} {h_\epsilon}}{J {h_\epsilon}}\,\right|\,\right)\, \right|\,d\sigma(x) \\ 
\leq &\overline{K}_2\,\left\{\displaystyle\int_{\partial\Omega} [\,(g'(\gamma(u_{\epsilon}+ \xi_x)) - g'(\gamma(u_{\epsilon})))]^{\,2}[\gamma(w_1(x)-w_2(x))]^{2}\left[\gamma\left(\left|\displaystyle\frac{J_{\partial\Omega} {h_\epsilon}}{J {h_\epsilon}}\right|\right)
\right]^{}d\sigma(x)\right\}^\frac{1}{2}||\Phi|| {_{X^\frac{1}{2}}}\,,
\end{array}
$$
where $\overline{K}_2$ is the constant of the embedding
$X {^\frac{1}{2}}=H^1(\Omega)$  {into} $L^2(\partial\Omega)$ given
by Theorem \ref{trace_hk}, and  $ w_1(x) \leq \xi_x \leq w_2(x)$ or $ w_2(x) \leq \xi_x \leq w_1(x)$.
Now choosing  $\overline{q} = \displaystyle\frac{1}{1- 2 \eta}$, we have by Theorems \ref{trace_hk}
and \ref{imbed_xalpha} that $X^{\eta} $ is continuously embedded in
$L^{2\overline{q}}(\partial
\Omega)$. Therefore, if $\overline{p}$ is the conjugate exponent of $\overline{q}$ we have, by
H\"{o}lder's inequality

$$
\begin{array}{ll}
&\left\{\,\displaystyle\int_{\partial\Omega} [\,g'(\gamma(u_{\epsilon}+ \xi_x)) - g'(\gamma(u_{\epsilon}))]^{\,2}[\,\gamma(w_1(x)-w_2(x))\,]^{\,2} \left[\,\gamma\left(\,\left|\,\displaystyle\frac{J_{\partial\Omega} {h_\epsilon}}{J {h_\epsilon}}\,\right|\,\right)
\,\right]^{\,2}d\sigma(x)\, \right\}^\frac{1}{2}\\ 
\leq& \left\{\,\displaystyle\int_{\partial\Omega} [\,g'(\gamma(u_{\epsilon}+ \xi_x)) - g'(\gamma(u_{\epsilon}))]^{\,2\overline{p}}\left[\,\gamma\left(\,\left|\,\displaystyle\frac{J_{\partial\Omega} {h_\epsilon}}{J {h_\epsilon}}\,\right|\,\right)\,
\right]^{\,2\overline{p}}d\sigma(x)\,\right\}^\frac{1}{2\overline{p}} \\
&\cdot\,\left\{\,\displaystyle\int_{\partial\Omega}[\gamma(w_1(x)-w_2(x))]^{\,2\overline{q}}d\sigma(x) \,\right\}^\frac{1}{2\overline{q}} \\ 
\leq& \left\{\,\displaystyle\int_{\partial\Omega} [g'(\gamma(u_{\epsilon}+ \xi_x)) - g'(\gamma(u_{\epsilon}))]^{\,2\overline{p}}\left[\,\gamma\left(\,\left|\,\displaystyle\frac{J_{\partial\Omega} {h_\epsilon}}{J {h_\epsilon}}\,\right|\,\right)\,\right]^{\,2\overline{p}}d\sigma(x) \,\right\}^\frac{1}{2\overline{p}}||\,w_1-w_2\,||_{L^{2\overline{q}}(\partial\Omega)}\\ 
\leq &\overline{K}_3\,\left\{\,\displaystyle\int_{\partial\Omega} [\,g'(\gamma(u_{\epsilon}+ \xi_x)) - g'(\gamma(u_{\epsilon}))]^{\,2\overline{p}} \left[\,\gamma\left(\,\left|\,\displaystyle\frac{J_{\partial\Omega} {h_\epsilon}}{J {h_\epsilon}}\,\right|\,\right)\,\right]^{\,2\overline{p}}d\sigma(x) \,\right\}^\frac{1}{2\overline{p}}||\,w_1-w_2\,||_{X^{\eta}}\,,
\end{array}
$$
where $\overline{K}_3$ is the embedding constant of
${X^{\eta}}$  {into} $L^{2\overline{q}}(\partial\Omega$).
Therefore, we have  
$$
\begin{array}{ll}
&||\,G(u_{\epsilon} + w_1\,,\,\epsilon) - G(u_{\epsilon}\,,\,\epsilon) - G_u(u_{\epsilon}\,,\,\epsilon)w_1 -\,G(u_{\epsilon} + w_2\,,\,\epsilon) + G(u_{\epsilon}\,,\,\epsilon) + G_u(u_{\epsilon}\,,\,\epsilon)w_2\,||_{X^{\beta}} \\
\leq &\overline{K}_2\,\overline{K}_3\,\left\{\,\displaystyle\int_{\partial\Omega} [\,g'(\gamma(u_{\epsilon}+ \xi_x)) - g'(\gamma(u_{\epsilon}))]^{\,2\overline{p}}\left[\,\gamma\left(\,\left|\,\displaystyle\frac{J_{\partial\Omega} {h_\epsilon}}{J {h_\epsilon}}\,\right|\,\right)
\,\right]^{\,2\overline{p}}d\sigma(x) \,\right\}^\frac{1}{2\overline{p}}||\,w_1-w_2\,||_{X^{\eta}}\,.
\end{array}
$$

Now the integrand above is bounded by  $4^{\,\overline{p}}||\,g'\,||_\infty^{\,2\overline{p}}||\,\theta\,||_\infty^{\,2\overline{p}}$, where $\theta$ is the bounded function
given by (\ref{theta_1})  {and} (\ref{theta_2}), and goes to $0$ a.e. as $\rho \to 0$, since $||\,w_1\,||_{X^{\eta}}\leq \rho$,  $||\,w_2\,||_{X^{\eta}}\leq \rho$  and  $w_1(x)\leq \xi_x\leq w_2(x)$. 
 Thus, the integral goes to $0$ by Lebesgue's dominated convergence Theorem.
 \eproof
 
We are now in a position to prove the main result of this section

\begin{teo}
Assume the hypotheses of Theorem \ref{equicon} hold.  Then the family of attractors
     \{$\mathcal{A}_{\epsilon }\,{|}\, 0 \leq \epsilon \leq \epsilon_0\}$, of
     the problem (\ref{abstract_scale}), whose existence is guaranteed by
     Theorem \ref{global_attract}  is lower semicontinuous in   $X^\eta$.
     (We observe that the conditions on  $\eta$ hold if
     $\frac{1}{2} -\delta \eta< \frac{1}{2}$, with $\delta$ sufficiently small.)
  \end{teo}
\proof

The system generated by
(\ref{abstract_scale}) is gradient for any $\epsilon$ and its equilibria are all hyperbolic for $\epsilon $ in a neighborhood of $0$. Also, 
 the equilibria are continuous in $\epsilon$ by Theorem \ref{equicon}, the  linearization is  continuous
 in $\epsilon$ as shown during  the proof of  Theorem \ref{equicon} and  the local unstable manifolds of the equilibria are continuous
 in $\epsilon$, by Theorem \ref{manifcont}. The result follows then from
Theorem \ref{lower}.
\eproof

\end{document}